\documentclass[12pt, a4paper]{article}
\usepackage{fullpage,amsthm,amsmath}
\usepackage{harrymacros}
\usepackage{enumitem}

\usepackage[T1]{fontenc}
\usepackage{bigfoot} 
\usepackage[numbered,framed]{matlab-prettifier}

\usepackage{filecontents}

\lstMakeShortInline"

\newcommand{\wnorm}[1]{\abs{#1}}
\newcommand{\esigma}{\sigma_{ess}}
\newcommand{\erho}{\rho_{ess}}

\newcommand{\ip}{\langle \cdot, \cdot \rangle}

\DeclareMathOperator{\Id}{Id}

\DeclareMathOperator{\Leb}{Leb}

\title{
Fourier approximation of the statistical properties of Anosov maps on tori.
}
\author{Harry Crimmins\footnote{h.crimmins@unsw.edu.au} and Gary Froyland\footnote{g.froyland@unsw.edu.au} \\ \\ School of Mathematics and Statistics \\ University of New South Wales \\ Sydney NSW 2052, Australia}

\begin{document}

\maketitle

\begin{abstract}
  We study the stability of statistical properties of Anosov maps on tori by examining the stability of the spectrum of an analytically twisted Perron-Frobenius operator on the anisotropic Banach spaces of Gou{\"e}zel and Liverani \cite{gouezel2006banach}.
  By extending our previous work in \cite{cfstability}, we obtain the stability of various statistical properties (the variance of a CLT and the rate function of an LDP) of Anosov maps to general perturbations, including new classes of numerical approximations. In particular, we obtain new results on the stability of the rate function under deterministic perturbations.
  As a key application, we focus on perturbations arising from numerical schemes and develop two new Fourier-analytic method for efficiently computing approximations of the aforementioned statistical properties.
  This includes the first example of a rigorous scheme for approximating the peripheral spectral data of the Perron-Frobenius operator of an Anosov map without mollification.
  Using the two schemes we obtain the first rigorous estimates of the variance and rate function for Anosov maps.
\end{abstract}

\newpage

\tableofcontents

\newpage

\section{Introduction}

We consider the stability of various statistical objects associated with Anosov maps $T:\mathbb{T}^d\circlearrowleft$ on tori, with a view to building rigorous numerical estimates of these objects.
In particular, we consider the Sinai-Ruelle-Bowen (SRB) measure $\mu$ on $\mathbb{T}^d$, the limiting variance of scaled Birkhoff sums of a smooth observation function $g:\mathbb{T}^d\to\mathbb{R}$ guaranteed by the central limit theorem (CLT), and the rate function associated with large deviations of Birkhoff sums of $g$ from $\mathbb{E}_\mu(g)$ as guaranteed by a large deviation principle (LDP).
We use the naive Nagaev-Guivarc'h method (see e.g.\ \cite{gouezel2015limit}) to relate the spectral data of an analytically twisted Perron-Frobenius operator to statistical properties of the system.
The stability these properties then follows from the stability of the spectrum of the twisted Perron-Frobenius operator \cite{cfstability}.

The stability of the SRB measure to deterministic and stochastic perturbations has been well studied.
Early results of Kifer \cite{Kifer_1974} and Young \cite{young1986stochastic} established the stability of SRB measures under small stochastic perturbations.
Differentiability of the SRB measure under sufficiently smooth deterministic perturbations was proven for Anosov maps \cite{katok1989differentiability}, followed by a more detailed analysis in \cite{ruelle1997differentiation, ruelle2008differentiation}.
Spectral approaches for higher order differentiability of the SRB measure are developed in \cite{gouezel2006banach, butterley2007smooth}.

The rigorous numerical approximation of SRB measures of Anosov maps has been considered in \cite{F95} where Ulam's method is applied on Markov partitions, and in \cite{DJ99} via a two-step process, whereby one first convolves with a locally supported stochastic kernel, and then applies Ulam.
The paper \cite{BKL02} improves on the approach of \cite{DJ99} by linking the support radius of the kernel to the size of the Ulam partition elements, resulting in a single limit process;  unfortunately, the support radius scaling is very large compared to the partition element diameter.
Each of the above approaches has significant shortcomings:  while \cite{F95} avoids convolution, computing Markov partitions is usually impractical; in practice, one usually does not implement the convolution in the methods of \cite{DJ99,BKL02} because of a large computational cost.

In the present paper, we develop two new rigorous approaches that are computationally practical.
Each is based on Fourier approximation, which is a natural basis for the periodic domain $\mathbb{T}^d$, and can exploit the smoothness in the map $T$.
Furthermore, the spatial Fourier basis need not be adapted to the unstable and stable directions of $T$;  we use a standard Cartesian coordinate system.
Our first scheme builds finite-rank approximations of the Perron-Frobenius operator in a two-step process involving convolution, but the convolution is cheaply implementable via Fourier methods.
Our second scheme builds a single sequence of finite-rank approximations of the Perron-Frobenius operator that avoids convolution altogether.
The SRB measure approximation is found as the leading eigenvector of these finite-rank estimates of the Perron-Frobenius operator (in frequency space), and is a real-analytic function close to the true SRB measure in the anisotropic norms of \cite{gouezel2006banach}.
We implement both approaches and find rigorous numerical approximations of the SRB measure for an example Anosov map.

In the Anosov setting the stability of the variance has been established for specific perturbations, usually via a combination of the Green-Kubo formula and spectral stability of the (untwisted) transfer operator. Our approach is essentially the same, although by considering the stability of the twisted transfer operator we additionally obtain the stability of the rate function. Stability of the variance is noted in \cite{BKL02} for stochastic perturbations with smooth kernels, and in \cite{gouezel2006banach} for smooth deterministic perturbations and a broad class of random perturbations. The variance has also been shown to be stable to deterministic perturbations of Lorenz flows \cite{bahsoun2018variance}.

Rigorous numerical estimation of the variance has so far been considered for smooth expanding circle maps \cite{pollicott2017rigorous, wormell2017spectral}, Lasota-Yorke maps \cite{bahsoun2016rigorous,cfstability}, the intermittent Liverani-Saussol-Vaienti interval map \cite{bahsoun2016rigorous}, and piecewise expanding multidimensional maps \cite{cfstability}.
We will show that our new Fourier-based constructions and the anisotropic spaces of \cite{gouezel2006banach} fit into the general stability framework developed in \cite{cfstability}, which yields stability of the CLT variance for a broad class of deterministic, stochastic, and numerical perturbations.
We thus obtain the first rigorous numerical scheme for estimating the variance for Anosov systems, and implement this scheme to produce numerical estimates for an example Anosov map.
Our implementation constructs a finite-rank approximation of the Perron-Frobenius operator (in frequency space) and exploits the smoothness of spectral data to only require the solution of a single linear equation to estimate the variance.

In the context of Anosov flows, the stability of the rate function with respect to smooth deterministic perturbations of the velocity field is noted in \cite{butterley2007smooth}.
We extend the general stability framework \cite{cfstability} for statistical properties of multidimensional piecewise expanding maps to the setting of Anosov diffeomorphisms, proving stability of the rate function with respect to standard classes of deterministic, stochastic, and numerical perturbations.
This yields a proof of uniform convergence of the rate function with respect to smooth deterministic perturbations of an Anosov map.
For each of the two new Fourier-based schemes, we prove uniform convergence of the rate function estimates arising from finite-rank approximations of twisted Perron-Frobenius operators.
We implement these schemes to provide rigorous numerical estimates of a rate function for an example observation and Anosov map.

An overview of the paper is as follows.
Section \ref{sec:twisted_spectral_stability} briefly reviews results from \cite{cfstability} concerning the spectral stability of twisted quasi-compact operators in the sense of \cite{keller1999stability}.
In Section \ref{sec:nagaev_guivarch} we state the abstract requirements for the Nagaev-Guivarc'h method to guarantee a central limit theorem and large deviation principle, and state an abstract hypothesis that yields stability of the variance and rate function using a modification of the abstract stability result from \cite{cfstability}.
In Section \ref{sec:stability_anosov} we briefly introduce the functional analytic setting of Gou{\"e}zel and Liverani \cite{gouezel2006banach} and verify the abstract hypothesis for Anosov maps in this setting.
We then obtain stability of the rate function for Anosov diffeomorphisms under deterministic perturbations; see Theorem \ref{thm:anosov_stat_stability}.
Section \ref{sec:approximation_anosov} introduces our first Fourier-analytic scheme, which has wide applicability due to additional mollification, allowing the scheme to ``smooth away'' the complications of hyperbolic dynamics.
Proposition \ref{prop:convolution_fourier_approx} verifies that there is a sequence of Fourier projections of mollified transfer operators that satisfy the abstract stability hypotheses and Theorem \ref{thm:twisted_anosov_numeric_stability} states the corresponding stability results.
In Section \ref{sec:fejer_anosov} we introduce our second scheme, which removes the mollification step for Anosov maps with  approximately constant stable and unstable directions;  we obtain a pure (weighted) Fourier projection method based on the Fej{\'e}r kernel.
We may subject the transfer operator to stochastic kernels with relaxed regularity requirements (compared to Section \ref{sec:approximation_anosov}) and full (but sufficiently decaying) support and still obtain the abstract stability hypotheses;  see Proposition \ref{prop:stochastic_perturbation_kl}.
Theorem \ref{thm:stochastic_perturb_stability} summarises the stability results that are obtained, and Corollary \ref{cor:fejer_stability} addresses the perturbation arising from approximation with Fej{\'e}r kernels.
We conclude in Section \ref{sec:est_stat_prop_anosov} by implementing the above schemes on a nonlinear perturbation of Arnold's cat map.
We compare our Fourier-based estimates of the SRB measure and the variance with estimates from a (non-rigorous) pure Ulam method, and compute Fourier-based estimates of the rate function.
To the authors' knowledge, this constitutes the first rigorous computation of the variance and rate function of an Anosov map.

\section{A review of the spectral stability of twisted quasi-compact operators}\label{sec:twisted_spectral_stability}

In this section we review the theory of spectral stability of quasi-compact operators from \cite{keller1999stability} and the corresponding results for twisted operators from \cite{cfstability}. Let $(E, \norm{\cdot})$ be a (complex) Banach space, denote by $L(E)$ the bounded linear operators on $E$, and let $A \in L(E)$. Denote the spectrum of $A$ by $\sigma(A)$. The essential spectrum of $A$ is
\begin{equation*}
  \esigma(A) = \{ \omega \in \sigma(A) : \text{$\omega$ is not an eigenvalue of $A$ of finite algebraic multiplicity}  \}.
\end{equation*}
Denote the spectral radius and essential spectral radius of $A$ by $\rho(A)$ and $\erho(A)$, respectively. We say that $A$ is quasi-compact if $\erho(A) < \rho(A)$. If $A$ is quasi-compact and $\sigma(A) \cap \{ \omega : \abs{\omega} = \rho(A)\}$ consists of a single simple eigenvalue $\lambda$ then we call $A$ \emph{simple} and $\lambda$ the \emph{leading eigenvalue} of $A$. In this case $A$ has decomposition \cite[III.6.4-5]{kato1966perturbation}
\begin{equation}\label{eq:quasi-compact_decomp}
  A = \lambda \Pi + N,
\end{equation}
where $\Pi$ is the rank-one eigenprojection corresponding to $\lambda$, $N \in L(E)$ is such that $\rho(N) < \rho(A)$, and $N\Pi = \Pi N = 0$. We call \eqref{eq:quasi-compact_decomp} the quasi-compact decomposition of $A$.

Let $\wnorm{\cdot}$ be a norm on $E$ such that the closed, unit ball in $(E, \norm{\cdot})$ is relatively compact in the topology of $\wnorm{\cdot}$. After possibly scaling $\wnorm{\cdot}$, we may assume that $\wnorm{\cdot} \le \norm{\cdot}$.
Define the norm $\tnorm{\cdot}$ on $L(E)$ by
\begin{equation*}
  \tnorm{A} = \sup_{\norm{f} = 1} \wnorm{Af}.
\end{equation*}
It is classical that if $A$ is a simple quasi-compact operator and $A'$ is sufficiently close to $A$ in $L(E)$ then $A'$ is also a simple quasi-compact operator with leading eigenvalue close to  $\lambda$ (see e.g. \cite[IV.3.5]{kato1966perturbation}).
However, the condition of closeness in $L(E)$ is seldom satisfied by applications in dynamical systems. In \cite{keller1999stability}, Keller and Liverani showed that if $A'$ is close to $A$ in the weaker topology of $\tnorm{\cdot}$, both operators obey a Lasota-Yorke inequality, and growth restrictions are placed on the (various) operator norms of iterates of $A$ and $A'$, then one can recover appropriately modified versions of the spectral stability results from operator norm based perturbation theory. We now detail the requirements for these results, referring to \cite{keller1999stability} for exact statements.

\begin{definition}\label{KLdefn}
  A family of operators $\{ A_\epsilon \}_{\epsilon \ge 0} \subseteq L(E)$ satisfies the Keller-Liverani (KL) condition if each of the following is satisfied:
  \begin{enumerate}[label=(KL\arabic*)]
    \item \label{en:kl_conv} There exists a monotone upper-semicontinuous function $\tau: [0, \infty) \to [0, \infty)$ such that $\tau(\epsilon) > 0$ whenever $\epsilon > 0$, $\tnorm{A_\epsilon - A_0} \le \tau(\epsilon)$, and $\lim_{\epsilon \to 0} \tau(\epsilon) = 0$.
    \item \label{en:kl_l1_bound} There exists $C_1, K_1 > 0$ such that $\wnorm{A_\epsilon^n} \le C_1 K_1^n$ for every $\epsilon \ge 0$ and $n \in \N$.
    \item \label{en:kl_ly_bound} There exists $C_2,C_3, K_2 > 0$ and $\alpha \in (0,1)$ such that
    \begin{equation}\label{eq:kl_ly}
      \norm{A_\epsilon^n f} \le C_2 \alpha^n \norm{f} + C_3 K_2^n \wnorm{f}
    \end{equation}
    for every $\epsilon \ge 0$, $f \in E$ and $n \in \N$.
  \end{enumerate}
\end{definition}

\begin{remark}
  If rates of convergence are not required, then it suffices to prove that $\tnorm{A_\epsilon - A_0} \to 0$ instead of constructing the function $\tau$ in (KL1). After possibly passing to a sub-family $\{ A_\epsilon\}_{\epsilon \in [0, \epsilon')}$ for some $\epsilon ' > 0$, the two conditions are equivalent.
  We will use this fact frequently without further comment.
\end{remark}

Instead of considering a single simple quasi-compact operator, one sometimes considers an analytic operator-valued map $A(\cdot) : D \to L(E)$, where $D \subseteq \C$ is an open neighbourhood of 0 and $A(0)$ is quasi-compact and simple. In this case, classical analytic perturbation theory for linear operators \cite{kato1966perturbation} posits the existence of some $\delta> 0$ such that $A(z)$ is a simple quasi-compact operator for each $z \in D_\delta = \{ \omega : \abs{\omega} < \delta\}$.
Moreover, the quasi-compact decomposition of $A(z)$ depends analytically on $z$ i.e. there are analytic maps $\lambda(\cdot): D_\delta \to \C$, $\Pi(\cdot): D_\delta \to L(E)$, and $N(\cdot): D_\delta \to L(E)$ such that $A(z)$ has quasi-compact decomposition $A(z) = \lambda(z) \Pi(z) + N(z)$. In \cite{cfstability} the authors considered the question of spectral stability of such analytic operator-valued maps under conditions similar to (KL) and when the analytic families are induced by a `twist'.

\begin{definition}\label{def:twist}
  If $M: D \to L(E)$ is analytic on an open neighbourhood $D \subseteq \C$ of 0 and $M(0)$ is the identity, then we call $M$ a \emph{twist}. If $A \in L(E)$ then the operators $A(z) := A M(z)$ are said to be \emph{twisted} by $M$. We say that $M$ is compactly $\wnorm{\cdot}$-bounded if for every compact $V \subseteq D$ we have
  \begin{equation*}
    \sup_{z \in V} \wnorm{M(z)} < \infty.
  \end{equation*}
\end{definition}

We state a version of \cite[Theorem 2.6]{cfstability} concerning the spectral stability of twisted quasi-compact operators.

\begin{theorem}[{\cite[Theorem 2.6]{cfstability}}]
  \label{theorem:convergence_of_eigendata_derivs}
  Let $\{ A_\epsilon \}_{\epsilon \ge 0}$ satisfy (KL), where $A_0$ is a simple quasi-compact operator with leading eigenvalue $\lambda_0$ satisfying $\alpha < \abs{\lambda_0}$, and let $M : D \to \C$ be a compactly $\wnorm{\cdot}$-bounded twist.
  Then there exists $\theta, \epsilon' > 0$ and, for each $\epsilon \in [0, \epsilon']$, analytic functions $\lambda_\epsilon(\cdot): D_\theta \to \C$, $\Pi_\epsilon(\cdot): D_\theta \to L(E)$, and $N_\epsilon(\cdot) : D_\theta \to L(E)$ such that $A_\epsilon(z)$ is a simple quasi-compact operator with decomposition $A_\epsilon(z) = \lambda_\epsilon(z) \Pi_\epsilon(z) + N_\epsilon(z)$ whenever $z \in D_\theta$.
  Additionally, for each $n \in \N$ we have the following convergence as $\epsilon \to 0$ on $D_\theta$:
  \begin{enumerate}
    \item $\lambda_\epsilon^{(n)}(\cdot)$ converges uniformly to $\lambda_0^{(n)}(\cdot)$.
    \item $\Pi_\epsilon^{(n)}(\cdot)$ converges uniformly to $\Pi_0^{(n)}(\cdot)$ in $\tnorm{\cdot}$.
    \item $N_\epsilon^{(n)}(\cdot)$ converges uniformly to $N_0^{(n)}(\cdot)$ in $\tnorm{\cdot}$.
  \end{enumerate}
\end{theorem}

We finish this section with a result concerning the robustness of the condition (KL) to perturbations that are simultaneously small in the operator norms $\norm{\cdot}$ and $\wnorm{\cdot}$, which we will frequently use in the sequel.

\begin{proposition}\label{prop:kl_norm_perturb}
  Suppose that $\{A_\epsilon\}_{\epsilon \ge 0}$ satisfies (KL) and $\{B_\epsilon\}_{\epsilon \in [0,\epsilon_1)} \subseteq L(E)$ satisfies $B_0 = 0$, $\lim_{\epsilon \to 0} \norm{B_\epsilon} = 0$ and $\sup_{\epsilon \in [0,\epsilon_1)} \wnorm{ B_\epsilon} < \infty$.
  Then there exists $\epsilon_2 \in (0,\epsilon_1)$ so that $\{A_\epsilon + B_\epsilon\}_{\epsilon \in [0, \epsilon_2)}$ satisfies (KL).
  \begin{proof}
    We prove (KL1), (KL2) and (KL3) separately.
    \paragraph{(KL1)}
    As $\norm{B_\epsilon} \to 0$, there exists $\epsilon' > 0$ so that $\sup_{\epsilon \in [0, \epsilon')} \norm{B_\epsilon} < \infty$.
    As $\{A_\epsilon\}_{\epsilon \ge 0}$ satisfies (KL1), $\{B_\epsilon \}_{\epsilon \in [0, \epsilon')}$ is bounded in $L(E)$, $\norm{B_\epsilon} \to 0$, and
    \begin{equation*}
      \tnorm{A_\epsilon + B_\epsilon - A_0} \le \tnorm{A_\epsilon - A_0} + \norm{B_\epsilon},
    \end{equation*}
    it is clear that (KL1) is satisfied.
    \paragraph{(KL2)}
    As $\{A_\epsilon\}_{\epsilon \ge 0}$ satisfies (KL2) and $\sup_{\epsilon \in [0,\epsilon_1)} \wnorm{B_\epsilon} < \infty$, we have
    \begin{equation*}
      \sup_{\epsilon \in [0,\epsilon_1)} \wnorm{A_\epsilon + B_\epsilon} \le C_1 K + \sup_{\epsilon \in [0,\epsilon_1)} \wnorm{B_\epsilon} < \infty.
    \end{equation*}
    The required bound follows by iterating this inequality.
    \paragraph{(KL3)}
    By expanding $(A_{\epsilon} + B_{\epsilon})^n$, applying a counting argument, and using (KL3) for $\{ A_{\epsilon}\}_{\epsilon \ge 0}$ we have
    \begin{equation}\begin{split}\label{eq:kl_norm_perturb_1}
      \norm{(A_{\epsilon} + B_{\epsilon})^n f} &\le \norm{A_{\epsilon}^n f} + \sum_{k=0}^{n-1} \frac{n!}{k!(n-k)!} \norm{A_{\epsilon}}^k \norm{B_{\epsilon}}^{n-k} \norm{f} \\
      &\le C_2 \alpha^n \norm{f} + C_3 K_2^n\wnorm{f} + 2^{n} \sum_{k=0}^{n-1} \norm{A_{\epsilon}}^k \norm{B_{\epsilon}}^{n-k} \norm{f}.
    \end{split}\end{equation}
    Let $\beta \in (\alpha, 1)$ and choose $n \in \N$ so that $C_1 \alpha^n < \beta^n$. As (KL3) holds for $\{A_\epsilon\}_{\epsilon \ge 0}$, it follows that $\{A_\epsilon\}_{\epsilon \ge 0}$ is bounded in $L(E)$.
    Hence, as $\norm{B_\epsilon} \to 0$, there exists $\epsilon'' > 0$ so that for every $\epsilon \in [0, \epsilon'')$ we have
    \begin{equation*}\label{eq:kl_norm_perturb_3}
      C_1 \alpha^n + 2^{n} \sum_{k=0}^{n-1} \norm{A_{\epsilon}}^k \norm{B_{\epsilon}}^{n-k} \le \beta^n.
    \end{equation*}
    By applying \eqref{eq:kl_norm_perturb_3} to \eqref{eq:kl_norm_perturb_1} for all $\epsilon \in [0, \epsilon'')$ we have
    \begin{equation}\label{eq:kl_norm_perturb_2}
      \norm{(A_{\epsilon} + B_{\epsilon})^n f} \le \beta^n \norm{f} + C_3 K_2^n\wnorm{f}.
    \end{equation}
    Using (KL2) for $\{ A_\epsilon + B_\epsilon\}_{\epsilon \in [0,\epsilon_1)}$ one may iterate \eqref{eq:kl_norm_perturb_2} to obtain for all $\epsilon \in [0, \epsilon'')$ and $k \in \Z^+$ that
    \begin{equation}\label{eq:kl_norm_perturb_4}
      \norm{(A_{\epsilon} + B_{\epsilon})^{nk} f} \le \beta^{nk} \norm{f} + C_4 K_3^{nk}\wnorm{f},
    \end{equation}
    where $C_4, K_3 > 0$ are independent of $k$ and $\epsilon$.
    In the present setting, $\sup_{\epsilon \in [0, \epsilon'')} \norm{A_\epsilon + B_\epsilon}$ is finite.
    Using this fact one easily obtains (KL3) from \eqref{eq:kl_norm_perturb_4}.
  \end{proof}
\end{proposition}

\section{Stability of the statistical properties of dynamical systems via the naive Nagaev-Guivarc'h method}\label{sec:nagaev_guivarch}

In \cite{cfstability} it is shown that when a dynamical system satisfies a central limit theorem or large deviation principle via the Nagaev-Guivarc'h method, then the variance and rate function associated with these statistical laws are stable under perturbations satisfying (KL). However, the the method as formulated in \cite{cfstability} (see also \cite{hennion2001limit, Aimino2015}) is occasionally inapplicable due to some technical requirements; for example, it is required that the Perron-Frobenius operator is quasi-compact on a Banach algebra, and that the Banach space in question is a subspace of $L^1$. Neither condition is generally verified by the transfer operator associated to Anosov maps on any known Banach space.
In \cite{gouezel2015limit} Gou{\"e}zel details the naive Nagaev-Guivarc'h method, which circumvents some of these shortcomings with a more abstract formulation of the method. In this section we recall how a CLT and LDP may be obtained via the naive Nagaev-Guivarc'h method.
In Theorem \ref{thm:stability_of_statistical_data} we show that when the method is compatible with the stability theory for twisted quasi-compact operators detailed in Section \ref{sec:twisted_spectral_stability}, then the variance of the CLT and the rate function of the LDP are stable to perturbations of the type (KL).
In what follows $\{Y_k\}_{k \in \N}$ is a sequence of real random variables with partial sums $S_n = \sum_{k=0}^{n-1} Y_k$ satisfying $\lim_{n \to \infty} \E(S_n)/n = 0$.

\begin{theorem}[Central limit theorem {\cite[Theorem 2.4]{gouezel2015limit}}]\label{thm:naive_ng_clt}
  If there exist a Banach space $E$, operator-valued map $A(\cdot) : I \to L(E)$, where $I$ is a real open neighbourhood of 0, and $\zeta \in E$, $\nu \in E^*$ such that:
  \begin{enumerate}
    \item $A(0)$ is a simple quasi-compact operator with $\rho(A(0)) = 1$.
    \item The mapping $t \mapsto A(t)$ is $\mathcal{C}^2$.
    \item $\E(e^{i t S_n}) = \nu A(t)^n \zeta$ for all $n \in \N$ and $t \in I$.
  \end{enumerate}
  Then $\{Y_k\}_{k \in \N}$ satisfies a CLT: there exists $\sigma^2 \ge 0$ such that $S_n/\sqrt{n}$ converges in distribution to a $N(0, \sigma^2)$ random variable as $n \to \infty$.
\end{theorem}

\begin{theorem}[{Large deviation principle \cite[Remark 2.3]{Fernando2018}}]\label{thm:naive_ng_ldp}
  If there exist a Banach space $E$, operator-valued map $A(\cdot) : I \to L(E)$, where $I$ is a real open neighbourhood of 0, and $\zeta \in E$, $\nu \in E^*$ such that:
  \begin{enumerate}
    \item $A(0)$ is a simple quasi-compact operator with $\rho(A(0)) = 1$.
    \item The mapping $t \mapsto A(t)$ is $\mathcal{C}^1$ and $t \mapsto \ln\rho(A(t))$ is strictly convex in some neighbourhood of 0.
    \item $\E(e^{t S_n}) = \nu A(t)^n \zeta$ for all $n \in \N$ and $t \in I$.
  \end{enumerate}
  Then $\{Y_k\}_{k \in \N}$ satisfies a LDP: there exists a non-negative, continuous and convex rate function $r : J \to \R$, where $J$ is an open neighbourhood of $0$, such that for every $\epsilon \in J \cap (0, \infty)$ we have
  \begin{equation*}
    \lim_{n \to \infty} \frac{1}{n} \ln \Prob(S_n \ge n \epsilon) = - r(\epsilon).
  \end{equation*}
\end{theorem}

\begin{remark}
  As stated, Theorem \ref{thm:naive_ng_ldp} differs slightly to the result in \cite[Remark 2.3]{Fernando2018}, however a straightforward modification of the arguments from \cite{Fernando2018} readily yields Theorem \ref{thm:naive_ng_ldp}.
\end{remark}

Both the CLT and LDP are parameterised, and under the settings of Theorems \ref{thm:naive_ng_clt} and \ref{thm:naive_ng_ldp} these parameters are determined by the spectral data of $A(t)$ as follows.
As $t \mapsto A(t)$ is $\mathcal{C}^k$ ($k = 1, 2$) and $A(0)$ is a simple quasi-compact operator, by \cite[Proposition 2.3]{gouezel2015limit} there exists $\theta > 0$ and $\mathcal{C}^k$ maps $\lambda(\cdot): (-\theta, \theta) \to \C$, $\Pi(\cdot): (-\theta, \theta) \to L(E)$, and $N(\cdot): (-\theta, \theta) \to L(E)$ such that for for $t \in (-\theta, \theta)$ the operator $A(t)$ is quasi-compact and simple with leading eigenvalue $\lambda(t)$ and decomposition $A(t) = \lambda(t) \Pi(t) + N(t)$.
The variance of CLT is
\begin{equation}\label{eq:var_formula}
  \sigma^2 = \lambda^{(2)}(0),
\end{equation}
and the rate function of the LDP is
\begin{equation*}
  r(s) = \sup_{t \in (-\theta, \theta)} (st - \ln\abs{\lambda(t)}).
\end{equation*}
Moreover, due to the strict convexity and continuous  differentiability of $t \mapsto \ln\abs{\lambda(t)}$ on $(-\theta,\theta)$, and the application of the local Gartner-Ellis Theorem \cite[Lemma XIII.2]{hennion2001limit} used to obtain Theorem \ref{thm:naive_ng_ldp}, we have that the domain of the rate function is
\begin{equation}\label{eq:domain_of_rate_function}
    \left(\restr{\frac{\mathrm{d}}{\mathrm{d}z} \ln\abs{\lambda(z)}}{z = -\theta}, \restr{\frac{\mathrm{d}}{\mathrm{d}z} \ln\abs{\lambda(z)}}{z = \theta}\right).
\end{equation}

In Theorem \ref{thm:naive_ng_clt} the characteristic function of $S_n$ is encoded by $\nu A(t)^n \zeta$ whilst in Theorem \ref{thm:naive_ng_ldp} it is the moment-generating function of $S_n$ that is encoded.
These settings are frequently unified by the following hypothesis:
\begin{hypothesis}\label{hypothesis:analytic_coding}
  Suppose that $\{Y_k\}_{k \in \N}$ is a sequence of real random variables with partial sums $S_n = \sum_{k=0}^{n-1} Y_k$ satisfying $\lim_{n \to \infty} \E(S_n)/n = 0$.
  We say that $\{Y_k\}_{k \in \N}$ satisfies Hypothesis \ref{hypothesis:analytic_coding} if there exists a Banach space $(E, \norm{\cdot})$, $\zeta \in E$, $\nu \in E^*$, and an analytic operator-valued map $A : D \mapsto L(E)$, where $D \subseteq \C$ is an open neighbourhood of $0$, such that $A(0)$ is a simple quasi-compact operator, $\rho(A(0)) = 1$, $t \mapsto \rho(A(t))$ is strictly convex in some real neighbourhood of 0, and
  \begin{equation*}
    \E(e^{z S_n}) = \nu A(z)^n \zeta
  \end{equation*}
  for every $z \in D$ and $n \in \N$.
\end{hypothesis}

It is clear that if Hypothesis \ref{hypothesis:analytic_coding} holds then Theorems \ref{thm:naive_ng_clt} and \ref{thm:naive_ng_ldp} both hold i.e. $\{Y_k\}_{k \in \N}$ satisfies a CLT and LDP. Hypothesis \ref{hypothesis:analytic_coding} is frequently verified by applications of the naive Nagaev-Guivarc'h method to dynamical systems. In addition, in these applications the map $z \mapsto A(z)$ arises from a twist i.e. there exists a twist $M : D \mapsto L(E)$ such that $A(z) = A(0) M(z)$. In this case we can apply Theorem \ref{theorem:convergence_of_eigendata_derivs} to obtain stability of the variance and rate function with respect to perturbations satisfying (KL).

\begin{theorem}\label{thm:stability_of_statistical_data}
  Let $\{Y_k\}_{k \in \N}$ be a sequence of real random variables satisfying Hypothesis \ref{hypothesis:analytic_coding}.
  Suppose that $\wnorm{\cdot}$ is a second norm on $E$ so that the closed, unit ball in $(E,\norm{\cdot})$ is relatively compact with respect to $\wnorm{\cdot}$ and that there exists a compactly $\wnorm{\cdot}$-bounded twist $M: D \to L(E)$ such that $A(z) = A(0)M(z)$.
  If $\{A_\epsilon\}_{\epsilon \ge 0}$ satisfies (KL), where $A_0 = A(0)$, then there exists $\theta, \epsilon' > 0$ and, for every $\epsilon \in [0 , \epsilon']$, analytic maps $\lambda_\epsilon(\cdot) : D_\theta \to \C$, $\Pi_\epsilon(\cdot) : D_\theta \to L(E)$ and $N_\epsilon(\cdot) : D_\theta \to L(E)$ such that for every $\epsilon \in [0, \epsilon']$ and $z \in D_\theta$ the operator $A_\epsilon(z)$ has quasi-compact decomposition $A_\epsilon(z) = \lambda_\epsilon(z) \Pi_\epsilon(z) + N_\epsilon(z)$ and as $\epsilon \to 0$ the maps $\lambda_\epsilon(\cdot)$, $\Pi_\epsilon(\cdot)$ and $N_\epsilon(\cdot)$ converge as in Theorem \ref{theorem:convergence_of_eigendata_derivs}.
  Moreover, we have stability of the parameters of the CLT and LDP for $\{Y_k\}_{k \in \N}$ in the following sense:
  \begin{enumerate}
    \item The variance is stable: $\lim_{\epsilon \to 0} \lambda_\epsilon^{(2)}(0) = \sigma^2$.
    \item The rate function is stable: for each sufficiently small compact subset $U$ of the domain of the rate function $r_g$ there exists an interval $V \subseteq (-\theta, \theta)$ so that
    \begin{equation*}
      \lim_{\epsilon \to 0} \sup_{z \in V} (sz -\log \abs{\lambda_\epsilon(z)}) = r(s)
    \end{equation*}
    uniformly on $U$.
  \end{enumerate}
  \begin{proof}
    By Theorem \ref{theorem:convergence_of_eigendata_derivs} there exists $\theta, \epsilon' > 0$ and, for each $\epsilon \in [0,\epsilon']$, maps $\lambda_\epsilon(\cdot), \Pi_\epsilon(\cdot)$ and $N_\epsilon(\cdot)$ as required for Theorem \ref{thm:stability_of_statistical_data}.
    The proof of the stability of the variance follows from \eqref{eq:var_formula} and Theorem \ref{theorem:convergence_of_eigendata_derivs}.
    After possibly reducing the value of $\theta$, so that the domain of the rate function is the one given in Theorem \ref{thm:naive_ng_ldp}, the stability of the rate function follows from \cite[Theorem 3.8 and Proposition 3.9]{cfstability} with some minor modifications that we will now discuss. In \cite{cfstability} it is required that $A_\epsilon(z)$ is positive for each $\epsilon \ge 0$ and $z \in \R$ so that $\lambda_\epsilon(z)$ is positive too; we remove this assumption and consequently deal with $\abs{\lambda_\epsilon(z)}$ instead.
    By Theorem \ref{theorem:convergence_of_eigendata_derivs} and the reverse triangle inequality we have uniform convergence of $\abs{\lambda_\epsilon(\cdot)}$ to $\abs{\lambda_0(\cdot)}$ on each compact subset of $D_\theta$.
    Although \cite[Theorem 3.8]{cfstability} proves a H{\"o}lder estimate for the convergence of the rate functions, straightforward modifications to the proof yields uniform convergence as is required for Theorem \ref{thm:stability_of_statistical_data}.
    The proof of \cite[Proposition 3.10]{cfstability} holds upon replacing \cite[Lemma 3.9]{cfstability} with the condition that $z \mapsto \ln\rho(\LL(z))$ is strictly convex in a real neighbourhood of $0$,
    and noting that $z \mapsto \ln \abs{\lambda_0(z)}$ is $\mathcal{C}^1$ on $(-\theta, \theta)$.
  \end{proof}
\end{theorem}

\section{Stability of statistical limit laws for Anosov maps}\label{sec:stability_anosov}

It is classical that topologically transitive Anosov diffeomorphisms satisfy a central limit theorem and large deviation principle for sufficiently smooth observables: such results were first established by using Markov partitions to reduce to the case of subshifts of finite type \cite{orey1989, ratner1973central}.
In order to apply the stability results from Section \ref{sec:nagaev_guivarch} to Anosov maps we require that these limit laws hold due to the naive Nagaev-Guivarc'h method; verifying this is the main point of this section.
In particular, using the functional analytic setup in \cite{gouezel2006banach}, we confirm that $\{g \circ T^k \}_{k \in \N}$ satisfies Hypothesis \ref{hypothesis:analytic_coding} for $T$ a $\mathcal{C}^{r+1}$ Anosov map, $r > 1$, and appropriate observables $g$.
Theorem \ref{thm:stability_of_statistical_data} then yields stability of the variance and rate function to perturbations of type (KL), which forms the basis for our results in Sections \ref{sec:approximation_anosov}, \ref{sec:fejer_anosov} and \ref{sec:est_stat_prop_anosov}.

\subsection{The functional analytic setup of Gou{\"e}zel and Liverani}
\label{sec:func_analytic_anosov}

We review the functional analytic setup of \cite{gouezel2006banach}.
Let $d > 1$ and $X$ be a $d$-dimensional, $\mathcal{C}^\infty$, compact, connected Riemannian manifold and $T \in \mathcal{C}^{r+1}(X,X)$, $r > 1$, be an Anosov map. By using an adapted metric we may assume that $\nu_u > 1$ is everywhere less than the local expansion of $T$ in the unstable direction, and $\nu_s < 1$ is everywhere greater than the local contraction of $T$ in the stable direction (see \cite[Proposition 5.2.2]{brin2002introduction}).
We review the construction of such a metric in Section \ref{sec:cond_stochastic_stab_anosov}.
When $T$ has a unique SRB measure we denote it by $\mu$.
Let $\Omega$ denote both the Riemannian measure on $X$ and the linear functional $f \mapsto \intf f d\Omega$. The transfer operator $\LL : \mathcal{C}^{r}(X, \R) \to \mathcal{C}^{r}(X, \R)$ associated with $T$ is defined by
\begin{equation}\label{eq:pf_definition}
  \intf (\LL h) \cdot u d\Omega = \intf h \cdot (u \circ T) d\Omega,
\end{equation}
where $u,h \in \mathcal{C}^{r}(X, \R)$.
For $\LL$ to have `good' spectral properties it is necessary to consider it as an operator on an appropriately chosen anisotropic Banach space.
We now describe the construction of such a space from \cite{gouezel2006banach}. Denote by $\Sigma$ the set of admissible leaves (see \cite[Section 3]{gouezel2006banach} for the full definition).
For each $W \in \Sigma$ we denote the collection of $\mathcal{C}^{r}$ vector fields that are defined on a neighbourhood of $W$ by $\mathcal{V}^r(W)$, and by $\mathcal{C}^q_0(W, \R)$ the set of functions in $\mathcal{C}^q(W, \R)$ that vanish on a neighbourhood of $\partial W$.
For $h \in \mathcal{C}^r (X, \R)$, $q > 0$, $p \in \N$ with $p \le r$ let
\begin{equation*}
  \norm{h}^{-}_{p,q} =
  \sup_{W \in \Sigma} \,
  \sup_{\substack{v_1, \dots v_p \in \mathcal{V}^r(W) \\ \abs{v_i}_{\mathcal{C}^r} \le 1}} \,
  \sup_{\substack{\varphi \in \mathcal{C}^q_0(W, \R) \\ \abs{\varphi}_{\mathcal{C}^q} \le 1}} \,
  \int_W (v_1 \dots v_p h) \cdot \varphi .
\end{equation*}
Then
\begin{equation}\label{eq:real_norm}
 \norm{h}_{p,q} = \sup_{0 \le k \le p} \norm{h}^{-}_{k,q+k} = \sup_{p'\le p, q' \ge q + p'} \norm{h}^{-}_{p',q'}
\end{equation}
is a norm on $\mathcal{C}^{r}(X, \R)$. Denote by $B^{p,q}$ the completion of $\mathcal{C}^r (X, \R)$ under this norm. As the naive Nagaev-Guivarc'h method requires a complex Banach space, we consider the complexification $B^{p,q}_\C$ of the spaces $B^{p,q}$. When endowed with the norm\footnote{We abuse notation and denote the norm on $B^{p,q}_\C$ by $\norm{\cdot}_{p,q}$.}
\begin{equation}\label{eq:complex_norm}
  \norm{h_r + i h_i}_{p,q} =  \max \{\norm{h_r}_{p,q}, \norm{h_i}_{p,q} \},
\end{equation}
$B^{p,q}_\C$ is a complex Banach space. It is on this space that the operator $\LL$ is quasi-compact.

\begin{theorem}[{\cite[Theorem 2.3]{gouezel2006banach}}]
  If $p \in \Z^+$ and $q > 0$ satisfy $q + p < r$ then the operator $\LL: B^{p,q}_\C \to B^{p,q}_\C$ has spectral radius one. In addition, $\LL$ is quasi-compact with $\erho(\LL) \subseteq \{ \omega \in \C : \abs{\omega} \le \max\{\nu_u^{-p}, \nu_s^{q} \} \}$.
  Moreover, the eigenfunctions corresponding to eigenvalues of modulus 1 are distributions of order 0, i.e., measures. If the map is topologically transitive, then 1 is a simple eigenvalue and no other eigenvalues of modulus one are present.
\end{theorem}

\subsection{The naive Nagaev-Guivarc'h method for Anosov maps}

The main technical result of this section is Proposition \ref{prop:anosov_hypothesis}, which verifies Hypothesis \ref{hypothesis:analytic_coding} in the setting of Section \ref{sec:func_analytic_anosov}.
As $1 \in B^{p,q}_\C$, for any $g \in \mathcal{C}^r(X, \R)$ we may define $e^{zg}$ by the power series $\sum_{k=0}^\infty z^k g^k /k!$.
We define $M_g(\cdot): \C \to L(B^{p,q}_\C)$ by setting $M_g(z)(h) = e^{zg} h$ for $h \in \mathcal{C}^r(X,\C)$ and then passing to $B^{p,q}_\C$ by density.

\begin{proposition}\label{prop:defn_of_twist_anosov}
  Let $p \in \Z^+$, $q > 0$ satisfy $p + q < r$. If $g \in \mathcal{C}^{r}(X, \R)$ and $M_g : \C \to L(B^{p,q}_\C)$ is defined by $M_g(z)(f) = e^{zg} f$, then $M_g$ is a compactly $\norm{\cdot}_{p-1, q+1}$-bounded twist.
\end{proposition}

By \cite[Proposition 4.1]{gouezel2006banach} there is a continuous injection from $B^{p,q}_\C$ into $\mathcal{D}_q'(X)$, the distributions of order at most $q$. With this injection in mind, we can consider some elements of $B^{p,q}_\C$ as probability measures (exactly when they are probability measures in $\mathcal{D}_q'(X)$).
Recall that $g \in \mathcal{C}^r(X, \R)$ is an $L^2(\mu)$-coboundary if there is $\phi \in L^2(\mu)$ so that $g = \phi - \phi \circ T$.

\begin{proposition}\label{prop:anosov_hypothesis}
  Assume that $\LL$ is a simple quasi-compact operator. Let $p \in \Z^+$ and $q > 0$ satisfy $q + p < r$. Let $m \in B^{p,q}$ be a probability measure, and suppose $g \in \mathcal{C}^{r}(X,\R)$ satisfies $\intf g d\mu = 0$ and is not an $L^2(\mu)$-coboundary.
  Define $\LL(\cdot) : \C \to L(B^{p,q}_\C)$ by $\LL(z) = \LL M_g(z)$.
  Then $\{ g \circ T^k \}_{k \in \N}$, when considered on the probability space $(X,m)$, satisfies Hypothesis \ref{hypothesis:analytic_coding} with Banach space $B^{p,q}_\C$ and operator-valued map $\LL(\cdot)$.
\end{proposition}

With Proposition \ref{prop:anosov_hypothesis} in hand, by Theorems \ref{thm:naive_ng_clt} and \ref{thm:naive_ng_ldp} we immediately obtain a CLT and LDP (on appropriate probability spaces) for $\{g \circ T^k\}_{k \in \N}$ whenever $g \in \mathcal{C}^{r}(X, \R)$ satisfies $\intf g d\mu = 0$ and is not an $L^2(\mu)$-coboundary.
Moreover, as $M_g$ is a compactly $\norm{\cdot}_{p-1,q+1}$-bounded twist and the unit ball in $B^{p,q}_\C$ is relatively compact in $\norm{\cdot}_{p-1,q+1}$ \cite[Lemma 2.1]{gouezel2006banach}, Theorem \ref{thm:stability_of_statistical_data} yields stability of the variance and rate function with respect to perturbations of the type (KL).
In Section \ref{sec:est_stat_prop_anosov} we use this fact to numerically approximate the statistical properties of Anosov maps on tori. For the moment we state following application, which extends \cite[Theorem 2.8, Remark 2.11]{gouezel2006banach}. For $y \in \R$ we denote by $\tau_y$ the map $x \mapsto x+y$.

\begin{theorem}[Stability of the rate function under deterministic perturbations]\label{thm:anosov_stat_stability}
  Let $T(\cdot) \in \mathcal{C}^1([0, 1], \mathcal{C}^{r+1}(X,X))$ be such that $T(0)$ is a topologically transitive Anosov diffeomorphism. Let $p \in \Z^+$ and $q > 0$ satisfy $p + q < r$.
  Fix a probability measure $m \in B^{p,q}$ and suppose $g \in \mathcal{C}^r(X,\R)$ satisfies $\intf g d\mu = 0$ and is not an $L^2(\mu)$-coboundary.
  There exists $\epsilon > 0$ and, for each $t \in [0, \epsilon]$, a number $A_t$ and map $r_t : J - A_t \to \R$, where $J$ is an open real neighbourhood of 0, so that $\{ g \circ T(t)^k - A_t \}_{k\in \N}$ satisfies a LDP on $(X,m)$ with rate function $r_t$, $A_t \to A_0 = 0$ and $r_t \circ \tau_{-A_t} \to r_0$ compactly on $J$.
\end{theorem}

We defer the proofs of Propositions \ref{prop:defn_of_twist_anosov} and \ref{prop:anosov_hypothesis}, and Theorem \ref{thm:anosov_stat_stability} to Appendix \ref{sec:proofs_sec_4}.

\section{Approximating the statistical data of Anosov maps}\label{sec:approximation_anosov}

In this section we introduce a scheme for approximating the spectrum of the Perron-Frobenius operator associated to an Anosov map on the $d$-dimensional torus $\mathbb{T}^d$, which we identify with $\R^d / \Z^d$.
The scheme proceeds by convolving the Perron-Frobenius operator with a compactly supported mollifier, and then approximating the `smoothened' operator using Fourier series. A similar idea is developed in \cite{BKL02}, where Ulam's method is considered instead of Fourier series and convergence of the SRB measure and variance are obtained under this scheme. We note that \cite{BKL02} did not include any computations (as we do, in Section \ref{sec:est_stat_prop_anosov}) nor did they consider the stability of the rate function.

Throughout this section we adopt the setting, assumptions and notation of Section \ref{sec:func_analytic_anosov}, and fix $p \in \Z^+$ and $q > 0$ satisfying $p + q < r$.
Let $\Leb$ denote the normalised Haar measure on $\mathbb{T}^d$.
For some $\epsilon_1 > 0$, suppose that the family of stochastic kernels $\{q_\epsilon \}_{\epsilon \in (0, \epsilon_1)} \subseteq \mathcal{C}^\infty(\mathbb{T}^d, \R)$ satisfies the following conditions:
\begin{enumerate}[label=(S\arabic*)]
  \item \label{en:S1} $q_\epsilon \ge 0$ and $\intf q_\epsilon d\Leb = 1$;
  \item \label{en:S2} The support of $q_\epsilon$ is contained in $B_\epsilon(0)$.
\end{enumerate}
For such a family we define operators $Q_\epsilon : \mathcal{C}^r(\mathbb{T}^d, \C) \to \mathcal{C}^r(\mathbb{T}^d, \C)$ by $Q_\epsilon f = f * q_\epsilon$. Recall that convolution is defined with respect to the Haar measure $\Leb$ on $\mathbb{T}^d$, which may differ from the measure $\Omega$ that is induced by the adapted metric.
It is evident, however, that the Radon-Nikodym derivatives $\frac{\mathrm{d}\Leb}{\mathrm{d}\Omega}$ and $\frac{\mathrm{d}\Omega}{\mathrm{d}\Leb}$ both exist, and are elements of $\mathcal{C}^{\infty}(\mathbb{T}^d, \R)$. As a consequence we obtain the following characterisation of $Q_\epsilon$:

\begin{lemma}\label{lemma:compact_convolution}
  $Q_\epsilon$ extends to a bounded operator $Q_\epsilon : B^{p,q}_\C \to \mathcal{C}^{\infty}(\mathbb{T}^d,\C)$. Consequently, $Q_\epsilon$ is compact as an element of $L(B^{p,q}_\C, \mathcal{C}^{m}(\mathbb{T}^d,\C))$ for every $m \in \Z^+$.
  Moreover, $Q_\epsilon$ is compact as an element of $L(B^{p,q}_\C)$ and, for each $m_1,m_2 \in \Z^+$, as an element of $L(\mathcal{C}^{m_1}(\mathbb{T}^d,\C), \mathcal{C}^{m_2}(\mathbb{T}^d,\C))$.
\end{lemma}

Let $\LL_0 = \LL$ and, for each $\epsilon \in (0, \epsilon_1)$, let $\LL_\epsilon = Q_\epsilon \LL_0$, which is in $L(B^{p,q}_\C)$ by virtue of the previous lemma.

\begin{lemma}\label{lemma:convolution_kl}
  There exists $\epsilon_2 \in (0, \epsilon_1)$ so that $\{ \LL_{\epsilon} \}_{\epsilon \in [0, \epsilon_2)}$ satisfies (KL) on $B^{p,q}_\C$ with $\wnorm{\cdot} = \norm{\cdot}_{p-1,q+1}$.
\end{lemma}

By Lemma \ref{lemma:compact_convolution}, for each $\epsilon \in (0, \epsilon_1)$ the operator $\LL_\epsilon$ is compact and, for every $m \in \Z^+$, maps the unit ball of $B^{p,q}_\C$ into a bounded subset of $\mathcal{C}^m(\mathbb{T}^d, \C)$.
For this reason $\LL_\epsilon$ may be approximated with Fourier series for each $\epsilon > 0$.
For $k = (k_1, \dots, k_d) \in \Z^d$ we set $\norm{k}_\infty= \max_i \abs{k_i}$ and $\norm{k}_1 = \sum_{i} \abs{k_i}$. For each $n \in \Z^+$ define $\Pi_n: \mathcal{C}(\mathbb{T}^d, \C) \to \mathcal{C}(\mathbb{T}^d, \C)$ by
\begin{equation*}
  (\Pi_{n}f)(x) = \sum_{\substack{k \in \Z^d \\ \norm{k}_\infty \le n}} \hat{f}(k) e^{2 \pi i x \cdot k},
\end{equation*}
where $\hat{f}$ denotes the Fourier transform\footnote{Specifically, $\hat{f}(k) = \intf_{\mathbb{T}^d} f(x) e^{-2 \pi i x \cdot k} d\Leb(x)$. } of $f$. For every $\epsilon \in (0, \epsilon_1)$ and $n \in \Z^+$ let $\LL_{\epsilon, n} = \Pi_n \LL_\epsilon$.
To simplify our notation, we set $\LL_{\epsilon, \infty} = \LL_\epsilon$. Our main technical result for this section is the following.

\begin{proposition}\label{prop:convolution_fourier_approx}
  There exists $\epsilon_3 \in (0, \epsilon_2)$ and a map $N: [0, \epsilon_3) \to \N \cup \{\infty\}$ with $N^{-1}(\infty) = \{0\}$, so that for any map $n: [0, \epsilon_3) \to \N \cup \{\infty\}$ with $n \ge N$ the family of operators $\{ \LL_{\epsilon, n(\epsilon)}\}_{\epsilon \in [0, \epsilon_3)}$ satisfies (KL) on $B^{p,q}_\C$ with $\wnorm{\cdot} = \norm{\cdot}_{p-1,q+1}$.
\end{proposition}

\begin{remark}\label{remark:peripheral_eigenvalues}
  By Proposition \ref{prop:convolution_fourier_approx} we may apply the results in \cite{keller1999stability} to $\{\LL_{\epsilon,n(\epsilon)}\}_{\epsilon \in [0, \epsilon_3)}$.
  Hence, all the isolated eigenvalues of $\LL$ with modulus strictly greater than the constant $\alpha$ in \eqref{eq:kl_ly} are approximated by eigenvalues of $\LL_{\epsilon,n(\epsilon)}$, with error vanishing as $\epsilon \to 0$.
  When such an eigenvalue of $\LL$ is simple, we additionally have that the corresponding eigenprojection and eigenvector are approximated by those of $\LL_{\epsilon,n(\epsilon)}$ in $\tnorm{\cdot}$ and $\wnorm{\cdot}$, respectively.
\end{remark}

Propositions \ref{prop:convolution_fourier_approx} and \ref{prop:anosov_hypothesis} allow us to apply Theorem \ref{thm:stability_of_statistical_data} to obtain the stability of the invariant measure, variance and rate function.
We note that since Anosov diffeomorphisms on tori are mixing \cite[Proposition 18.6.5]{katok1997introduction}, it follows that $1$ is the only eigenvalue of $\LL$ of modulus 1, and is a simple eigenvalue. In particular $T$ has a unique SRB measure $\mu \in B^{p,q}$, as required by Theorem \ref{thm:stability_of_statistical_data}.

\begin{theorem}\label{thm:twisted_anosov_numeric_stability}
  Suppose that $g \in \mathcal{C}^{r}(\mathbb{T}^d, \R)$ satisfies $\intf g d\mu = 0$ and is not a $L^2(\mu)$-coboundary. Denote by $M_g(\cdot) : \C \to L(B^{p,q}_\C)$ the map $M_g(z)f = e^{zg} f$.
  Let $N, n$ satisfy the conditions of Proposition \ref{prop:convolution_fourier_approx}.
  There exists $\theta, \epsilon' > 0$ so that for each $\epsilon \in [0, \epsilon')$ and $z \in D_\theta$ the operator $\LL_{\epsilon,n(\epsilon)}(z)$ is quasi-compact and simple with leading eigenvalue $\lambda_{\epsilon}(z)$ depending analytically on $z$.
  Moreover, we have stability of the following statistical data associated to $T$ and $\{ g \circ T^k\}_{k \in \N}$:
  \begin{enumerate}
    \item The invariant measure is stable: there exists eigenvectors $v_{\epsilon} \in B^{p,q}_\C$ of $\LL_{\epsilon,n(\epsilon)}$ for the eigenvalue $\lambda_{\epsilon}(0)$ for which $\lim_{\epsilon \to 0} \norm{v_{\epsilon} - \mu}_{p-1,q+1} = 0$.
    \item The variance is stable: $\lim_{\epsilon \to 0} \lambda_{\epsilon}^{(2)}(0) = \sigma^2$.
    \item The rate function is stable: For each sufficiently small compact subset $U$ of the domain of the rate function $r$ there exists an interval $V \subseteq (-\theta, \theta)$ so that
    \begin{equation*}
      \lim_{\epsilon \to 0} \sup_{z \in V} (sz -\log \abs{\lambda_\epsilon(z)}) = r(s)
    \end{equation*}
    uniformly on $U$.
  \end{enumerate}
\end{theorem}

\begin{remark}\label{remark:omega_v_leb}
  In Section \ref{sec:est_stat_prop_anosov} we aim to estimate the statistical properties of an Anosov map $T$ using
  Lemma \ref{lemma:convolution_kl}, Proposition \ref{prop:convolution_fourier_approx} and Theorem \ref{thm:twisted_anosov_numeric_stability}.
  However, these results concern the stability of the spectral data of the transfer operator associated to an Anosov map $T$ on the $d$-dimensional torus $\mathbb{T}^d$ \emph{equipped with an adapted metric}.
  In particular, this operator, say $\LL_\Omega$, is defined by duality with respect to the adapted Riemannian measure $\Omega$.
  From a computational perspective one would much rather approximate the transfer operator $\LL_{\Leb}$ that is defined by duality with respect to $\Leb$, the Haar probability measure on $\mathbb{T}^d$, since this removes the need to compute any quantities that depend on the adapted metric.
  Luckily, the relationship between these operators (and their twists) is simple: they are conjugate and therefore have the same spectrum (see Proposition \ref{prop:conjugacy}).
  Hence an approximation of the spectrum of $\LL_{\Leb}(z)$ is also an approximation of the spectrum of $\LL_{\Omega}(z)$.
  However, it is not clear from the proofs in this section that if $\{\Pi_{n(\epsilon)}Q_\epsilon \LL_{\Omega}\}_{\epsilon \in [0,\epsilon')}$ satisfies (KL) due to Proposition \ref{prop:convolution_fourier_approx} then so too does $\{\Pi_{n(\epsilon)}Q_\epsilon \LL_{\Leb}\}_{\epsilon \in [0,\epsilon')}$ i.e. a numerical scheme that is valid for $\LL_{\Omega}$ may not be valid for $\LL_{\Leb}$.
  In Appendix \ref{app:stability_for_leb_pf} we show that this obstruction does not occur, at least not in the current setting; the relevant results are Propositions \ref{prop:conjugacy}, \ref{prop:Omega_to_Leb_KL} and \ref{prop:mollifier_scheme_leb}.
\end{remark}

The remainder of this section is dedicated to the proofs of the aforementioned results.

\begin{proof}[{The proof of Lemma \ref{lemma:compact_convolution}}]
  Let $\epsilon \in (0, \epsilon_1)$, $k \in \N^d$ and $f \in \mathcal{C}^r(\mathbb{T}^d,\C)$. Denote $\frac{\partial^{\norm{k}_1}}{\partial x_1^{k_1} \dots \partial x_d^{k_d}}$ by $\partial_k$.
  From the beginning of \cite[Section 4]{gouezel2006banach} the map $h \mapsto \intf h d\Omega$ is bounded on $B^{p,q}_\C$. Moreover, multiplication by $\mathcal{C}^r$ functions is bounded on $B^{p,q}_\C$ by \cite[Lemma 3.2]{gouezel2006banach}.
  These facts, together with standard properties of convolutions, imply that there exists a $C$ independent of $f, k$ and $\epsilon$ such that
  \begin{equation*}\begin{split}
    \sup_{x \in \mathbb{T}^d} \abs{(\partial_k(q_\epsilon * f))(x)} &= \sup_{x \in \mathbb{T}^d} \abs{\intf (\partial_k q_\epsilon)(x-y )  f(y) d\Leb(y)} \\
    &= \sup_{x \in \mathbb{T}^d} \abs{\intf (\partial_k q_\epsilon)(x-y ) f(y) \frac{\mathrm{d}\Leb}{\mathrm{d}\Omega}(y) d\Omega(y)}
     \le C \norm{f \frac{\mathrm{d}\Leb}{\mathrm{d}\Omega}}_{p,q} \norm{\partial_k q_\epsilon}_{\mathcal{C}^q}.
  \end{split}\end{equation*}
  Since $\frac{\mathrm{d}\Leb}{\mathrm{d}\Omega} \in \mathcal{C}^\infty(\mathbb{T}^d, \R)$, using the continuity of multiplication by $\mathcal{C}^r$ functions again yields
  \begin{equation*}
    \norm{ Q_\epsilon f}_{\mathcal{C}^{\norm{k}_1}} \le C' \norm{ \frac{\mathrm{d}\Leb}{\mathrm{d}\Omega}}_{\mathcal{C}^{r}}\norm{f}_{p,q} \norm{q_\epsilon}_{\mathcal{C}^{q + \norm{k}_1}}
  \end{equation*}
  for some appropriate constant $C'$. Hence, $Q_\epsilon$ extends to a bounded operator $Q_\epsilon : B^{p,q}_\C \to \mathcal{C}^{m} (\mathbb{T}^d, \C)$ for every $m \in \Z^+$.
  It follows that $Q_\epsilon$ also extends to a bounded operator $Q_\epsilon: B^{p,q}_\C \to \mathcal{C}^{\infty}(\mathbb{T}^d, \C)$.
  As bounded sets in $\mathcal{C}^{\infty}(\mathbb{T}^d, \C)$ are compact in $\mathcal{C}^{m} (\mathbb{T}^d, \C)$ for every $m \in \Z^+$, each operator $Q_\epsilon: B^{p,q}_\C \to \mathcal{C}^{m} (\mathbb{T}^d, \C)$ is therefore compact.
  That $Q_\epsilon:B^{p,q}_\C \to B^{p,q}_\C$ is compact follows from the continuous embedding of $\mathcal{C}^r(\mathbb{T}^d, \C)$ into $B^{p,q}_\C$ \cite[Remark 4.3]{gouezel2006banach}.
  It is standard that $Q_\epsilon \in L(\mathcal{C}^{m_1}(\mathbb{T}^d,\C), \mathcal{C}^{m_2}(\mathbb{T}^d,\C))$ for each $m_1,m_2 \in \Z^+$.
\end{proof}

\begin{proof}[{The proof of Lemma \ref{lemma:convolution_kl}}]
  For each $y \in \mathbb{T}^d$ let $T_y : \mathbb{T}^d \to \mathbb{T}^d$ be defined by $T_y(x) = T(x) + y$, and let $\LL_{T_y}$ denote the transfer operator associated with $T_y$ (defined by duality as in \eqref{eq:pf_definition}).
  Let $h \in \mathcal{C}^r(\mathbb{T}^d,\C)$ and $x,y \in \mathbb{T}^d$.
  Let $\tau_y$ denote the translation map induced by $y$. As $D_x T_y = (D_{T(x)}\tau_y )(D_x T)$ we have
  \begin{equation*}\begin{split}
    (\LL h)(x-y) = (h \abs{\det D T}^{-1} ) \circ T^{-1}(x-y) &= h\circ T_y^{-1}(x) \cdot \abs{\det D_{T_y^{-1}(x)} T}^{-1} \\
    &= h\circ T_y^{-1}(x) \cdot \abs{\det D_{T_y^{-1}(x)} T_y}^{-1} \abs{\det (D_{x-y}\tau_y)^{-1}}^{-1} \\
    &=(\LL_{T_y} h)(x) \abs{\det D_{x}\tau_{-y}}.
  \end{split}\end{equation*}
  Hence,
  \begin{equation}\label{eq:convolution_kl_1}
    (\LL_\epsilon h)(x) = \intf q_\epsilon(y) (\LL h)(x-y)d\Leb(y) = \intf q_\epsilon(y) \abs{\det D_{x}\tau_{-y}} (\LL_{T_y}h)(x) d\Leb(y).
  \end{equation}
  For $\epsilon \in (0, \epsilon_1)$ let $A_{\epsilon}$ be defined by
  \begin{equation*}
    (A_{\epsilon}h)(x) = \intf q_\epsilon(y) (\LL_{T_y}h)(x) d\Leb(y).
  \end{equation*}
  Set $A_0 = \LL$ and, for each $\epsilon \in [0, \epsilon_1)$, let $F_{\epsilon} = \LL_{\epsilon} - A_\epsilon$.
  We aim to apply Proposition \ref{prop:kl_norm_perturb} to $\{A_{\epsilon}\}_{\epsilon \in [0,\epsilon_1)}$ and $\{F_\epsilon\}_{\epsilon \in [0, \epsilon_1)}$, which would imply the required result as $\LL_\epsilon = A_\epsilon + F_\epsilon$.

  We begin by proving that $\{A_\epsilon\}_{\epsilon \in [0,\epsilon_1)}$ satisfies (KL).
  Note that $\{A_\epsilon\}_{\epsilon \in [0,\epsilon_1)}$ is a perturbation of the kind considered in \cite{gouezel2006banach}. Specifically, we take $\mathbb{T}^d$ to be their $\Omega$, the Haar measure $\Leb$ to be their $\mu$, and set $g(\omega, x) = q_\epsilon(\omega)$.
  Therefore, by the discussion between Corollary 2.6 and Theorem 2.7 in \cite{gouezel2006banach}, there exists some $\epsilon' \in (0, \epsilon_1)$ such that $\{A_\epsilon \}_{\epsilon \in [0, \epsilon')}$ satisfies (KL) on $B^{p,q}_\C$ with $\wnorm{\cdot} = \norm{\cdot}_{p-1,q+1}$ provided that
  \begin{enumerate}[label=(C\arabic*)]
    \item \label{en:C1} For a fixed, small, open (in the $\mathcal{C}^{r+1}(\mathbb{T}^d,\mathbb{T}^d)$ topology) neighbourhood $U$ of $T$ we have $T_y \in U$ whenever $y \in B_{\epsilon'}(0)$; and
    \item \label{en:C2} $\lim_{\epsilon \to 0} \intf q_\epsilon(y) {d}_{\mathcal{C}^{r+1}}(T_y, T) d\Leb(y) = 0$.
  \end{enumerate}
  The condition \ref{en:C2} is derived from \cite[equation (2.5)]{gouezel2006banach} by setting $g(\omega, x) = q_\epsilon(\omega)$, and observing that $x \mapsto q_\epsilon(\omega)$ is constant and so has $\mathcal{C}^{p+q}$ norm $\abs{q_\epsilon(\omega)} = q_\epsilon(\omega)$.
  The other term in \cite[equation (2.5)]{gouezel2006banach} is 0 since $\intf q_\epsilon(\omega) d\Leb(\omega) = 1$.
  As $\mathbb{T}^d$ is compact and $T \in \mathcal{C}^{r+1}(\mathbb{T}^d,\mathbb{T}^d)$, the map $x \mapsto D_x^k T$ is uniformly continuous for each $0 \le k \le r+1$.
  It then follows from the definition of $T_y$ that
  \begin{equation*}
    \lim_{\epsilon \to 0} \sup_{y \in B_\epsilon(0)} d_{\mathcal{C}^{r+1}}(T_y, T) = 0.
  \end{equation*}
  Recalling that $q_\epsilon$ satisfies \ref{en:S1} and \ref{en:S2}, it is clear that there exists $\epsilon' \in (0, \epsilon_1)$ so that $\{A_{\epsilon}\}_{\epsilon \in [0,\epsilon')}$ satisfies both \ref{en:C1} and \ref{en:C2}, and therefore also (KL).

  We will now prove that $\{F_{\epsilon}\}_{\epsilon \in [0, \epsilon')}$ satisfies the requirements of Proposition \ref{prop:kl_norm_perturb}.
  For $y \in \mathbb{T}^d$ let $f_y: \mathbb{T}^d \to \R$ be defined by $f_y(x) = 1 - \det D_{x}\tau_{-y}$.
  From the definition of $\norm{\cdot}_{p,q}$ and \ref{en:S1} we obtain
  \begin{equation*}
    \norm{F_\epsilon h}_{p,q} \le \intf q_{\epsilon}(y) \norm{f_y \cdot (\LL_{T_y} h)}_{p,q} d\Leb(y).
  \end{equation*}
  As multiplication by $\mathcal{C}^r$ functions is continuous on $B^{p,q}_\C$ (\cite[Lemma 3.2]{gouezel2006banach}), using \ref{en:S1} there is some $C >0$ such that
  \begin{equation}\label{eq:convolution_kl_2}
    \norm{F_\epsilon }_{p,q} \le C \sup_{y \in \supp q_\epsilon} \norm{\LL_{T_y} }_{p,q} \intf q_\epsilon(y) \norm{f_y}_{\mathcal{C}^r} d\Leb(y).
  \end{equation}
  As mentioned at the beginning of \cite[Section 7]{gouezel2006banach}, the estimates in \cite[Lemma 2.2]{gouezel2006banach} apply uniformly to every map in $U$ and so there exists some $\eta > 0$ such that $\sup_{y \in B_\eta(0)} \norm{\LL_{T_y}}_{p,q} < \infty$.
  Since $\det D_x\tau_y = \frac{\mathrm{d}\Omega}{\mathrm{d}\Leb}(x + y) \frac{\mathrm{d}\Leb}{\mathrm{d}\Omega}(x)$, and  $\frac{\mathrm{d}\Omega}{\mathrm{d}\Leb},  \frac{\mathrm{d}\Leb}{\mathrm{d}\Omega}  \in \mathcal{C}^\infty(\mathbb{T}^d, \R)$, we have $\lim_{y \to 0} f_y = 0$ in $\mathcal{C}^r(\mathbb{T}^d, \R)$.
  Applying these facts and \ref{en:S1} to \eqref{eq:convolution_kl_2} yields
  \begin{equation}\label{eq:convolution_kl_3}
    \lim_{\epsilon \to 0} \norm{F_\epsilon}_{p,q} \le C \limsup_{\epsilon \to 0} \sup_{y \in B_\epsilon(y)}\norm{f_y}_{\mathcal{C}^r} \sup_{y \in B_\epsilon(0)} \norm{\LL_{T_y}}_{p,q} = 0.
  \end{equation}
  The same argument applies when estimating $\norm{F_\epsilon }_{p-1,q+1}$, and so there exists $\epsilon'' \in (0, \epsilon')$ so that $\{F_\epsilon\}_{\epsilon \in [0, \epsilon'')}$ satisfies the requirements for Proposition \ref{prop:kl_norm_perturb}.

  Hence Proposition \ref{prop:kl_norm_perturb} applies to $\{A_{\epsilon}\}_{\epsilon \in [0,\epsilon'')}$ and $\{F_\epsilon\}_{\epsilon \in [0, \epsilon'')}$.
  Namely, there exists $\epsilon_2 \in (0, \epsilon'')$ so that $\{A_{\epsilon} + F_\epsilon\}_{\epsilon \in [0,\epsilon_2)}$ satisfies (KL). Since $A_{\epsilon} + F_\epsilon = \LL_\epsilon$, this completes the proof.
\end{proof}

We require the following classical result on the convergence of Fourier series on $\mathbb{T}^d$ (see e.g. \cite[Proposition 5.6 and the proof of Theorem 5.7]{roe1988elliptic}).
\begin{proposition}\label{prop:convergence_fourier}
  For each $m \in \N$ we have $\Pi_n \to \Id$ strongly in $L(\mathcal{C}^{m + \ceil*{\frac{d+1}{2}}}, \mathcal{C}^m)$.
\end{proposition}

\begin{proof}[{The proof of Proposition \ref{prop:convolution_fourier_approx}}]
  By Lemma \ref{lemma:convolution_kl} the family of operators $\{ \LL_\epsilon\}_{\epsilon \in [0, \epsilon_2)}$ satisfies (KL) on $B^{p,q}_\C$ with $\wnorm{\cdot} = \norm{\cdot}_{p-1,q+1}$.
  We plan to find $N: (0, \epsilon_2) \to \N$ so that we may apply Proposition \ref{prop:kl_norm_perturb} with $A_{\epsilon} = \LL_\epsilon$ and $B_\epsilon = \LL_{\epsilon, N(\epsilon)} - \LL_\epsilon$.

  By Proposition \ref{prop:convergence_fourier}, $\Pi_n \to \Id$ strongly in $L(\mathcal{C}^{r + \ceil*{\frac{d+1}{2}}}, \mathcal{C}^{r})$.
  As the unit ball of $\mathcal{C}^{r + 1 + \ceil*{\frac{d+1}{2}}}$ is compact in $\mathcal{C}^{r + \ceil*{\frac{d+1}{2}}}$, Proposition \ref{prop:convergence_fourier}, the uniform boundedness principle and standard estimates imply that $\Pi_n \to \Id$ in $L(\mathcal{C}^{r + 1+ \ceil*{\frac{d+1}{2}}}, \mathcal{C}^{r })$.
  As $\mathcal{C}^{r}$ embeds continuously into $B^{p,q}_\C$ \cite[Remark 4.3]{gouezel2006banach}, there exists $C > 0$ so that, for each $\epsilon \in [0, \epsilon_2)$ and $n \in \N$, we have
  \begin{equation*}
    \norm{\LL_{\epsilon,n} - \LL_\epsilon}_{p,q} \le C\norm{ \Pi_n - \Id}_{ L\left(\mathcal{C}^{r + 1 + \ceil*{\frac{d+1}{2}}}, \mathcal{C}^{r}\right)} \norm{Q_\epsilon}_{L\left(B^{p,q}_\C, \mathcal{C}^{r + 1 + \ceil*{\frac{d+1}{2}}}\right)} \norm{\LL}_{p,q}.
  \end{equation*}
  Hence, as $\norm{Q_\epsilon}_{L\left(B^{p,q}_\C, \mathcal{C}^{r + 1 + \ceil*{\frac{d+1}{2}}}\right)}$ is finite by Lemma \ref{lemma:compact_convolution}, for each $\epsilon \in (0, \epsilon_2)$ there exists $N_1(\epsilon)$ so that $\norm{\LL_{\epsilon,n} - \LL_\epsilon}_{p,q} \le \epsilon$ whenever $n > N_1(\epsilon)$.
  The same argument produces for each $\epsilon \in (0, \epsilon_2)$ an $N_2(\epsilon)$ so that $\norm{\LL_{\epsilon,n} - \LL_\epsilon}_{p-1,q+1} \le \epsilon$ whenever $n \ge N_2(\epsilon)$.
  To summarise, if $N(\epsilon) : = \max\{ N_1(\epsilon), N_2(\epsilon)\}$ and $n : (0, \epsilon_2) \to \N$ is such that $n \ge N$, then $\lim_{\epsilon \to 0} \norm{\LL_{\epsilon,n(\epsilon)} - \LL_\epsilon}_{p,q} = 0$ and
  \begin{equation}\label{eq:convolution_fourier_approx_1}
    \sup_{\epsilon \in (0, \epsilon_2)} \norm{\LL_{\epsilon,n(\epsilon)} - \LL_\epsilon}_{p-1,q+1} < \infty.
  \end{equation}
  Hence for each map $n \ge N$ we may apply Proposition \ref{prop:kl_norm_perturb} as planned to produce an $\epsilon_n \in (0, \epsilon_2)$ so that $\{ \LL_{\epsilon,n(\epsilon)} \}_{\epsilon \in [0, \epsilon_n)}$ satisfies (KL) on $B^{p,q}_\C$ with $\wnorm{\cdot} = \norm{\cdot}_{p-1,q+1}$.
  Examining the proof of Proposition \ref{prop:kl_norm_perturb}, we observe that $\epsilon_n$ may be chosen independently of $n$ since $\lim_{\epsilon \to 0} \sup_{\ell \ge N(\epsilon)} \norm{\LL_{\epsilon,\ell} - \LL_\epsilon}_{p,q} = 0$ and
  \begin{equation*}
    \sup_{\epsilon \in (0, \epsilon_2)} \sup_{\ell \ge N(\epsilon)} \norm{\LL_{\epsilon,\ell} - \LL_\epsilon}_{p-1,q+1} < \infty.
  \end{equation*}
\end{proof}

\begin{proof}[{The proof of Theorem \ref{thm:twisted_anosov_numeric_stability}}]
  Hence, by Proposition \ref{prop:anosov_hypothesis}, the sequence $\{g \circ T^k\}_{k \in \N}$ satisfies Hypothesis \ref{hypothesis:analytic_coding} with operator-valued map $z \mapsto \LL M_g(z)$ and Banach space $B^{p,q}_\C$.
  The closed, unit ball in $B^{p,q}_\C$ is relatively compact with respect to $\norm{\cdot}_{p-1,q+1}$ by \cite[Lemma 2.1]{gouezel2006banach} and $M_g$ is a compactly $\norm{\cdot}_{p-1,q+1}$-bounded twist by Proposition \ref{prop:defn_of_twist_anosov}.
  Since $\{\LL_{\epsilon, n(\epsilon)}\}_{\epsilon \in [0, \epsilon_3)}$ satisfies (KL), we have verified all the requirements of Theorem \ref{thm:stability_of_statistical_data}, and so all the claims in the statement of Theorem \ref{thm:twisted_anosov_numeric_stability} follow, with the exception of the stability of the invariant measure.
  This claim follows from \cite[Proposition 2.4, Remark 2.5]{froyland2014detecting}, whose hypotheses are verified due to the convergence of eigenprojections in Theorem \ref{theorem:convergence_of_eigendata_derivs}.
\end{proof}

\section{Statistical stability for some Anosov maps under stochastic perturbations}
\label{sec:fejer_anosov}

In this section we show that the statistical properties of some Anosov maps on the $d$-dimensional torus $\mathbb{T}^d$ may be approximated using (weighted) Fourier series by realising the Fej{\'e}r kernel as a stochastic perturbation.
More generally, our results expand upon the stability results in \cite[Theorem 2.7]{gouezel2006banach} by allowing the stochastic kernel to be supported on all of $\mathbb{T}^d$ at the cost of additional requirements on the dynamics.
As a consequence of this generalisation, we \emph{no longer require a mollifier} to estimate the spectral and statistical properties of Anosov maps as in Section \ref{sec:approximation_anosov}, improving the computational aspect of our theory.

We adopt the setting, assumptions and notation from Section \ref{sec:func_analytic_anosov} and fix $p \in \Z^+$ and $q > 0$ satisfying $p + q < r$.
Our main assumption on the dynamics is that the associated transfer operator $\LL$ satisfies
\begin{equation}\label{eq:translate_bound}
  \sup_{y \in \mathbb{T}^d} \max\{ \norm{\tau_y\LL}_{p,q}, \norm{\tau_y\LL}_{p-1,q+1} \} < \infty,
\end{equation}
where $\tau_y$ denotes the translation operator induced by $y \in \mathbb{T}^d$.
In Section \ref{sec:cond_stochastic_stab_anosov} we provide some conditions for a map to satisfy \eqref{eq:translate_bound}.
For example, Proposition \ref{prop:translate_bound_cond} implies that \eqref{eq:translate_bound} holds for an iterate of $T$ if $T$ is close to a hyperbolic linear toral automorphism. Note that throughout this section we treat $\tau_y$ as both a composition operator and map.

For some $\epsilon_1 > 0$, suppose that $\{q_\epsilon  \}_{\epsilon \in (0, \epsilon_1)} \subseteq L^1(\Leb)$ is a family of stochastic kernels satisfying \ref{en:S1} and \begin{enumerate}[label=(S3)]
  \item \label{en:S3} For every $\eta > 0$, we have $\lim_{\epsilon \to 0} \intf_{\mathbb{T}^d \setminus B_\eta(0)} q_\epsilon d\Leb = 0$.
\end{enumerate}
\ref{en:S3} replaces \ref{en:S2} from Section \ref{sec:approximation_anosov} and allows the support of each $q_\epsilon$ to be all of $\mathbb{T}^d$.
Also note that we place no regularity requirements on the kernels $q_\epsilon$.
For $\epsilon > 0$ define $\LL_\epsilon := q_\epsilon * \LL$, and let $\LL_0 := \LL$. Our main technical result for this section is the following.

\begin{proposition}\label{prop:stochastic_perturbation_kl}
  If \eqref{eq:translate_bound} holds then there exists $\epsilon_2 \in (0 ,\epsilon_1)$ so that $\{\LL_\epsilon \}_{\epsilon \in [0, \epsilon_2)}$ satisfies (KL) on $B^{p,q}_\C$ with $\wnorm{\cdot} = \norm{\cdot}_{p-1,q+1}$.
\end{proposition}

\begin{remark}
  The same comments as those in Remark \ref{remark:peripheral_eigenvalues} apply to the family of operators $\{\LL_{\epsilon}\}_{\epsilon \in [0, \epsilon_2)}$.
  That is, the spectral data associated to the eigenvalues of $\LL_0$ with modulus greater than the constant $\alpha$ appearing in \eqref{eq:kl_ly} are well approximated by the spectral data of associated to eigenvalues of $\LL_{\epsilon}$, with error vanishing as $\epsilon \to 0$.
\end{remark}

As noted before the statement of Theorem \ref{thm:twisted_anosov_numeric_stability}, since $T$ is an Anosov diffeomorphism on a torus, $\LL$ is a simple quasi-compact operator. Thus, if \eqref{eq:translate_bound} holds then by Proposition \ref{prop:stochastic_perturbation_kl} the peripheral spectral data, invariant measure, variance and rate function are stable with respect to the class of stochastic perturbations in consideration.
The proof is the same as that of Theorem \ref{thm:twisted_anosov_numeric_stability}.

\begin{theorem}\label{thm:stochastic_perturb_stability}
  Assume that $\LL$ satisfies \eqref{eq:translate_bound}. Suppose that $g \in \mathcal{C}^{r}(\mathbb{T}^d, \R)$ satisfies $\intf g d\mu = 0$ and is not an $L^2(\mu)$-coboundary. Denote by $M_g : \C \to L(B^{p,q}_\C)$ the map defined by $M_g(z)f = e^{zg} f$.
  There exists $\theta, \epsilon' > 0$ so that for each $\epsilon \in [0, \epsilon')$ and $z \in D_\theta$ the operator $\LL_{\epsilon}(z)$ is quasi-compact and simple with leading eigenvalue $\lambda_{\epsilon}(z)$ depending analytically on $z$. Moreover, we have stability of the following statistical data associated to $T$ and $\{ g \circ T^k\}_{k \in \N}$:
  \begin{enumerate}
    \item The invariant measure is stable: there exists eigenvectors $v_{\epsilon} \in B^{p,q}_\C$ of $\LL_{\epsilon}$ for the eigenvalue $\lambda_{\epsilon}(0)$ for which $\lim_{\epsilon \to 0} \norm{v_{\epsilon} - \mu}_{p-1,q+1} = 0$.
    \item The variance is stable: $\lim_{\epsilon \to 0} \lambda_{\epsilon}^{(2)}(0) = \sigma^2$.
    \item The rate function is stable: For each sufficiently small compact subset $U$ of the domain of the rate function $r_g$ there exists an interval $V \subseteq (-\theta, \theta)$ so that
    \begin{equation*}
      \lim_{\epsilon \to 0} \sup_{z \in V} (sz -\log \abs{\lambda_\epsilon(z)}) = r(s)
    \end{equation*}
    uniformly on $U$.
  \end{enumerate}
\end{theorem}

A key application of the results in this section is the rigorous approximation of the spectral and statistical data of some Anosov maps using Fourier series.
We define the $n$th 1-dimensional Fej{\'e}r kernel $K_{n,1}$ on $\mathbb{T}^1 = \R / \Z$ by
\begin{equation*}
  K_{n,1}(t) = \sum_{k=-n}^n \left(1 - \frac{\abs{k}}{n+1} \right)e^{2 \pi i k t}.
\end{equation*}
The $n$th $d$-dimensional Fej{\'e}r kernel $K_{n,d}$ on $\mathbb{T}^d$ is defined by
\begin{equation*}
  K_{n,d}(t) = \prod_{i=1}^d K_{n,1}(t_i),
\end{equation*}
where $t_i$ is the $i$th component of $t \in \mathbb{T}^d$. It is well known that the 1-dimensional Fej{\'e}r kernels satisfy \ref{en:S1} and \ref{en:S3} as $n \to \infty$ (they are a \emph{summability kernel}; see \cite[Section 2.2 and 2.5]{katznelson2002introduction}).
It is routine to verify that the $d$-dimensional kernels consequently satisfy the same conditions, and so we may apply Proposition \ref{prop:stochastic_perturbation_kl} with $q_{1/n} = K_{n,d}$.
A straightforward computation yields
\begin{equation*}
  K_{n,d}(t) = \sum_{\substack{k \in \Z^d \\ \norm{k}_\infty \le n}} \prod_{i=1}^d \left(1 - \frac{\abs{k_i}}{n+1}\right) e^{2\pi i k \cdot t},
\end{equation*}
and, therefore, convolution with the Fej{\'e}r kernel may be represented using weighted Fourier series:
\begin{equation}\label{eq:fejer_fourier_series}
  (K_{n,d} * f)(x) = \sum_{\substack{k \in \Z^d \\ \norm{k}_\infty \le n}} \prod_{i=1}^d \left(1 - \frac{\abs{k_i}}{n+1}\right) \hat{f}(k) e^{2\pi i k \cdot x}.
\end{equation}

By the above considerations and Proposition \ref{prop:stochastic_perturbation_kl} and Theorem \ref{thm:stochastic_perturb_stability} we obtain the following stability result for stochastic perturbations induced by the  Fej{\'e}r kernel.

\begin{corollary}\label{cor:fejer_stability}
  Assume that $\LL$ satisfies \eqref{eq:translate_bound} and that $g \in \mathcal{C}^{r}(\mathbb{T}^d, \R)$ satisfies $\intf g d\mu = 0$ and is not an $L^2(\mu)$-coboundary.
  For $n \in \N$ let $\LL_{1/n} := K_{n,d} * \LL$, where $K_{n,d}$ is the $d$-dimensional Fej{\'e}r kernel, and let $\LL_0 = \LL$.
  There exists $N > 0$ so that the family of operators $\{ \LL_{1/n}\}_{n \ge N}$ satisfies (KL) on $B^{p,q}_\C$ with $\wnorm{\cdot} = \norm{\cdot}_{p-1,q+1}$.
  Consequently, we have stability of the invariant measure, variance, and rate function associated to $T$ and $\{ g \circ T^k\}_{k \in \N}$ as in Theorem \ref{thm:stochastic_perturb_stability}.
\end{corollary}

The operators $\LL_{1/n}$ are finite-dimensional and leave the span of $\{ e^{2 \pi i k \cdot x} : \norm{k}_\infty \le n \}$ invariant.
Therefore, we could compute all of the spectral data of $\LL_{1/n}$ via its matrix representation with respect to the basis $\{ e^{2 \pi i k \cdot x} : \norm{k}_\infty \le n \}$ and use Corollary \ref{cor:fejer_stability} to estimate the statistical properties of $T$.

\begin{remark}\label{remark:omega_v_leb_fejer}
  The comments made in Remark \ref{remark:omega_v_leb} also apply here: the stability results in Proposition \ref{prop:stochastic_perturbation_kl}, Theorem \ref{thm:stochastic_perturb_stability} and Corollary \ref{cor:fejer_stability} apply to the transfer operator $\LL_{\Omega}$ that is associated to $T$ via duality with respect to $\Omega$, rather than the transfer operator $\LL_{\Leb}$ that is associated to $T$ via duality with respect to $\Leb$.
  In Appendix \ref{app:stability_for_leb_pf} we show that these operators, and their twists, are conjugate, and therefore have the same spectrum, and that if Proposition \ref{prop:stochastic_perturbation_kl}, Theorem \ref{thm:stochastic_perturb_stability} and Corollary \ref{cor:fejer_stability} apply to $\LL_{\Omega}$ then they also hold true when $\LL_{\Omega}$ is replaced by $\LL_{\Leb}$.
\end{remark}

\begin{proof}[{Proof of Proposition \ref{prop:stochastic_perturbation_kl}}]
  For each $\eta \in (0,1/2)$ let $p_\eta$ be the characteristic function of $B_\eta(0)$.
  For $\epsilon \ge 0$ let $a_{\epsilon, \eta} = \intf q_\epsilon p_\eta d\Leb$, $A_{\epsilon,\eta} = a_{\epsilon, \eta}^{-1} (q_\epsilon p_\eta) * \LL$ and $B_{\epsilon, \eta} = (q_\epsilon (1- a_{\epsilon, \eta}^{-1} p_\eta)) * \LL$.
  Set $A_{0,\eta} = \LL$ and $B_{0, \eta} = 0$.
  Our goal is to apply Proposition \ref{prop:kl_norm_perturb} by showing that (i) for some $\eta > 0$ there exists $\epsilon' \in (0, \epsilon_1)$ so that $\{A_{\epsilon, \eta}\}_{\epsilon \in [0, \epsilon')}$ satisfies (KL), and then proving (ii) that $\{ B_{\epsilon, \eta}\}_{\epsilon \in [0, \epsilon')}$ satisfies the necessary requirements of Proposition \ref{prop:kl_norm_perturb}.
  As $\LL_\epsilon = A_{\epsilon,\eta} + B_{\epsilon, \eta}$, this yields the required statement.

  Since $\{a_{\epsilon, \eta}^{-1} q_\epsilon p_\eta \}_{\epsilon \in (0, \epsilon_1)}$ satisfies \ref{en:S1} the perturbation $\{A_{\epsilon, \eta}\}_{\epsilon \in [0, \epsilon_1)}$ is similar to the (convolution type) perturbation considered in Lemma \ref{lemma:convolution_kl}.
  As such we will explain how to modify the proof of Lemma \ref{lemma:convolution_kl} to obtain the required result.
  Examining the proof of Lemma \ref{lemma:convolution_kl}, we note that it was not important that the family of kernels was in  $\mathcal{C}^\infty$, indeed it is sufficient for the kernels to be contained in $L^1(\Leb)$.
  We must verify the conditions \ref{en:C1} and \ref{en:C2}, and a different argument is required here since $\{a_{\epsilon, \eta}^{-1} q_\epsilon p_\eta \}_{\epsilon \ge 0}$ does not satisfy \ref{en:S2}.
  By choosing $\eta$ sufficiently small we may make the support of every $a_{\epsilon, \eta}^{-1} q_\epsilon p_\eta$ small enough so that for every $\epsilon \in [0, \epsilon_1)$ and $y \in \supp q_\epsilon p_\eta$ the map $T_y(x) := T(x) + y$ is in the set $U$ from \ref{en:C1} in Lemma \ref{lemma:convolution_kl}.
  This verifies \ref{en:C1} in our setting; we will now verify \ref{en:C2}.
  By \ref{en:S1} and \ref{en:S3} for $\{q_\epsilon \}_{\epsilon \ge 0}$ we have $\lim_{\epsilon \to 0} a_{\epsilon, \eta'}^{-1} = 1$ for every $\eta' \in (0, \eta)$.
  Hence, for every $\eta' \in (0, \eta)$ we have
  \begin{equation}\begin{split}\label{eq:stochastic_perturbation_kl_1}
    \lim_{\epsilon \to 0} \int a_{\epsilon, \eta}^{-1}|q_\epsilon(y) &p_\eta(y)| {d}_{\mathcal{C}^{r+1}}(T_y, T) \mathrm{d}\Leb(y) = \lim_{\epsilon \to 0} \intf_{B_\eta(0)} q_\epsilon(y) {d}_{\mathcal{C}^{r+1}}(T_y, T) d\Leb(y)\\
    =& \lim_{\epsilon \to 0} \left(\intf_{B_\eta(0) \setminus B_{\eta'}(0)} q_\epsilon(y) {d}_{\mathcal{C}^{r+1}}(T_y, T) d\Leb(y) + \intf_{B_{\eta'}(0)} q_\epsilon(y) {d}_{\mathcal{C}^{r+1}}(T_y, T) d\Leb(y) \right)\\
    \le& \lim_{\epsilon \to 0} \left(\sup_{y \in B_\eta(0) \setminus B_{\eta'}(0)} {d}_{\mathcal{C}^{r+1}}(T_y, T)\right) \intf_{\mathbb{T}^d \setminus B_{\eta'}(0)} q_\epsilon(y) d\Leb(y) + \sup_{y \in B_{\eta'}(0)} {d}_{\mathcal{C}^{r+1}}(T_y, T) \\
    \le& \sup_{y \in B_{\eta'}(0)} {d}_{\mathcal{C}^{r+1}}(T_y, T),
  \end{split}\end{equation}
  where we have used both \ref{en:S1} and \ref{en:S3}.
  Using the fact that $\lim_{\eta' \to 0}\sup_{y \in B_{\eta'}(0)} {d}_{\mathcal{C}^{r+1}}(T_y, T) = 0$ we obtain \ref{en:C2} from \eqref{eq:stochastic_perturbation_kl_1}.
  It remains to provide an alternative proof of \eqref{eq:convolution_kl_3}. With $f_y = 1 - \det D_x \tau_{-y}$ as in the proof of Lemma \ref{lemma:convolution_kl}, using \ref{en:S3} we have for each $\eta' \in (0, \eta)$ that
  \begin{equation}\begin{split}\label{eq:stochastic_perturbation_kl_2}
    \lim_{\epsilon \to 0} \int a_{\epsilon, \eta}^{-1}q_\epsilon(y) p_\eta(y)& \norm{f_y}_{\mathcal{C}^r} \, \mathrm{d}\Leb(y) \\
    &\le \lim_{\epsilon \to 0} \left( \intf_{B_{\eta'}(0)} q_\epsilon(y) \norm{f_y}_{\mathcal{C}^r} d\Leb(y) + \intf_{B_\eta(0) \setminus B_{\eta'}(0)} q_\epsilon(y) \norm{f_y}_{\mathcal{C}^r} d\Leb(y) \right)\\
    &\le  \sup_{y \in B_{\eta'}(0)}\norm{f_y}_{\mathcal{C}^r} + \lim_{\epsilon \to 0} \left(\intf_{B_\eta(0) \setminus B_{\eta'}(0)} q_\epsilon(y) d\Leb(y)\right) \sup_{y \in B_\eta(0)} \norm{f_y}_{\mathcal{C}^r} \\
    &\le \sup_{y \in B_{\eta'}(0)}\norm{f_y}_{\mathcal{C}^r}.
  \end{split}\end{equation}
  Since $\lim_{y \to 0} \norm{f_y}_{\mathcal{C}^r} = 0$, by letting $\eta' \to 0$ in \eqref{eq:stochastic_perturbation_kl_2} we can conclude that the left-hand side of \eqref{eq:stochastic_perturbation_kl_2} is 0.
  Recalling that $T_y \in U$ for every $y \in \supp q_\epsilon p_\eta$, it follows that $\sup_{y \in \supp q_\epsilon p_\eta} \norm{\LL_{T_y}}_{p,q} < \infty$.
  Hence
  \begin{equation*}
    \lim_{\epsilon \to 0} \sup_{y \in \supp q_\epsilon p_\eta} \norm{\LL_{T_y}}_{p,q} \int a_{\epsilon, \eta}^{-1}q_\epsilon(y) p_\eta(y) \norm{f_y}_{\mathcal{C}^r} \, \mathrm{d}\Leb(y) = 0,
  \end{equation*}
  which proves \eqref{eq:convolution_kl_3} in our setting. The same comments made in the sentence following \eqref{eq:convolution_kl_3} apply here too.
  Hence the arguments in Lemma \ref{lemma:convolution_kl} apply to $\{A_{\epsilon, \eta}\}_{\epsilon \in [0, \epsilon_1)}$, and so
  $\{A_{\epsilon, \eta}\}_{\epsilon \in [0, \epsilon_1)}$ satisfies (KL) on $B^{p,q}_\C$ with $\wnorm{\cdot} = \norm{\cdot}_{p-1,q+1}$.

  We will prove that $\lim_{\epsilon \to 0} B_{\epsilon,\eta} = 0$ in both $L(B^{p,q})$ and $L(B^{p-1,q+1})$, as this readily implies the same for $L(B^{p,q}_\C)$ and $L(B^{p-1,q+1}_\C)$.
  Let $h \in \mathcal{C}^r(\mathbb{T}^d, \R)$, $k \le p$ be a non-negative integer, $W \in \Sigma$, $\{v_i\}_{i=1}^k \subseteq \mathcal{V}^r(W)$ with $\norm{v_i}_{\mathcal{C}^r} \le 1$, and $\varphi \in \mathcal{C}^{q + k}_0(W, \R)$ with $\norm{\varphi}_{\mathcal{C}^{q + k}} \le 1$.
  With $s_{\epsilon, \eta} = q_\epsilon (1 - a_{\epsilon,\eta}^{-1}p_\eta)$, we have
  \begin{equation*}\begin{split}
    \big\lvert\int_W v_1 \dots v_k ( s_{\epsilon, \eta} * \LL h)(x) &\cdot \varphi(x) \,\mathrm{d}\Omega(x) \big\rvert \\
    &\le \intf_{\mathbb{T}^d} \abs{s_{\epsilon, \eta}(y)} \abs{\intf_W v_1 \dots v_k (\tau_{-y} \LL h)(x) \cdot\varphi(x) d\Omega(x)} d\Leb(y)\\
    &\le \left(\sup_{y \in \mathbb{T}^d} \norm{\tau_{y} \LL}_{p,q}\right) \norm{h}_{p,q} \intf_{\mathbb{T}^d} \abs{s_{\epsilon, \eta}} d\Leb(y).
  \end{split}\end{equation*}
  Hence,
  \begin{equation}\begin{split}\label{eq:stochastic_perturbation_kl_3}
    \norm{B_{\epsilon,\eta}}_{p,q} &\le \left(\sup_{y \in \mathbb{T}^d} \norm{\tau_{y} \LL}_{p,q}\right) \intf_{\mathbb{T}^d} \abs{s_{\epsilon, \eta} } d\Leb(y) \\
    &\le \left(\sup_{y \in \mathbb{T}^d} \norm{\tau_{y} \LL}_{p,q}\right) \intf q_\epsilon(y)\abs{ 1- a_{\epsilon,\eta}^{-1}p_\eta(y)} d\Leb(y).
  \end{split}\end{equation}
  Since $\lim_{\epsilon \to 0} a_{\epsilon,\eta} = 1$, using \ref{en:S3} we have
  \begin{equation}\begin{split}\label{eq:stochastic_perturbation_kl_4}
    \lim_{\epsilon \to 0}\int q_\epsilon(y)&\abs{ 1- a_{\epsilon,\eta}^{-1}p_\eta(y)}\, \mathrm{d}\Leb(y) \\
    &\le \lim_{\epsilon \to 0}\left(\intf q_\epsilon(y)(1- p_\eta(y)) d\Leb(y) + \intf q_\epsilon(y)p_\eta \abs{1- a_{\epsilon,\eta}^{-1}} d\Leb(y) \right) \\
    &\le \lim_{\epsilon \to 0} \left(\intf_{\mathbb{T}^d \setminus B_\eta(0)} q_\epsilon(y) d\Leb(y) + \abs{1- a_{\epsilon,\eta}^{-1}}\right) = 0.
  \end{split}\end{equation}
  Together, \eqref{eq:translate_bound}, \eqref{eq:stochastic_perturbation_kl_3} and \eqref{eq:stochastic_perturbation_kl_4} imply that $\lim_{\epsilon \to 0} \norm{B_{\epsilon,\eta}}_{p,q} = 0$.
  The same argument proves that $\lim_{\epsilon \to 0} \norm{B_{\epsilon,\eta}}_{p-1,q+1} = 0$, and so there exists $\epsilon' \in (0, \epsilon_1)$ so that $\sup_{\epsilon \in (0, \epsilon')} \norm{B_{\epsilon,\eta}}_{p-1,q+1} < \infty$.
  We have verified the conditions of Proposition \ref{prop:kl_norm_perturb}, concluding the proof.
\end{proof}

\subsection{A class of maps satisfying Proposition \ref{prop:stochastic_perturbation_kl}, Theorem \ref{thm:stochastic_perturb_stability}, and Corollary \ref{cor:fejer_stability}}\label{sec:cond_stochastic_stab_anosov}

Our main result for this section, Proposition \ref{prop:translate_bound_cond}, gives conditions for $T$ to have an iterate satisfying \eqref{eq:translate_bound}. For instance, we will deduce that Proposition \ref{prop:translate_bound_cond} applies to Anosov maps that are sufficiently close to hyperbolic linear toral automorphisms.

In \cite[Section 3]{gouezel2006banach} the usual Euclidean metric on $\mathbb{T}^d$ is replaced by an equivalent \emph{adapted} metric.
The choice of adapted metric will be crucial to our arguments in this section, so we begin by reviewing the construction of such metrics, following \cite[Proposition 5.2.2]{brin2002introduction}. Let $\wnorm{\cdot}$ denote the usual Euclidean norm on the tangent space of $\mathbb{T}^d$.
Let $E^s(x)$ and $E^u(x)$ denote the stable and unstable directions, respectively, of $T$ at $x$, and let $\pi^s_x$ and $\pi_x^u$ be the projections induced by the splitting $T_x \mathbb{T}^d = E^s(x) \oplus E^u(x)$. Let $d_s = \dim E^s(x)$ and $d_u = \dim E^u(x)$.
Since $T$ is Anosov, there exists $C > 0$, $\lambda_s \in (0,1)$, and $\lambda_u > 1$ so that $\wnorm{\restr{D_x T^n}{E^s(x)}} \le C\lambda_s^n$ and $\wnorm{\restr{D_x T^{-n}}{E^u(x)}} \le C\lambda_u^{-n}$ for every $n \in \N$.
Recall $\nu_u$ and $\nu_s$ from Section \ref{sec:func_analytic_anosov}, and that $\lambda_s < \nu_s < 1<  \nu_u < \lambda_u$. Let $N$ be such that $\max\{\frac{\lambda_s^{N+1}}{\nu_s^{N+1}}, \frac{\nu_u^{N+1}}{\lambda_u^{N+1}}\} <1/C$.
For $v^s \in E^s(x)$ and $v^u \in E^u(x)$ we define
\begin{equation*}
  \norm{v^s}_0 = \sum_{k=0}^N \nu_s^{-k} \wnorm{D_x T^k v^s}, \text{ and } \norm{v^u}_0 = \sum_{k=0}^N \nu_u^{k} \wnorm{D_x T^{-k} v^u}.
\end{equation*}
Note that $\norm{\restr{D_x T}{E^s(x)}}_0 < \nu_s$ and $\norm{\restr{D_x T^{-1}}{E^u(x)}}_0 < \nu_u^{-1}$.
For $v \in T_x\mathbb{T}^d$ we define
\begin{equation*}
  \norm{v}_0 = \sqrt{\norm{\pi^s_x v}_0^2 + \norm{\pi^u_x v}_0^2},
\end{equation*}
and, since $\norm{\cdot}_0$ satisfies the parallelogram law, we may recover a metric $\langle \cdot ,\cdot \rangle_0$ via the polarisation identity.
Note that $E^s(x) \perp E^u(x)$ with respect to $\langle \cdot ,\cdot \rangle_0$.
Because $x \mapsto \langle \cdot ,\cdot \rangle_0$ is not necessarily smooth (and so $\mathbb{T}^d$ equipped with $\langle \cdot ,\cdot \rangle_0$ would not be a $\mathcal{C}^\infty$ Riemannian manifold), for each sufficiently small $\xi > 0$ we instead consider a smooth metric $\langle \cdot, \cdot \rangle_\xi$ (with corresponding norm denoted $\norm{\cdot}_\xi$) such that
\begin{enumerate}[label=(M\arabic*)]
\item \label{en:M1} $\sup_{x \in \mathbb{T}^d}
  \sup_{\substack{v,w \in T_x \mathbb{T}^d \\ \norm{v}_0, \norm{w}_0 \le 1}}
  \abs{\langle v, w \rangle_\xi - \langle v, w \rangle_0} < \xi$;
\item \label{en:M2}  $\norm{\restr{D_x T}{E^s(x)}}_\xi < \nu_s$ and $\norm{\restr{D_x T}{E^u(x)}}_\xi < \nu_u^{-1}$; and
\item \label{en:M3}  $E^{s}(x)$ and $E^{u}(x)$ are $\xi$-orthogonal: for $w_s \in E^s(x)$ and $w_u \in E^u(x)$ with $\norm{w_s}_\xi, \norm{w_u}_\xi \le 1$ we have $\abs{\langle w_s, w_u\rangle_\xi} < \xi$.
\end{enumerate}
A metric is called adapted if it satisfies \ref{en:M2}.
For sufficiently small $\xi$, metrics satisfying \ref{en:M1}-\ref{en:M3} can be constructed by approximating $\ip_0$ (see e.g. \cite{hirsch2012differential}).

Let $\xi \ge 0$.
For $x \in \mathbb{T}^d$ we denote by $\Gamma^s_x$ and $\Gamma^u_x$ the orthogonal (with respect to $\ip_\xi$) projections onto $E^{s}(x)$ and $E^{u}(x)$, respectively. Although $\Gamma^s_x$ clearly depends on $\xi$, we suppress this from our notation.
Define
\begin{equation}\label{eq:almost_flat_constants}
  C_{\tau,\xi} = \sup_{x,y \in \mathbb{T}^d} \norm{D_x \tau_y}_\xi \quad \text{and} \quad \Theta_{T,\xi} = \sup_{x,y \in \mathbb{T}^d} \norm{\Gamma^s_{x+y} D_x \tau_y - (D_x \tau_y) \Gamma^s_{x}}_\xi.
\end{equation}
Note that both $C_{\tau,\xi}$ and $\Theta_{T,\xi}$ are finite. The key hypothesis for this section's main result is that $C_{\tau,0}^{-1} > \Theta_{T,0}$.
Roughly speaking, this condition ensures that translated leaves never lie in the unstable direction of $T$ (recall that leaves are approximately parallel to the stable directions).
One way to see this is by computing the quantities $C_{\tau,0}$ and $\Theta_{T,0}$ when $\norm{\cdot}_0$ is the usual Euclidean norm (ignoring the issue of whether the Euclidean norm is adapted). In this case, $C_{\tau,0} = 1$ and $\Theta_{T,0}$ measures the angle between $E^s(x+y)$ and $E^s(x)$. If this angle is everywhere close to 0 then the translate of a leaf will be approximately parallel to the stable direction of $T$ regardless of the translate.

We will select a specific adapted metric to both streamline our arguments and strengthen our results; although in doing so we impact the definition of the set of leaves $\Sigma$, and hence also the spaces $B^{p,q}_\C$ and $B^{p-1,q+1}_\C$.
Our main technical result for this section is the following; its proof constitutes Appendix \ref{app:proof_of_long_prop}.

\begin{proposition}\label{prop:translate_bound_cond}
   If
   \begin{equation}\label{eq:translate_bound_cond_0}
     C_{\tau,0}^{-1} > \Theta_{T,0},
   \end{equation}
  then there exists $N \in \Z^+$, an adapted metric $\langle \cdot, \cdot \rangle$, and a set of leaves $\tilde{\Sigma}$, inducing spaces $\tilde{B}^{p,q}_\C$ and $\tilde{B}^{p-1,q+1}_\C$, so that $\LL$ is quasi-compact on $\tilde{B}^{p,q}_\C$, with the same spectral data associated to eigenvalues outside of the ball of radius $\max\{\nu_u^{-p}, \nu_s^{q} \}$ as when considered as an operator on $B^{p,q}_\C$, and so that
   \begin{equation}\label{eq:translate_bound_cond_1}
     \sup_{y \in \mathbb{T}^d} \max\left\{ \norm{\tau_y\LL^N}_{p,q}, \norm{\tau_y \LL^N}_{p-1,q+1} \right\} < \infty.
   \end{equation}
\end{proposition}

We make two comments regarding the applicability of Proposition \ref{prop:translate_bound_cond}.
Firstly, maps satisfying \eqref{eq:translate_bound_cond_0} exist as $\Theta_{T,0} = 0$ whenever $T$ is a linear hyperbolic toral automorphism.
Secondly, the condition \eqref{eq:translate_bound_cond_0} is open in $\mathcal{C}^{r+1}(\mathbb{T}^d,\mathbb{T}^d)$.
To see this, suppose that $T$ satisfies \eqref{eq:translate_bound_cond_0}, and let $\ip$ be the metric one obtains by applying Proposition \ref{prop:translate_bound_cond} to $T$.
The following comments apply to all $T'$ in a sufficiently small $\mathcal{C}^{r+1}$-neighbourhood of $T$: $T'$ is Anosov map, $\ip$ is an adapted metric for $T'$, and the stable and unstable directions for $T$ and $T'$ are everywhere close in the Grassmanian.
It follows that $T'$ also satisfies \eqref{eq:translate_bound_cond_0} provided that it is sufficiently close to $T$ in the $\mathcal{C}^{r+1}$ topology.

We now give a concrete example of maps satisfying the conditions of Proposition \ref{prop:translate_bound_cond}.
For $\delta \in \R$ let $T_\delta : \mathbb{T}^2 \to \mathbb{T}^2$ be defined by
\begin{equation*}
  T_\delta(x_1, x_2) = (2x_1 + x_2, x_1 + x_2) + \delta(\cos(2\pi x_1), \sin(4 \pi x_2 +1)).
\end{equation*}
Note that $T_0$ is Arnold's `cat map' - a linear hyperbolic toral automorphism - and so it satisfies Proposition \ref{prop:translate_bound_cond}.
Moreover, since $\delta \mapsto T_\delta$ is smooth, it follows from the discussion in the previous paragraph that $T_\delta$ is Anosov and satisfies the conditions of Proposition \ref{prop:translate_bound_cond} for sufficiently small $\delta$.
In Appendix \ref{sec:perturbed_cat} we find an explicit range of $\delta$ satisfying these conditions, yielding the following result.

\begin{proposition}\label{thm:perturbed_cat_props}
  If $0\le \delta < 0.0108$ then $T_\delta$ is an Anosov diffeomorphism and satisfies the conditions of Proposition \ref{prop:translate_bound_cond}.
\end{proposition}

\section{Estimation of the statistical properties of Anosov maps}\label{sec:est_stat_prop_anosov}

In this section we implement the numerical schemes described in Sections \ref{sec:approximation_anosov} and \ref{sec:fejer_anosov}.
These are, respectively:

\paragraph{Fourier approximation of mollified transfer operators.}
Proposition \ref{prop:convolution_fourier_approx} says that if we convolve the (possibly twisted) transfer operator $\LL_{\Omega}$ with a locally supported stochastic kernel (parameterised by $\epsilon$), Fourier approximations (of order $n=n(\epsilon)$) of this mollified transfer operator satisfy the $(\|\cdot\|_{p,q},\|\cdot\|_{p-1,q+1})$ (KL) conditions as a family of operators in $\epsilon$.
That the same holds for $\LL_{\Leb}$ follows from Proposition \ref{prop:mollifier_scheme_leb}.
The Fourier approximations of $\LL_{\Leb}$ are numerically accessible and
Theorem \ref{thm:twisted_anosov_numeric_stability} and Proposition \ref{prop:mollifier_scheme_leb} then guarantees convergence of the SRB measure (in the $\|\cdot\|_{p-1,q+1}$ norm), convergence of the variance of a $\mathcal{C}^r$ observable, and uniform convergence of the rate function for $\mathcal{C}^r$ observables, as $\epsilon \to 0$.

\paragraph{Direct Fourier approximation via Fej{\'e}r kernels.}
Corollary \ref{cor:fejer_stability} states that if we convolve the (possibly twisted) transfer operator $\LL_\Omega$ with a Fej{\'e}r kernel (parameterised by $n$), this sequence of operators in $n$ satisfy the $(\|\cdot\|_{p,q},\|\cdot\|_{p-1,q+1})$ (KL) conditions.
That the same holds for $\LL_{\Leb}$ follows from Proposition \ref{prop:fejer_scheme_leb}.
The Fej{\'e}r kernels directly arise from Fourier projections and this second numerical scheme requires only direct Fourier approximation of the transfer operators $\LL_{\Leb}$.
Theorem \ref{thm:stochastic_perturb_stability} and Proposition \ref{prop:fejer_scheme_leb} guarantees convergence of the SRB measure (in the $\|\cdot\|_{p-1,q+1}$ norm), convergence of the variance of a $\mathcal{C}^r$ observable, and uniform convergence of the rate function for $\mathcal{C}^r$ observables, as $n\to\infty$.

For the remainder of this section we will only deal with the operator $\LL_{\Leb}$, it's twists $\LL_{\Leb}(z)$ and approximations of both these two operators i.e. $\LL_{\Leb,\epsilon}$ and $\LL_{\Leb,\epsilon}(z)$.
To simplify notation we drop the reference to $\Leb$.

\subsection{General setup}

We note that $\mathcal{L}_\epsilon(z)$ arising from both (i) the convolution with a locally supported stochastic kernel $q_\epsilon$ and (ii) convolution with a Fej{\'e}r kernel, can be considered as operators on $L^2(\mathbb{T}^2)$.
A numerical approximation $\mathcal{L}_\epsilon(z)$ of the twisted transfer operator $\mathcal{L}(z)$ can be formed in a number of ways, detailed below, but each of these will be based on Fourier approximation.
This is a natural approach as we have a periodic spatial domain and the map and observable are smooth.
First, we set up the Fourier function basis and $L^2$-orthogonal projection of the action of $\mathcal{L}_\epsilon$ on these basis functions.
Using the usual $L^2$ inner product $\langle f,g\rangle =\intf_{\mathbb{T}^2} f\cdot\overline{g} d\Leb$, for $x\in \mathbb{T}^2$ and $\mathbf{j}\in\mathbb{Z}^2$, define a complex Fourier basis $f_{\mathbf{j}}(x)=e^{2\pi i\mathbf{j}\cdot x}$, so that $g=\sum_{\mathbf{j}\in\mathbb{Z}^2} \langle g,f_\mathbf{j}\rangle f_\mathbf{j}:=\sum_{\mathbf{j}\in\mathbb{Z}^2} \hat{g}(\mathbf{j}) f_\mathbf{j}$.
To obtain a representation of $\mathcal{L}_\epsilon$ in this basis, we compute:
\begin{eqnarray}
\nonumber L(z)_{\epsilon,\mathbf{k}\mathbf{j}}:=\langle \mathcal{L}(z)_\epsilon f_\mathbf{k},f_\mathbf{j}\rangle&=&\intf_{\mathbb{T}^2}\left(\intf_{\mathbb{T}^2} q_\epsilon(x-y)\mathcal{L}(e^{z g(y)}f_\mathbf{k}(y)) d\Leb(y)\right)\overline{f_\mathbf{j}(x)} d\Leb(x)\\
\nonumber&=&\intf_{\mathbb{T}^2}\left(\intf_{\mathbb{T}^2} q_\epsilon(x-Ty)e^{z g(y)}f_\mathbf{k}(y) d\Leb(y)\right)\overline{f_\mathbf{j}(x)} d\Leb(x)\\
\nonumber&=&\intf_{\mathbb{T}^2}\left(\intf_{\mathbb{T}^2} q_\epsilon(x-Ty)\overline{f_\mathbf{j}(x)} d\Leb(x)\right)e^{z g(y)}f_\mathbf{k}(y) d\Leb(y)\\
\nonumber&=&\intf_{\mathbb{T}^2}\left(\intf_{\mathbb{T}^2} q_\epsilon(x)e^{-2\pi i\mathbf{j}\cdot(x+Ty)} d\Leb(x)\right)e^{z g(y)}f_\mathbf{k}(y) d\Leb(y)\\
\nonumber&=&\intf_{\mathbb{T}^2} \langle q_\epsilon,f_\mathbf{j}\rangle e^{-2\pi i\mathbf{j}\cdot (Ty)}e^{z g(y)}f_\mathbf{k}(y) d\Leb(y)\\
\nonumber&=&\langle q_\epsilon,f_\mathbf{j}\rangle\langle f_{-\mathbf{j}}\circ T\cdot{e^{z g}},f_\mathbf{-k}\rangle \\
\label{used}&=&\hat{q_\epsilon}(\mathbf{j})\cdot\widehat{\left({f_{-\mathbf{j}}}\circ T\cdot e^{z g}\right)}({\mathbf{-k}})
\end{eqnarray}
Notice that \eqref{used} only involves Fourier transforms of trivial objects (e.g.\ composition of a basis function with the map, exponential functions, the stochastic kernel, and the basis functions themselves).
To obtain spectral information for $\mathcal{L}_\epsilon(z):L^2\to L^2$ we may solve the generalised eigenvalue problem $L_\epsilon(z)v_\epsilon(z)=\lambda_\epsilon(z) v_\epsilon(z)$.

\subsection{Discrete Fourier transform}
To numerically approximate the above Fourier transforms, we first truncate the Fourier modes so that $\mathbf{j}\in\{-n/2+1,\ldots,-1,0,1,\ldots,n/2\}^2$, where $n=2^{n'}$ for some $n'\in\mathbb{Z}^+$.
Corresponding to this frequency grid is a regular spatial grid on $\mathbb{T}^2$ of the same cardinality; we call these frequency and spatial grids ``coarse grids''.
The $L^2$ inner products are estimated using MATLAB's two-dimensional discrete fast Fourier transform (DFT) \verb"fft2" on equispaced spatial and frequency grids with cardinalities $N=2^{N'}$ for some integer $N'\ge n'$;  these grids will be referred to as ``fine grids''.
The DFT is a collocation process, and by using $N\ge n$, we evaluate our functions on a finer spatial grid and produce more accurate estimates of the (lower) frequencies in the coarse grid.
One may also think of the DFT as a type of interpolation;  for fixed $n$, as $N$ increases we achieve increasingly accurate estimates of the $L^2$ inner products.
The cardinality $n^2$ of the coarse grid determines the size of the $n^2\times n^2$ matrix $L_{\epsilon,n(\epsilon)}(z)$ (if convolving with stochastic kernels) or $L_n(z)$ (if convolving with Fej{\'e}r kernels), while the cardinality $N^2$ of the fine grid determines the computation effort put into estimating the inner products via the DFT.
In our experiments we will use $n=32, 64, 128$ and $N=512$.

The kernel $q_\epsilon$ will be either:
\begin{enumerate}
\item A stochastic kernel given by an $L^1$-normalised $\mathcal{C}^\infty$ bump function with support restricted to the disk of radius $\epsilon$ centred at 0.  The particular bump function we use in the numerics is a well-known transformed version of a Gaussian given by $q_\epsilon(x)=(C/\epsilon^2)\exp(-1/(1-\|x/\epsilon\|^2)$ for $x\in B_\epsilon(0)$, where $C$ is a fixed $L^1$-normalising constant.
\item The square Fej{\'e}r kernel of order $n$.
Because of the special form of the square Fej{\'e}r kernel we have that ${\hat{q_n}}({\mathbf{j}})=(1-|\mathbf{j}_1|/(n/2+1))(1-|\mathbf{j}_2|/(n/2+1))$, which may be inserted directly into (\ref{used}).
Another advantage of the Fej{\'e}r kernel is that no explicit mollification is required, with the ``$\epsilon$'' slaved to the coarse resolution $n$.
\end{enumerate}
In our experiments, given a coarse frequency resolution $n$, we will try to select $\epsilon$ so that the stochastic kernel ``matches'' the Fej{\'e}r kernel.
We do this by choosing $\epsilon$ so that $\min_{\mathbf{j}\in\{-n/2+1,\ldots,n/2\}^2} |\hat{q}_\epsilon(\mathbf{j})|\approx \min_{\mathbf{j}\in\{-n/2+1,\ldots,n/2\}^2} |\hat{q}_n(\mathbf{j})|$.

\subsection{Numerical results}

\begin{figure}[hbt]
  \centering
  \includegraphics[width=15cm]{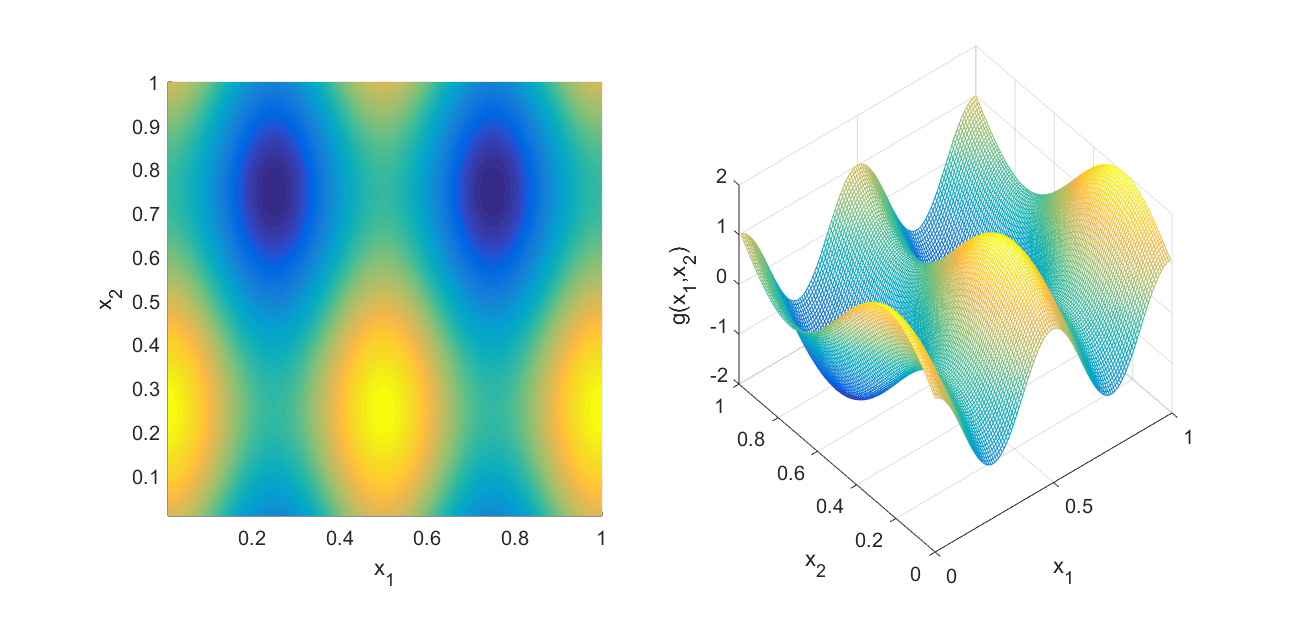}
  \caption{Graph of the observable used in the variance and rate function calculations.  Left: view from above. Right:  view from the side with the fine $N\times N$ spatial mesh  visible ($N=512$).}\label{fig:obs}
\end{figure}

The specific map $T:\mathbb{T}^2\circlearrowleft$ on which we carry out our numerics is a small perturbation of a linear toral automorphism:
$$T(x_1,x_2)=(2x_1+x_2+2\delta \cos(2\pi x_1),x_1+x_2+\delta\sin(4\pi x_2+1)),$$
with $\delta=0.01$.
By Proposition \ref{thm:perturbed_cat_props} we have that $T$ is Anosov and satisfies the conditions of Proposition \ref{prop:translate_bound_cond}, which is required in order to rigorously estimate the statistical properties of $T$ using the Fej{\'e}r kernel as per Corollary \ref{cor:fejer_stability}.
The observable we use when computing the variance and the rate function is $g(x_1,x_2)=\cos(4\pi x_1)+\sin(2\pi x_2)$, displayed in Figure \ref{fig:obs}.

\subsubsection{Estimating the SRB measure}
\begin{figure}[hbt]
  \centering
  \includegraphics[height=6cm]{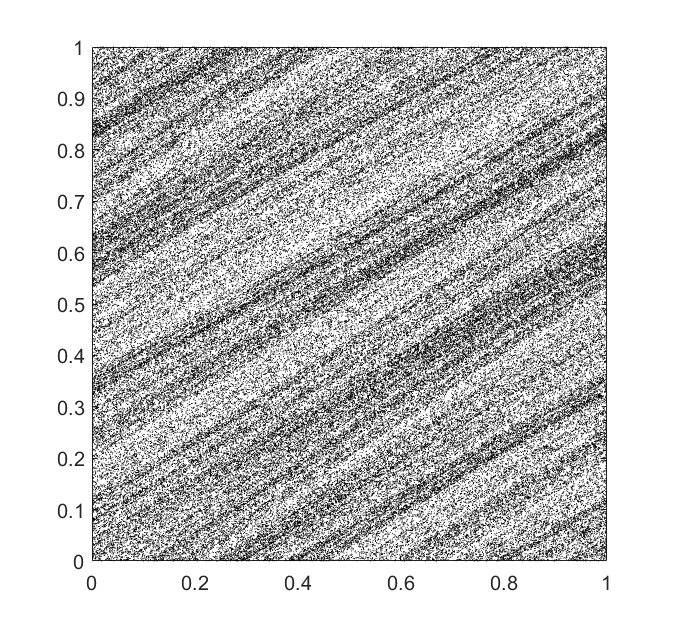}
  \caption{A trajectory of length $1.5\times 10^5$ initialised at a random location.}\label{fig:traj}
\end{figure}

\begin{figure}[hbt]
  \centering
  \includegraphics[height=6cm]{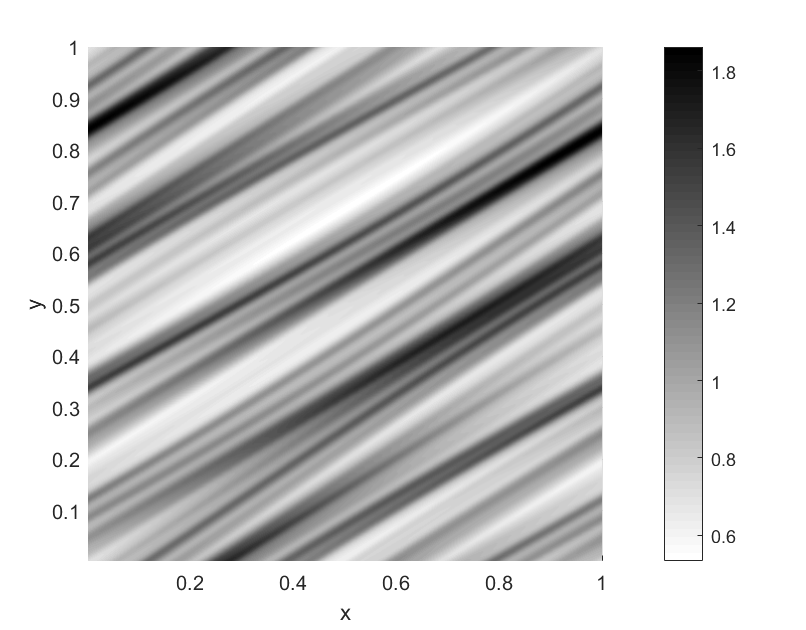}
  \includegraphics[height=7cm]{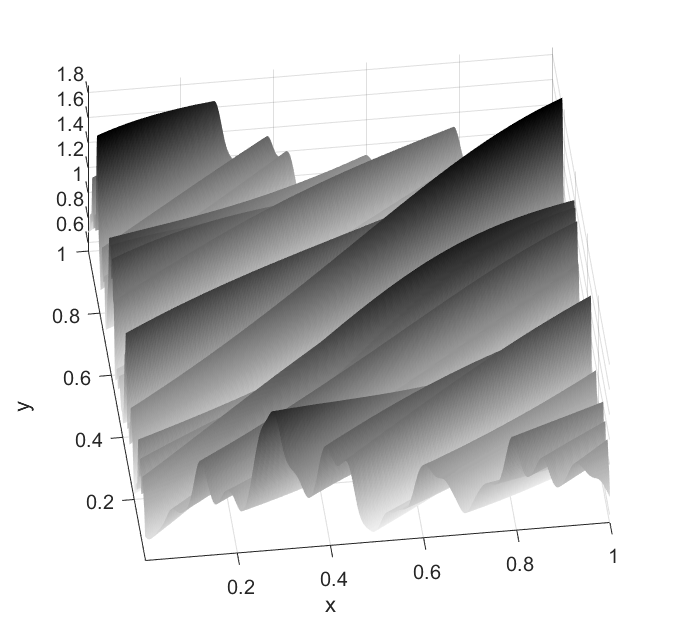}
  \caption{Approximations of the SRB measure computed as the leading right eigenvector of $L_{128}$ (Fej{\'e}r kernel of order 64), using a fine grid cardinality of $N=512$ to evaluate the Fourier transforms. Left: Darker regions indicate higher ``density''. Right: The same image rotated and represented in three dimensions with the vertical axis indicating ``density''.}\label{fig:fejersbr}
\end{figure}

Transitive Anosov systems possess a unique Sinai-Ruelle-Bowen (SRB) measure \cite[Theorem 1]{young2002srb}, which is exhibited by trajectories beginning in a full Lebesgue measure subset of $\mathbb{T}^2$.
A trajectory of length $1.5\times 10^5$ initialised at a random location is shown in Figure \ref{fig:traj}.
To create a numerical approximation of the SRB measure we compute the leading eigenvector of
$L_n$ (the matrix associated with the Fej{\'e}r kernel).
Figure \ref{fig:fejersbr} illustrates the results of using $n=128, N=512$.
The left panel of Figure \ref{fig:fejersbr} is shaded so that higher ``density'' is indicated by darker shading.
Note that this compares very well with the density of points in the trajectory shown in Figure \ref{fig:traj}, and that Figure \ref{fig:fejersbr} (left) captures many structures more clearly than the trajectory image.
The right panel of Figure \ref{fig:fejersbr} shows the same image as the left panel, but rotated and with the density plotted along the vertical axis.
The high degree of smoothness of the estimate of the SRB measure along unstable directions is evident.
Reducing $n$ from 128 to $n=64$ or $n=32$ has little effect on the image in unstable directions as these slow oscillations are still well-captured by lower order Fourier modes, but the higher frequency oscillations in stable directions will not be captured as well and the image will be ``smoothed'' in the stable directions.

\begin{figure}[hbt]
  \centering
  \includegraphics[height=6cm]{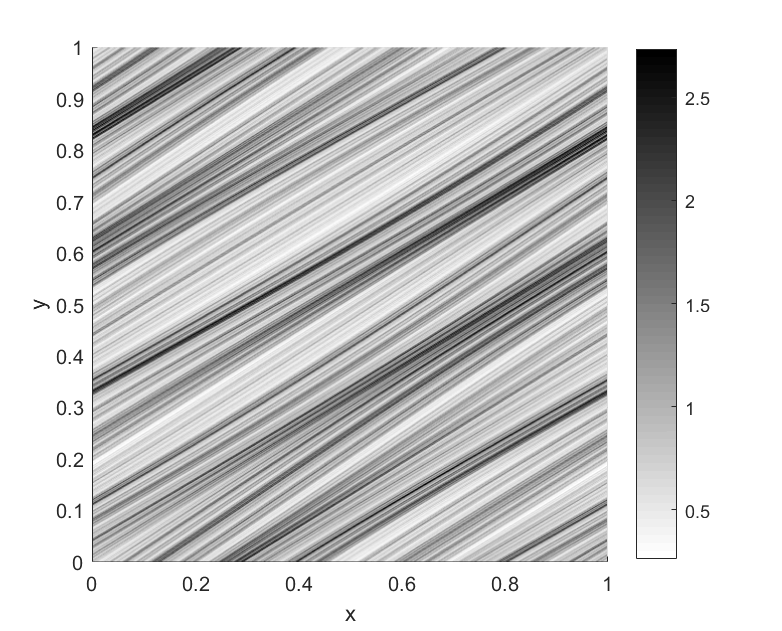}
  \includegraphics[height=7cm]{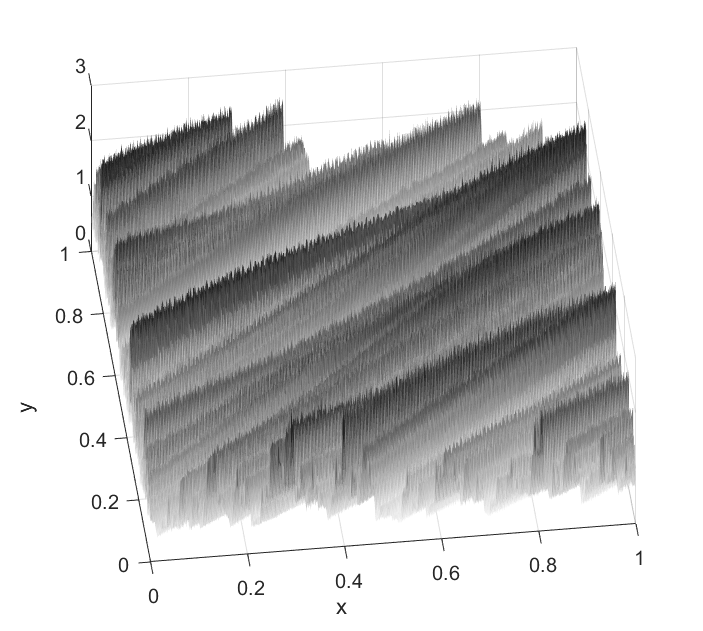}
  \caption{Approximations of the SRB measure computed as the leading left eigenvector of $P_{512}$ (Ulam grid of size $512\times 512$). Left: Darker regions indicate higher ``density''. Right: The same image rotated and represented in three dimensions with the vertical axis indicating ``density''.}\label{fig:ulamsbr}
\end{figure}

As a non-rigorous comparison, we form an Ulam matrix using a $512\times 512$ equipartition of boxes $\{B_1,\ldots,B_{2^{18}}\}$ on $\mathbb{T}^2$.
We compute a row-stochastic matrix $P_{512}$ as $P_{512,ij}=\Leb(B_i\cap T^{-1}B_j)/\Leb(B_i)$, where the entries $P_{512,ij}$ are estimated by uniformly sampling 1600 points in each box and counting the fraction of points initialised in $B_i$ that have their image in $B_j$.
The Ulam estimate of the SRB measure is then obtained as the leading left eigenvector of $P_{512}$.
The images corresponding to Figure \ref{fig:fejersbr} are shown in Figure \ref{fig:ulamsbr}.
In comparison to Figure \ref{fig:fejersbr} two things are noticeable.
Firstly, Figure \ref{fig:ulamsbr}(left) appears to produce a finer representation of the SRB measure than Figure \ref{fig:fejersbr}(left), and secondly, the estimate in Figure \ref{fig:ulamsbr}(right) is rougher in unstable directions than the estimate in Figure \ref{fig:fejersbr}.
Each of these observations is relatively easy to explain at a superficial level through the different approximation bases used.
In terms of regularity of approximation basis the Ulam method is very low order (piecewise constant) because it uses a basis of indicator functions on the $512\times 512$ grid.
On the other hand, the approximation basis for the Fourier approximation of very high order (analytic).
The Ulam basis is thus very flexible and can adapt well to the roughness of the SRB measure in stable directions, but has no apriori smoothness in unstable directions.
In contrast, the Fourier basis is less flexible in stable directions, requiring more modes to capture rapid oscillations, but is extremely efficient at approximating smooth functions and easily captures the smooth variation in unstable directions.
A recent alternative non-rigorous collocation-based method of SRB measure approximation has been explored in \cite{dynamicmode2019} for certain families (Blaschke products) of analytic Anosov maps. In the case of analytic expanding maps, \cite{dynamicmode2019} proves that this method produces the true absolutely continuous invariant measure in the limit of increasing numerical resolution.

\subsubsection{Estimating the variance}
To estimate the variance of a centred observable $g:\mathbb{T}^2\to\mathbb{R}$ we employ the method described in \cite[Section 4.2]{cfstability}.
We use the representation of the (approximate) variance
\begin{equation}
\label{vareqn}
\sigma^2_n:=\lambda^{(2)}_n(0)=\intf_{\mathbb{T}^2} g^2v_n(0)+2g(\mathrm{Id}-\mathcal{L}_n(0))^{-1}\mathcal{L}_n(0)(gv_n(0))\ d\Leb.
\end{equation}
The main difference to the calculations in \cite{cfstability} is that here we use Fourier approximation, whereas \cite{cfstability} used Ulam's method, which was better suited to the piecewise expanding maps considered there.
The key computational component in (\ref{vareqn}) is the solution of a single linear equation to obtain an estimate for $(\mathrm{Id}-\mathcal{L}_n(0))^{-1}\mathcal{L}_n(0)(gv_n(0))$, which is $\frac{d}{dz}v_n(z)$ at $z=0$.
Because our approximate transfer operator $L_n(z)$ is represented in frequency space, we set up and solve this linear equation in frequency space, yielding the DFT of $\frac{d}{dz}v_n(z)$ at $z=0$.
This Fourier transform is then inverted with the inverse DFT to produce the required spatial estimate.
Similarly, the DFT of $v_n(0)$ is computed as the leading right eigenvector of $L_n(0)$ and inverted with the inverse DFT to obtain a spatial estimate of $v_n(0)$.
These two spatial estimates (analytic functions consisting of a linear combination of Fourier modes) are then evaluated on the fine spatial grid and the integral in (\ref{vareqn}) is computed as a simple Riemann sum over the fine spatial grid.
The Ulam-based variance estimates are calculated identically to \cite{cfstability}.
In Table \ref{tab:variance} we report variance estimates over a range of coarse grid resolutions to roughly indicate the dependence of the estimates on grid resolution.
\begin{table}[hbt]
  \centering
  \begin{tabular}{|l|l|l|l|l|}
    \hline
    Coarse grid & & & & \\
     resolution $n$ & $n=32$ & $n=64$  & $n=128$  & $n=256$  \\ \hline\hline
    Stochastic kernel & 0.9359 ($\epsilon=0.0693$) &  0.9342 ($\epsilon=0.0378$)   &  0.9337 ($\epsilon=0.0210$)   &  \\
    \hline
    Fej{\'e}r kernel & 0.9447 & 0.9395 & 0.9366 &  \\ \hline
    Ulam  &  & 0.9320 & 0.9307 & 0.9348 \\
    \hline
  \end{tabular}
  \caption{Variance estimates from the two Fourier approximation approaches and Ulam's method. For the two Fourier approximation methods we use a fine frequency grid resolution of $N=512$ and for Ulam's method we use $1600$ sample points per box.}\label{tab:variance}
\end{table}

\subsubsection{Estimating the rate function}
\begin{figure}[hbt]
  \centering
  \includegraphics[width=13cm]{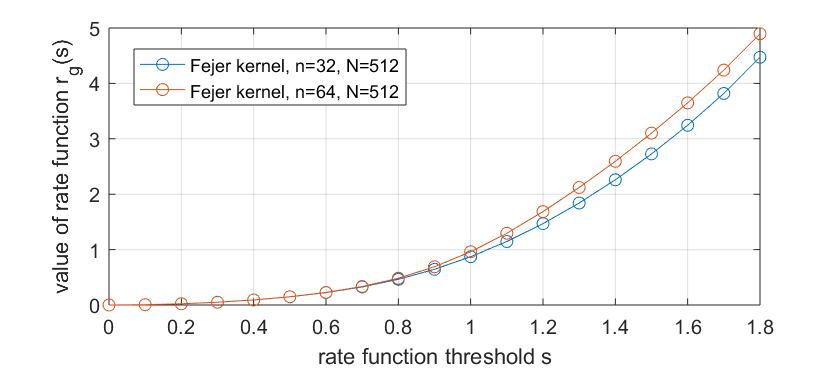}
  \caption{Estimates of the rate function $r_g$ using Fej{\'e}r kernels with $n=32, 64$ and $N=512$.}\label{fig:ratefunction}
\end{figure}

We numerically estimate the rate function $r_g(s)=\sup_{z\in V}(sz-\log|\lambda_\epsilon(z)|)$ for a centred observable $g:\mathbb{T}^2\to \mathbb{R}$ using the Fej{\'e}r kernel approach.
We create the Fej{\'e}r kernel estimate $L_n(z)$ of the twisted transfer operator and compute the leading eigenvalue $\lambda_n(z)$ and eigenvector $v_n(z)$ of $L_n(z)$.
The leading eigenvector $v_n(z)$ is converted from frequency space to a function on $\mathbb{T}^2$ by evaluating the linear combination (according to the entries of $v_n(z)$) of the $n$ associated Fourier modes on a fine $N\times N$ spatial grid for $N=512$.
For a given $s$, we are now in a position to find the minimum of $-(sz-\log|\lambda_\epsilon(z)|)$ as a function of $z$.
We used MATLAB's \verb"fminunc" routine (unconstrained function minimisation) with the default quasi-newton option, which takes around four to five iterates to converge to the minimum within a preset tolerance of $10^{-6}$.
We asked for the values of $r_g(s)$ for $s$ between 0 and 1.8 in steps of 0.1, and initialised the search for the optimal $z$ value using the optimal $z$ from the previous value of $s$.
The results are shown in Figure \ref{fig:ratefunction} for coarse grids of size $n=32$ and $n=64$, with fine grid collocation and function evalution using $N=512$.
Note that the range of $g$ is $[-2,2]$ (see also Figure \ref{fig:obs}), and that $g$ is already centred with respect to Lebesgue measure on the 2-torus.
In the rate function computations we centre $g$ according to our estimate of the SRB measure, but do not expect the range of $g$ to vary significantly.
The large values of $r_g(s)$ as $s$ approaches $2$ are consistent with this observation.

\appendix

\section{Proofs for Section \ref{sec:stability_anosov}}
\label{sec:proofs_sec_4}

In this appendix we prove Propositions \ref{prop:defn_of_twist_anosov} and \ref{prop:anosov_hypothesis}, and Theorem \ref{thm:anosov_stat_stability}.

\begin{proof}[{The proof of Proposition \ref{prop:defn_of_twist_anosov}}]
  It is clear that $M_g(0)$ is the identity.
  For each $k \in \N$ let $P_k : \mathcal{C}^r(X, \C) \to \mathcal{C}^r(X, \C)$ be defined by $P_k f = g^k f $.
  Multiplication by $g$ is continuous on $B^{p,q}_\C$ by \cite[Lemma 3.2]{gouezel2006banach} and so $P_k \in L(B^{p,q}_\C)$ for each $k \in \N$.
  Moreover, as $\norm{P_k}_{p,q} \le \norm{P_1}_{p,q}^k$, for each $z \in \C$ the series $\sum_{k=0}^\infty z^k P_k$ is absolutely convergent in the $\norm{\cdot}_{p,q}$ operator norm, with limit $M_g(z)$. Hence $z \mapsto M_g(z)$ is a well-defined analytic map taking values in $L(B^{p,q}_\C)$.
  The same argument holds when $B^{p,q}_\C$ is replaced with $B^{p-1,q+1}_\C$, and so $z \mapsto M_g(z)$ is analytic on $L(B^{p-1,q+1}_\C)$ too. In particular, it is compactly $\norm{\cdot}_{p-1,q+1}$-bounded.
\end{proof}

From the beginning of \cite[Section 4]{gouezel2006banach}, for each $h \in B^{p,q}$ and $\phi \in \mathcal{C}^q(X,\R)$ we have
\begin{equation*}
  \abs{\intf h \phi d \Omega} \le C \norm{h}_{p,q} \norm{\phi}_{\mathcal{C}^q},
\end{equation*}
for some $C > 0$ independent of $h$ and $\phi$. It is straightforward to show that the same inequality holds for $h \in B^{p,q}_\C$ and $\phi \in \mathcal{C}^q(X,\C)$ (although with a different $C$, which is inconsequential).
Hence, the functional $\Omega$ is in $(B^{p,q}_\C)^*$. Let $m \in B^{p,q}$ be a probability measure (i.e. the image of $m$ under the inclusion map from $B^{p,q}$ to $\mathcal{D}_q'(X)$ is a probability measure) and $h \in \mathcal{C}^r(X,\C)$.
Since $\mathcal{C}^r(X,\R)$ is dense in $B^{p,q}$, there exists $\{m_i\}_{i \in \Z^+} \subseteq \mathcal{C}^r(X,\R)$ such that $m_i \to m$ in $B^{p,q}$.
As $B^{p,q}$ is continuously injected into $\mathcal{D}_q'(X)$ it follows that $m_i \to m$ in $\mathcal{D}_q'(X)$.
Note that $\varphi \in \mathcal{C}^r(X, \R)$ naturally induces a measure, which we will also denote by $\varphi$, so that $\intf f d\varphi = \intf f \varphi d\Omega$ for each Borel measurable function $f:X \to \C$.
Hence,
\begin{equation}\label{eq:measure_identity}
  \Omega(h m) = \lim_{n \to \infty} \Omega(h m_n) = \lim_{n \to \infty} \intf h m_n d\Omega = \lim_{n \to \infty} \intf h dm_n  = \intf h dm.
\end{equation}

\begin{proposition}\label{prop:twist_coding}
  Let $g \in \mathcal{C}^{r}(X, \R)$, $S_n(g) = \sum_{k=0}^{n-1} g \circ T^k$ and $m \in B^{p,q}$ be a probability measure. Then for each $n \in \N$ and $z \in \C$ we have
  \begin{equation*}
     \intf e^{z S_n(g)} dm = \Omega(\LL(z)^n m).
  \end{equation*}
  \begin{proof}
    For $h \in \mathcal{C}^r(X,\C)$ we have $\LL h = (h \abs{\det T}^{-1}) \circ T^{-1}$. It is straightforward to verify that for every $f_1, f_2 \in \mathcal{C}^r(X,\C)$ we have $\LL( f_1 \circ T \cdot f_2) = f_1 \LL(f_2)$.
    By \cite[Lemma 3.2]{gouezel2006banach}, multiplication by $f_1$ is continuous on $B^{p,q}_\C$. Hence, by passing to the completion we may take $f_2 \in B^{p,q}_\C$.
    Setting $f_1 = e^{zg}$ and $f_2 = m$, and then inductively using this identity, it follows that for each $n \in \N$ and $z \in \C$ we have $\LL(z)^n m = \LL^n\left(e^{zS_n(g)} m\right)$.
    Upon integrating, and using \eqref{eq:measure_identity} and that $\LL$ preserves $\Omega$-integrals, we have
    \begin{equation*}
      \Omega(\LL(z)^n m) = \Omega\left(\LL^n\left(e^{zS_n(g)} m\right)\right) = \Omega\left(e^{zS_n(g)} m\right) = \intf e^{zS_n(g)} dm.
    \end{equation*}
  \end{proof}
\end{proposition}

In the setting of Proposition \ref{prop:anosov_hypothesis}, there is a well-known connection between the condition that $g$ is not an $L^2(\mu)$-coboundary and the map $z \mapsto \rho(\LL(z))$ being strictly convex in a real neighbourhood of $0$.

\begin{lemma}\label{lemma:coboundary_cond}
  If $g \in \mathcal{C}^r(X,\R)$ is not an $L^2(\mu)$-coboundary and satisfies $\intf g d\mu = 0$, then $z \mapsto \ln\rho(\LL(z))$ is strictly convex in a real neighbourhood of $0$.
  \begin{proof}
    We will reduce the statement to a well-known result.
    As $z \mapsto \LL(z)$ is analytic, standard analytic perturbation theory for linear operators implies that $\LL(z)$ is quasi-compact and simple in an open complex neighbourhood $U$ of $0$.
    Moreover, $\LL(z)$ has decomposition $\LL(z) = \lambda(z) \Pi(z) + N(z)$ where $\lambda(z)$, $\Pi(z)$ and $N(z)$ depend analytically on $z$.
    Let $\mathcal{P}$ denote the completion of $\{ f \in \mathcal{C}^r(X,\R) : f \ge 0\}$ in $B^{p,q}_\C$.
    Since $g$ is real-valued and $\LL$ preserves $\mathcal{P}$ in $\mathcal{C}^r(X,\R)$, it follows that $\LL(z)$ preserves $\mathcal{P}$ for each $z \in \R$.
    It is clear that $\mu \in \mathcal{P}$ and that $\Omega f \ge 0$ for all $f \in \mathcal{P}$.
    Since $\Omega (\Pi(0) \mu) = 1$, by possibly shrinking $U$ we may assume that $\Omega \Pi(z) \mu \ne 0$ for any $z \in U$.
    Hence $\lim_{n \to \infty} \lambda(z)^{-n} \Omega( \LL^n(z) \mu) = \Omega(\Pi(z)\mu)$ for every $z \in U$.
    As $\Omega (\LL^n(z) \mu) \in [0, \infty)$ for every $n \in \Z^+$ and $z \in U$, and $\mathcal{P}$ is closed, it follows that $\Omega(\Pi(z)\mu) \in (0, \infty)$.
    Hence there exists a branch of the complex argument $\arg : \C \to \R$ that is continuous at $\Omega(\Pi(z)\mu)$. It follows that $\arg(\lambda(z)^{-n}) = \arg(\lambda(z)^{-n} \Omega (\LL^n(z) \mu)) \to \arg(\Omega(\Pi(z)\mu))$, which is only possible if $\lambda(z)$ is real. Thus for $z \in \R \cap U$ we have $\rho(\LL(z)) = \abs{\lambda(z)} = \lambda(z)$.
    As $\lambda(0) = 1$, by possibly shrinking $U$ again we may assume that $\lambda(z)$ is bounded away from 0 when $z$ ranges over $U$, implying that $\ln \rho(\LL(z))$ is finite for every $z\in U \cap \R$.

    Differentiating twice yields
    \begin{equation}\label{eq:coboundary_cond_1}
      \frac{\mathrm{d}^2}{\mathrm{d}z^2} \ln \rho(\LL(z)) = \frac{\mathrm{d}^2}{\mathrm{d}z^2} \ln \lambda(z) =  \frac{\lambda''(z) \lambda(z) - (\lambda'(z))^2}{\lambda(z)^2}.
    \end{equation}
    After possibly shrinking $U$ there exists an analytic map $z \mapsto \mu(z)$, defined for $z \in U$, such that $\mu(z)$ is an eigenvector of $\LL(z)$ for the eigenvalue $\lambda(z)$. As $\Omega(\mu(0)) = 1$, after possibly shrinking $U$ we may choose each $\mu(z)$ so that $\intf \mu(z) d \Omega = 1$.
    By differentiating $(\lambda(z) - \LL(z))\mu(z) = 0$, we have
    \begin{equation}\label{eq:coboundary_cond_2}
      \lambda'(z)\mu(z) = \LL'(z)\mu(z) - (\lambda(z) - \LL(z))\mu'(z).
    \end{equation}
    Note that $\LL'(z)f = \LL(z) (gf)$. Evaluating \eqref{eq:coboundary_cond_2} at $z=0$ and integrating with respect to $\Omega$ yields
    \begin{equation*}
      \lambda'(0) = \intf \LL'(0)d\mu = \intf g d\mu = 0.
    \end{equation*}
    Evaluating \eqref{eq:coboundary_cond_1} at $z = 0$, we have
    $\restr{\frac{\mathrm{d}^2}{\mathrm{d}z^2} \ln \rho(\LL(z))}{z=0} = \lambda''(0)$.
    By \cite[Lemma IV.3]{hennion2001limit}\footnote{There is a discrepancy in sign between our expressions that is due to the map $z \mapsto \LL(iz)$ being used in \cite[Chapter IV]{hennion2001limit}, whereas we use the map $z \mapsto \LL(z)$.} we have
    \begin{equation*}
      \lambda''(0) = \lim_{n \to \infty} \intf \frac{S_n^2}{n} d\mu.
    \end{equation*}
    Hence $\lambda''(0) \ge 0$ and it suffices to show that if $\lambda''(0) = 0$ then $g$ is an $L^2(\mu)$-coboundary, as the contrapositive then implies that $z \mapsto \ln\rho(\LL(z))$ is strictly convex in a real neighbourhood of $0$. The proof of this is well-known; see e.g. \cite[Lemma 6]{rousseau1983theoreme}.
  \end{proof}
\end{lemma}

\begin{proof}[{The proof of Proposition \ref{prop:anosov_hypothesis}}]
  We consider $Y_k = g \circ T^k$ on the probability space $(X, m)$. The expectation of $Y_k$ with respect to $m$ is
  \begin{equation*}
    \intf g \circ T^k dm = \Omega\left((g \circ T^k)m \right) = \Omega\left((\LL^k m) g \right).
  \end{equation*}
  By our assumptions, $\LL = \LL(0)$ is a simple quasi-compact operator on $B^{p,q}_\C$ with $\rho(\LL) = 1$.
  Let $\LL^k = \Pi + N^k$ be the quasi-compact decomposition of $\LL^k$. As $\mu$ is $T$-invariant and $\Omega\LL = \Omega$, it follows that $\Pi(f) = \Omega(f) \mu$. Hence, as $N^k \to 0$ and $m$ is a probability measure, we have $\LL^k m = \Omega(m)\mu + N^k \to \mu$ in $B^{p,q}_\C$.
  Using \eqref{eq:measure_identity} and the fact that $\Omega \in (B^{p,q}_\C)^*$ we have
  \begin{equation*}
    \lim_{k \to \infty} \intf g \circ T^k dm = \lim_{k \to \infty} \intf g \cdot (\LL^k m) d\Omega = \Omega( \mu g ) = \intf g d\mu = 0.
  \end{equation*}
  It follows that $\lim_{n \to \infty} \E(S_n) / n = 0$, where the expectation is taken with respect to $m$. Let $(E, \norm{\cdot}) = (B^{p,q}_\C, \norm{\cdot}_{p,q})$.
  The map $\LL(\cdot) : \C \to L(B^{p,q}_\C)$ induced by the twist $M_g$ is analytic by Proposition \ref{prop:defn_of_twist_anosov}, and $z \mapsto \ln\rho(\LL(z))$ is strictly convex in a real neighbourhood of 0 by Lemma \ref{lemma:coboundary_cond}.
  Hence, in view of Proposition \ref{prop:twist_coding} and as $\Omega \in (B^{p,q}_\C)^*$, Hypothesis \ref{hypothesis:analytic_coding} holds.
\end{proof}

\begin{proof}[{The proof of Theorem \ref{thm:anosov_stat_stability}}]
  For $t \in [0,1]$ let $\LL_t$ denote the Perron-Frobenius operator induced by $T(t)$. By \cite[Theorem 2.3]{gouezel2006banach}, topological transitivity of $T(0)$ implies that $\LL_0$ is a simple quasi-compact operator on $B^{p,q}_\C$ with $\rho(\LL_0) = 1$.
  By \cite[Section 7]{gouezel2006banach}, there exists some $t' > 0$ for which $\{\LL_t\}_{t \in [0, t']}$ satisfies (KL) on $B^{p,q}_\C$ with $\wnorm{\cdot} = \norm{\cdot}_{p-1,q+1}$.
  Applying Theorem \ref{thm:stability_of_statistical_data} to $\{g \circ T^k \}_{k \in \N}$, which satisfies Hypothesis \ref{hypothesis:analytic_coding} by Proposition \ref{prop:anosov_hypothesis}, we obtain $\theta > 0$ and $\epsilon \in (0, t')$ so that whenever $t \in [ 0, \epsilon]$ and $z \in D_\theta$ the operator $\LL_t M_g(z)$ is quasi-compact and simple with leading eigenvalue $\lambda_t(z)$.
  In particular, $1$ is a simple eigenvalue of $\LL_t$ for $t \in [ 0, \epsilon]$ and so $T(t)$ has a unique SRB measure $\mu_t$ in $B^{p,q}_\C$.
  By \cite[Proposition 2.4, Remark 2.5]{froyland2014detecting}, for each $t \in (0, \epsilon)$
  there exists an eigenvector $v_t$ of $\LL_t$ associated to the eigenvalue 1 such that $v_t \to \mu$ in $B_\C^{p-1,q+1}$ as $t \to 0$.
  By simplicity of the eigenvalue 1, for sufficiently small $t$ we have $\mu_t = \frac{v_t}{\intf v_t d \Omega}$ and so the continuity of $f \mapsto \intf f d\Omega$ on $B_\C^{p-1,q+1}$ implies that
  $\mu_t \to \mu$ in $B_\C^{p-1,q+1}$ too.

  Fix $t \in [ 0, \epsilon]$. Let $A_t = \intf g d\mu_t$ and $g_t = g - A_t$.
  As multiplication by $g$ is continuous on $B^{p-1,q+1}_\C$ \cite[Lemma 3.2]{gouezel2006banach} and $\Omega \in (B^{p-1,q+1}_\C)^*$, we have $\lim_{t \to 0} A_t = \lim_{t \to 0} \Omega(g \mu_t) = \Omega(g\mu) = \intf g d\mu = 0$.
  Note that $e^{zA_t}\LL_t M_{g_t}(z) = \LL_t M_{g}(z)$ for every $z \in \C$, and so $\LL_t M_{g_t}(z)$ is quasi-compact exactly when $\LL_t M_{g}(z)$ is.
  In particular, $\LL_t M_{g_t}(z)$ is a simple quasi-compact operator for every $z \in D_\theta$ with leading eigenvalue $\kappa_t(z) = e^{-zA_t} \lambda_t(z)$.
  From the material in this section and the last, it is routine to verify that
  $\{g_t \circ T(t)^k\}_{k \in \N}$ satisfies Hypothesis \ref{hypothesis:analytic_coding} on $(X,m)$ with operator-valued map $z \mapsto \LL_t M_{g_t}(z)$ and Banach space $B^{p,q}_\C$.
  Hence, by Theorem \ref{thm:naive_ng_ldp} the sequence $\{g_t \circ T(t)^k\}_{k \in \N}$ satisfies a LDP on $(X,m)$ with rate function $r_t : J_t \to \R$ defined by
  \begin{equation*}
    r_t(s) = \sup_{z \in (-\theta, \theta)} (s z - \ln \abs{\kappa_t(z)}).
  \end{equation*}
  Recall from \eqref{eq:domain_of_rate_function} that $r_t$ has domain
  \begin{equation*}
    J_t = \left(\restr{\frac{\mathrm{d}}{\mathrm{d}z} \ln\abs{\kappa_t(z)}}{z = -\theta}, \restr{\frac{\mathrm{d}}{\mathrm{d}z} \ln\abs{\kappa_t(z)}}{z = \theta}\right).
  \end{equation*}
  As $\kappa_t(z) = e^{-zA_t} \lambda_t(z)$, we therefore have $J_t = J_0 - A_t$ and
  \begin{equation*}
    r_t(s) = \sup_{z \in (-\theta, \theta)} ((s + A_t)z  - \ln \abs{\lambda_t(z)}).
  \end{equation*}
  By Theorem \ref{thm:stability_of_statistical_data}, for each compact $U \subseteq J_0$ there is a closed interval $V \subseteq (-\theta,\theta)$ so that the map $s \mapsto \sup_{z \in V} (s z  - \ln \abs{\lambda_t(z)})$ converges uniformly to $r_0$ on $U$.
  Since the map $z \mapsto \ln \abs{\lambda_t(z)}$ is convex on $(-\theta,\theta)$, by the arguments from \cite[Proposition 3.10]{cfstability} we have
  \begin{equation*}
    r_t(s - A_t) = \sup_{z \in V} (s z  - \ln \abs{\lambda_t(z)}).
  \end{equation*}
  Hence, $r_t \circ \tau_{-A_t} \to r_0$ compactly on $J_0$.
\end{proof}

\section{Spectral stability results for $\LL_{\Leb}$}
\label{app:stability_for_leb_pf}

In this appendix we prove the claims made in Remarks \ref{remark:omega_v_leb} and \ref{remark:omega_v_leb_fejer} regarding the relationship between the operators $\LL_{\Omega}(z)$ and $\LL_{\Leb}(z)$, and that the spectral stability results for $\LL_{\Omega}(z)$ from Sections \ref{sec:approximation_anosov} and \ref{sec:fejer_anosov} imply the spectral stability results for $\LL_{\Leb}(z)$.

\begin{proposition}\label{prop:conjugacy}
  Let $R : B^{p,q}_\C \to B^{p,q}_\C$ be defined
  \begin{equation*}
    R h = h \cdot \frac{\mathrm{d}\Omega}{\mathrm{d}\Leb}.
  \end{equation*}
  Then $R, R^{-1} \in L(B^{p,q}_\C) \cap L(B^{p-1,q+1}_\C)$ and for every $z \in \C$ we have
  \begin{equation}\label{eq:conjugacy}
    \LL_{\Omega}(z) = R^{-1}\LL_{\Leb}(z) R.
  \end{equation}
  Hence $\sigma(\LL_{\Omega}(z)) = \sigma(\LL_{\Leb}(z))$ for every $z \in \C$.
  \begin{proof}
    The fact that $R \in L(B^{p,q}_\C) \cap L(B^{p-1,q+1}_\C)$ follows from multiplication by $\mathcal{C}^r(\mathbb{T}^d,\R)$ functions being continuous on $B^{p,q}_\C$ and $B^{p-1,q+1}_\C$ \cite[Lemma 3.2]{gouezel2006banach}.
    Since $\frac{\mathrm{d}\Leb}{\mathrm{d}\Omega} \in \mathcal{C}^r(\mathbb{T}^d,\R)$, the same argument implies that $R^{-1}$ exists and is an element of $L(B^{p,q}_\C)$ and $L(B^{p-1,q+1}_\C)$.

    Let $f, h \in \mathcal{C}^r(\mathbb{T}^d, \C)$. By definition we have
    \begin{equation*}\begin{split}
      \intf \LL_\Omega(f) \cdot h d\Omega = \intf f \cdot (h \circ T) d\Omega &= \intf \left(f \frac{\mathrm{d}\Omega}{\mathrm{d}\Leb}\right) \cdot h \circ T d\Leb \\
      &= \intf \LL_{\Leb}\left(f \frac{\mathrm{d}\Omega}{\mathrm{d}\Leb}\right) \cdot  \frac{\mathrm{d}\Leb}{\mathrm{d}\Omega} \cdot h d\Omega \\
      &= \intf (R^{-1}\LL_{\Leb}R)(f) \cdot h d\Omega.
    \end{split}\end{equation*}
    Hence $(R^{-1}\LL_{\Leb}R)f = \LL_{\Omega}f$ for all $f \in \mathcal{C}^r(\mathbb{T}^d, \C)$, and so the same identity holds on $B^{p,q}_\C$ by density.
    The conjugacy relation \eqref{eq:conjugacy} holds for the twisted operators due to the untwisted conjugacy relation and the definition of the twist $M_g$ (see Proposition \ref{prop:defn_of_twist_anosov}).
    One has $\sigma(\LL_{\Omega}(z)) = \sigma(\LL_{\Leb}(z))$ immediately from \eqref{eq:conjugacy}.
  \end{proof}
\end{proposition}

\begin{proposition}\label{prop:Omega_to_Leb_KL}
  Let $\{k_\epsilon\}_{\epsilon \in (0, \epsilon_0)} \subseteq L^1(\Leb)$ be an $L^1(\Leb)$-bounded family. Set $\LL_{\Omega, \epsilon} = k_\epsilon \ast \LL_{\Omega}$ and $\LL_{\Leb, \epsilon} = k_\epsilon \ast \LL_{\Leb}$.
  Let $\LL_{\Omega, 0}= \LL_{\Omega}$ and $\LL_{\Leb, 0}= \LL_{\Leb}$.
  Suppose that $\{\LL_{\Omega, \epsilon}\}_{\epsilon \in [0, \epsilon_0)}$ satisfies (KL) and that one of the following conditions holds.
  \begin{enumerate}[label=(K\arabic*)]
    \item \label{en:K1} For every $\eta > 0$ there exists $\epsilon_\eta$ such that $\supp k_\epsilon \subseteq B_\eta(0)$ for every $\epsilon \in (0, \epsilon_\eta)$.
    \item \label{en:K2} $\LL_\Omega$ satisfies \eqref{eq:translate_bound} and for every $\eta > 0$ we have $\lim_{\epsilon \to 0} \intf_{\mathbb{T}^d \setminus B_\eta(0)} \abs{k_\epsilon} d\Leb = 0$
  \end{enumerate}
  Then $\{\LL_{\Leb, \epsilon}\}_{\epsilon \in [0, \epsilon_1)}$ satisfies (KL) for some $\epsilon_1 \in (0, \epsilon_0)$.
  \begin{proof}
    We have
    \begin{equation*}
      \LL_{\Leb, \epsilon} = R \LL_{\Omega,\epsilon} R^{-1} + (\LL_{\Leb, \epsilon} - R \LL_{\Omega,\epsilon} R^{-1}) := A_{\epsilon} + F_{\epsilon}.
    \end{equation*}
    We will prove that $\{A_{\epsilon}\}_{\epsilon \in [0, \epsilon_0)}$ satisfies (KL), and that there exists $\epsilon' \in (0, \epsilon_0)$ such that $\{F_\epsilon \}_{\epsilon \in [0, \epsilon')}$ satisfies the conditions required by Proposition \ref{prop:kl_norm_perturb}, which will then imply that $\{\LL_{\Leb, \epsilon}\}_{\epsilon \in [0, \epsilon_1)}$ satisfies (KL) for some $\epsilon_1 \in (0, \epsilon')$.

    It is straightforward to confirm that $\{A_{\epsilon}\}_{\epsilon \in [0, \epsilon_0)}$ satisfies (KL) by using (KL) for $\{\LL_{\Omega, \epsilon}\}_{\epsilon \in [0, \epsilon_0)}$, the conjugacy identity \eqref{eq:conjugacy} and the properties of the map $R$ as given in Proposition \ref{prop:conjugacy}. For example (KL1) follows from the estimate
    \begin{equation*}
      \tnorm{A_{\epsilon} - \LL_{\Leb}} = \tnorm{R
      (\LL_{\Omega, \epsilon} - \LL_{\Omega})R^{-1}} \le \norm{R}_{L(B^{p-1,q+1}_\C)}
      \tnorm{\LL_{\Omega, \epsilon} - \LL_{\Omega}}\norm{R^{-1}}_{L(B^{p,q}_\C)}.
    \end{equation*}
    The proofs of (KL2) and (KL3) follow from similar arguments.

    We will now prove that there exists $\epsilon' \in [0, \epsilon_0)$ such that $\{F_\epsilon \}_{\epsilon \in [0, \epsilon')}$ satisfies the conditions required by Proposition \ref{prop:kl_norm_perturb}.
    For brevity let $s = \frac{\mathrm{d}\Omega}{\mathrm{d}\Leb}$. Let $h \in \mathcal{C}^r(\mathbb{T}^d, \R)$, $k \le p$ be a non-negative integer, $W \in \Sigma$, $\{v_i\}_{i=1}^k \subseteq \mathcal{V}^r(W)$ with $\norm{v_i}_{\mathcal{C}^r} \le 1$, and $\varphi \in \mathcal{C}^{q + k}_0(W, \R)$ with $\norm{\varphi}_{\mathcal{C}^{q + k}} \le 1$.
    Since
    \begin{equation*}
      F_\epsilon R h =
      R(k_\epsilon \ast (\LL_{\Omega}h)) - k_\epsilon \ast (R\LL_{\Omega}h),
    \end{equation*}
    we have
    \begin{equation*}\begin{split}
      (F_\epsilon R h)(x)
      &= s(x) \intf k_\epsilon(y)(\LL_{\Omega}h)(x-y)d\Leb(y) - \intf k_\epsilon(y) s(x-y)(\LL_{\Omega}h)(x-y) d\Leb(y)\\
      &= \intf k_\epsilon(y)  (\LL_\Omega h)(x-y) \left( s(x) -s(x-y) \right) d\Leb(y).
    \end{split}\end{equation*}
    Hence, as multiplication by $\mathcal{C}^r$ functions is continuous on $B^{p,q}_\C$ there exists a $C_0$ such that
    \begin{equation*}\begin{split}
      \int_{W} v_1 \dots v_k&(F_\epsilon R h)(x) \cdot \varphi(x) \, \mathrm{d}x \\
      &= \intf k_\epsilon(y) \left(\int_{W} v_1 \dots v_k(\tau_{-y}\LL_\Omega h \cdot (s - \tau_{-y}s))(x) \cdot \varphi(x) \, \mathrm{d}x \right) d\Leb(y) \\
      &\le \intf \abs{k_\epsilon(y)} \norm{\tau_{-y}\LL_\Omega h \cdot (s - \tau_{-y}s)}_{p,q} d\Leb(y)\\
      &\le C_0\intf \abs{k_\epsilon(y)} \norm{\tau_{-y}\LL_\Omega h}_{p,q} \norm{ s - \tau_{-y}s}_{\mathcal{C}^r} d\Leb(y),
    \end{split}\end{equation*}
    and so
    \begin{equation}\label{eq:Omega_to_Leb_KL_1}
      \norm{F_\epsilon R h}_{p,q} \le C_0\intf \abs{k_\epsilon(y)} \norm{\tau_{-y}\LL_\Omega h}_{p,q} \norm{ s - \tau_{-y}s}_{\mathcal{C}^r} d\Leb(y).
    \end{equation}
    We will bound the right hand side of \eqref{eq:Omega_to_Leb_KL_1} differently depending on whether \ref{en:K1} or \ref{en:K2} holds.

    \paragraph{The case where \ref{en:K1} holds.}
    Recall from \eqref{eq:convolution_kl_1} that $(\tau_{-y}\LL_\Omega h)(x) = (\det D_x \tau_{-y})\cdot(\LL_{\Omega, T_y}h)(x)$ where $\LL_{\Omega, T_y}$ denotes the transfer operator associated to $T_y := T + y$ by duality with respect to $\Omega$.
    If we denote $x\mapsto \det D_x \tau_{-y}$ by $t_y$ then $t_y(x) = \frac{\mathrm{d}\Omega}{\mathrm{d}\Leb}(x + y) \frac{\mathrm{d}\Leb}{\mathrm{d}\Omega}(x)$.
    Since $\frac{\mathrm{d}\Omega}{\mathrm{d}\Leb}, \frac{\mathrm{d}\Leb}{\mathrm{d}\Omega} \in \mathcal{C}^r(\mathbb{T}^d, \R)$ we have
    \begin{equation*}
      \sup_{y \in \mathbb{T}^d}\norm{t_y}_{\mathcal{C}^r} := C_1 < \infty.
    \end{equation*}
    As noted at the beginning of \cite[Section 7]{gouezel2006banach}, there is a $\mathcal{C}^{r+1}(\mathbb{T}^d, \mathbb{T}^d)$ open neighbourhood $U$ of $T$ such that \cite[Lemma 2.2]{gouezel2006banach} applies uniformly to every $S \in U$, and so
    \begin{equation*}
      \sup_{S \in U} \norm{\LL_{\Omega, S}}_{p,q} := C_2 < \infty.
    \end{equation*}
    Hence, by \ref{en:K1}, there exists an $\epsilon' \in (0, \epsilon_0)$ such that $\norm{\LL_{\Omega, T_y}}_{p,q} < C_2$ for every $\epsilon \in (0, \epsilon')$ and $y \in \supp k_\epsilon$.
    Since multiplication by $\mathcal{C}^r$ functions on $B^{p,q}_\C$ \cite[Lemma 3.2]{gouezel2006banach} is continuous there exists a $C_3$ independent of $h$ and $y$ such that
    \begin{equation*}
      \norm{\tau_{-y}\LL_{\Omega}h}_{p,q} = \norm{t_y \cdot (\LL_{\Omega, T_y} h)}_{p,q} \le C_3 \norm{t_y}_{\mathcal{C}^r} \norm{\LL_{\Omega, T_y} h}_{p,q}.
    \end{equation*}
    Setting $C' = C_0 C_1 C_2 C_3$ and applying these estimates to \eqref{eq:Omega_to_Leb_KL_1} yields
    \begin{equation*}
      \norm{F_\epsilon R }_{p,q} \le C' \left(\intf_{\supp k_\epsilon} \abs{k_\epsilon(y)}  d\Leb(y)\right)
      \sup_{y \in \supp k_\epsilon}\norm{ s - \tau_{-y}s}_{\mathcal{C}^r},
    \end{equation*}
    provided that $\epsilon \in (0, \epsilon')$.
    Since $t_y \in \mathcal{C}^\infty(\mathbb{T}^d, \R)$, by \ref{en:K1} we have
    \begin{equation*}
      \lim_{\epsilon \to 0}\sup_{y \in \supp k_\epsilon}\norm{ s - \tau_{-y}s}_{\mathcal{C}^r} =0.
    \end{equation*}
    Recalling that $\{k_\epsilon\}_{\epsilon \in (0, \epsilon_0)}$ is $L^1(\Leb)$-bounded, we therefore have
    \begin{equation*}
    \limsup_{\epsilon \to 0} \norm{F_\epsilon R}_{p,q}
    \le C' \limsup_{\epsilon \to 0} \left(\intf_{\supp k_\epsilon} \abs{k_\epsilon(y)}  d\Leb(y)\right)
      \sup_{y \in \supp k_\epsilon}\norm{ s - \tau_{-y}s}_{\mathcal{C}^r} = 0.
    \end{equation*}
    As $R \in L(B^{p,q}_\C)$ is invertible it follows that $\lim_{\epsilon \to 0} F_\epsilon = 0$ in $L(B^{p,q}_\C)$.
    The same argument can be used to conclude that $\lim_{\epsilon \to 0} F_\epsilon = 0$ in $L(B^{p-1,q+1}_\C)$, and so $\{F_\epsilon\}_{\epsilon \in [0, \epsilon')}$ satisfies the requirements of Proposition \ref{prop:kl_norm_perturb}.

    \paragraph{The case where \ref{en:K2} holds.}
    By \eqref{eq:translate_bound} we have
    \begin{equation*}
      \sup_{y \in \mathbb{T}^d}\max\{\norm{\tau_{y}\LL_\Omega}_{p,q}, \norm{\tau_{y}\LL_\Omega}_{p-1,q+1}\} := C_4 < \infty.
    \end{equation*}
    Applying this to \eqref{eq:Omega_to_Leb_KL_1} yields
    \begin{equation}\label{eq:Omega_to_Leb_KL_2}
        \norm{F_\epsilon R}_{p,q} \le C_0 C_4 \intf \abs{k_\epsilon(y)} \norm{ s - \tau_{-y}s}_{\mathcal{C}^r} d\Leb(y).
    \end{equation}
    Fix $\eta > 0$. By splitting the integral in \eqref{eq:Omega_to_Leb_KL_2} according to the partition $\mathbb{T}^d = B_\eta(0)\cup( \mathbb{T}^d \setminus B_\eta(0))$ we obtain
    \begin{equation*}\begin{split}
      \intf \abs{k_\epsilon(y)} \norm{ s - \tau_{-y}s}_{\mathcal{C}^r} d\Leb(y)
      & \le \norm{k_\epsilon}_{L^1} \sup_{y \in B_\eta(0)} \norm{ s - \tau_{-y}s}_{\mathcal{C}^r} \\
      &+ \sup_{y \not\in B_\eta(0)} \norm{ s - \tau_{-y}s}_{\mathcal{C}^r}\intf_{\mathbb{T}^d \setminus B_\eta(0)} \abs{k_\epsilon(y)} d\Leb(y).
    \end{split}\end{equation*}
    By \ref{en:K2} we have
    \begin{equation}\label{eq:Omega_to_Leb_KL_3}
      \limsup_{\epsilon \to 0}\norm{F_\epsilon R}_{p,q} \le C_0 C_4\left(\sup_{\epsilon \in (0, \epsilon_0)} \norm{k_\epsilon}_{L^1}\right) \left( \sup_{y \in B_\eta(0)} \norm{ s - \tau_{-y}s}_{\mathcal{C}^r}\right),
    \end{equation}
    where the right-hand side is always finite by virtue of the $L^1(\Leb)$ boundedness of $\{k_\epsilon\}_{\epsilon \in (0, \epsilon_0)}$.
    Since $s \in \mathcal{C}^\infty$, letting $\eta \to 0$ in \eqref{eq:Omega_to_Leb_KL_3} yields $\lim_{\epsilon \to 0}\norm{F_\epsilon R}_{p,q} = 0$ in $L(B^{p,q}_\C)$, which implies that $\lim_{\epsilon \to 0}\norm{F_\epsilon}_{p,q} = 0$ by the invertibility of $R$.
    As before, the same argument can be used to conclude that $\lim_{\epsilon \to 0} F_\epsilon = 0$ in $L(B^{p-1,q+1}_\C)$, and so there exists some $\epsilon' \in (0, \epsilon_0)$ such that $\{F_\epsilon\}_{\epsilon \in [0, \epsilon')}$ satisfies the requirements of Proposition \ref{prop:kl_norm_perturb}.
  \end{proof}
\end{proposition}

Using Proposition \ref{prop:Omega_to_Leb_KL} we may now confirm that our spectral stability results for $\LL_\Omega$ from Sections \ref{sec:approximation_anosov} and \ref{sec:fejer_anosov} also apply to $\LL_{\Leb}$.
\begin{proposition}\label{prop:mollifier_scheme_leb}
  If Proposition \ref{prop:convolution_fourier_approx} applies to $\LL_{\Omega}$ then it applies to $\LL_{\Leb}$ too.
  Hence Theorem \ref{thm:twisted_anosov_numeric_stability} holds verbatim if $\LL_{\Omega}$ is replaced by $\LL_{\Leb}$.
  \begin{proof}
    Suppose that $\{q_\epsilon\}_{\epsilon \in (0, \epsilon_0)} \subseteq \mathcal{C}^\infty(\mathbb{T}^d, \R)$ is a family of kernels satisfying \ref{en:S1} and \ref{en:S2}. Recall the definition of $Q_\epsilon$ from the beginning of Section \ref{sec:approximation_anosov}.
    By Lemma \ref{lemma:convolution_kl} there exists $\epsilon_1 \in (0, \epsilon_0)$ so that $\{ Q_\epsilon \LL_{\Omega} \}_{\epsilon \in [0, \epsilon_1)}$ satisfies (KL).
    Since $\{q_\epsilon\}_{\epsilon \in (0, \epsilon_0)}$ is $L^1(\Leb)$-bounded (by \ref{en:S1}) and satisfies \ref{en:S2}, by Proposition \ref{prop:Omega_to_Leb_KL} we may conclude that there exists $\epsilon_2 \in (0, \epsilon_1)$ such that $\{Q_\epsilon \LL_{\Leb}\}_{\epsilon \in [0, \epsilon_2)}$ satisfies (KL) too.
    The proof of Proposition \ref{prop:convolution_fourier_approx} holds verbatim with $\LL_{\Omega}$ replaced with $\LL_{\Leb}$, as does that of Theorem \ref{thm:twisted_anosov_numeric_stability}.
  \end{proof}
\end{proposition}

\begin{proposition}\label{prop:fejer_scheme_leb}
  If Proposition \ref{prop:stochastic_perturbation_kl} applies to $\LL_{\Omega}$ then it applies to $\LL_{\Leb}$ too.
  Hence Theorem \ref{thm:stochastic_perturb_stability} and Corollary \ref{cor:fejer_stability} hold verbatim if $\LL_{\Omega}$ is replaced by $\LL_{\Leb}$.
  \begin{proof}
    Let $\{q_\epsilon\}_{\epsilon \in (0, \epsilon_1)} \subseteq L^1(\Leb)$ be a family of stochastic kernels satisfying \ref{en:S1} and \ref{en:S3}.
    Since $\LL_{\Omega}$ satisfies \eqref{eq:translate_bound}, by Proposition \ref{prop:stochastic_perturbation_kl} there exists some $\epsilon_2 \in (0, \epsilon_1)$ such that $\{q_\epsilon \ast \LL_\Omega \}_{\epsilon \in [0, \epsilon_2)}$ satisfies (KL).
    The family $\{q_\epsilon\}_{\epsilon \in (0, \epsilon_1)}$ is $L^1(\Leb)$-bounded by \ref{en:S1}.
    This, together with \ref{en:S2} and the fact that \eqref{eq:translate_bound} holds for $\LL_{\Omega}$, means that we can apply Proposition \ref{prop:Omega_to_Leb_KL} to conclude that there exists $\epsilon_3 \in (0, \epsilon_2)$ such that $\{q_\epsilon \ast \LL_{\Leb} \}_{\epsilon \in [0, \epsilon_3)}$ satisfies (KL) too.
    The proofs of Theorem \ref{thm:stochastic_perturb_stability} and Corollary \ref{cor:fejer_stability} are the same as before.
  \end{proof}
\end{proposition}

\section{The proof of Proposition \ref{prop:translate_bound_cond}}
\label{app:proof_of_long_prop}

Throughout this appendix we adopt the notation and setting of Section \ref{sec:cond_stochastic_stab_anosov}.
The proof of Proposition \ref{prop:translate_bound_cond} is broken into three parts.
In step 1 we begin by constructing the adapted metric $\langle \cdot, \cdot \rangle$.
In step 2 we define a set of leaves $\tilde{\Sigma}$, which induce spaces $\tilde{B}^{p,q}_\C$ and $\tilde{B}^{p-1,q+1}_\C$, and prove the claim in Proposition \ref{prop:translate_bound_cond} regarding the spectral properties of $\LL$.
Finally, in step 3 we prove that translated leaves may be covered by the image under $T^N$ of finitely many leaves in $\tilde{\Sigma}$ for some large $N$; this result is the core of the proof, and \eqref{eq:translate_bound_cond_1} then easily follows from the adaption of arguments from \cite{gouezel2006banach}.
Steps 2 and 3 lean heavily on the setting in \cite[Section 3]{gouezel2006banach}. We have maintained the notation used in \cite{gouezel2006banach} whenever possible.

\paragraph{Step 1: Constructing the adapted metric $\ip$.}

For sufficiently small $\xi \ge 0$, as $T$ is a $ \mathcal{C}^{r+1}$ diffeomorphism and $\mathbb{T}^d$ is compact, the quantity $D_\xi := \sup_{x \in \mathbb{T}^d} \norm{D_x T^{-1}}_\xi $ is finite. Moreover, as $\norm{\cdot}_\xi \to \norm{\cdot}_0$ uniformly, it follows that $D_\xi \to D_0$.
As $E^s(x)$ and $E^u(x)$ are $\xi$-orthogonal with respect to $\ip_\xi$, one easily verifies that $\norm{\Gamma^s_x \Gamma^u_x}_\xi = \norm{\Gamma^u_x \Gamma^s_x}_\xi < \xi$. What is less obvious, but still true, however, is that
\begin{equation*}
  \norm{\Gamma^s_x \Gamma^u_x}_\xi = \norm{\Gamma^u_x \Gamma^s_x}_\xi = \norm{(\Id - \Gamma_x^s)(\Id - \Gamma^u_x)}_\xi = \norm{(\Id - \Gamma_x^u)(\Id - \Gamma^s_x)}_\xi < \xi.
\end{equation*}
We refer the reader to the proof of Theorem 2 in \cite{buckholtz2000hilbert} for details.
For $\xi > 0$, $x \in \mathbb{T}^d$ and $\kappa \in (0,1)$ we define the stable cone by
\begin{equation*}
  \mathcal{C}_\xi(x, \kappa) = \{ u \in T_x \mathbb{T}^d : \norm{(\Id - \Gamma^s_x) u}_\xi \le \kappa \norm{\Gamma^s_x u}_\xi\}.
\end{equation*}

The following lemma is classical; we reprove it here to emphasise the quantitative estimate \eqref{eq:cone_distortion_3_0}.
\begin{lemma}\label{lemma:cone_distortion_3}
  For every $J > 0$ there exists $\xi_J > 0$ so that for every $\xi \in (0, \xi_J)$, $\gamma \in [0,J)$ and $x \in \mathbb{T}^d$ we have
  \begin{equation}\label{eq:cone_distortion_3_0}
    D_x T^{-1}\mathcal{C}_\xi(x, \gamma) \subseteq \mathcal{C}_\xi\left(T^{-1}x, \gamma \nu_s\nu_u^{-1} \right).
  \end{equation}
  \begin{proof}
    Let $\xi > 0$ and suppose that $u \in \mathcal{C}_\xi(x, \gamma)$.
    As $D_xT^{-1}(E^{s}(x)) = E^{s}(T^{-1}x)$, it follows that
    \begin{equation*}
      (\Id - \Gamma^s_{T^{-1}(x)})D_xT^{-1} u = (\Id - \Gamma^s_{T^{-1}(x)})D_xT^{-1}(\Id - \Gamma^s_{x})u.
    \end{equation*}
    Consequently,
    \begin{equation*}\begin{split}
      \norm{(\Id - \Gamma^s_{T^{-1}(x)})D_xT^{-1} u}_\xi
      &\le \norm{D_xT^{-1}\Gamma^u_x(\Id - \Gamma^s_{x})u}_\xi + \norm{D_xT^{-1}(\Id - \Gamma^u_x)(\Id - \Gamma^s_{x})u}_\xi.
    \end{split}\end{equation*}
    As $\norm{(\Id - \Gamma_x^u)(\Id - \Gamma^s_x)}_\xi < \xi$ and $u \in \mathcal{C}_\xi(x, \gamma)$, it follows that
    \begin{equation*}\begin{split}
      \norm{(\Id - \Gamma^s_{T^{-1}(x)})D_xT^{-1} u}_\xi &\le
      \left(\norm{\restr{D_xT^{-1}}{E^u(x)}}_\xi + \xi \norm{D_x T^{-1}}_\xi \right)\norm{(\Id - \Gamma^s_{x}) u}_\xi\\
      &\le \gamma\left(\norm{\restr{D_xT^{-1}}{E^u(x)}}_\xi + \xi \norm{D_x T^{-1}}_\xi \right)\norm{\Gamma_x^s u}_\xi.
    \end{split}\end{equation*}
    Similarly,
    \begin{equation*}\begin{split}
      \norm{\Gamma_{T^{-1}(x)}^s D_xT^{-1}u}_\xi \ge &\norm{\Gamma_{T^{-1}(x)}^s D_xT^{-1} \Gamma^s_x u}_\xi - \norm{\Gamma_{T^{-1}(x)}^s D_xT^{-1}(\Id -\Gamma_x^{u})(\Id - \Gamma^s_{x}) u}_\xi\\
      &-\norm{\Gamma_{T^{-1}(x)}^s \Gamma_{T^{-1}(x)}^{u}D_xT^{-1}(\Id - \Gamma^s_{x}) u}_\xi\\
      \ge &\left(\norm{ \restr{D_x T}{E^{s}(x)} }_\xi^{-1} - 2\xi \gamma\norm{D_x T^{-1}}_\xi \right)\norm{\Gamma^s_x u}_\xi.
    \end{split}\end{equation*}
    For sufficiently small $\xi$ we have  $\norm{\restr{D_xT^{-1}}{E^u(x)}}_\xi < \nu_u^{-1}$ and $\norm{\restr{D_xT}{E^s(x)}}_\xi < \nu_s$.
    Hence, as $ \sup_{x \in \mathbb{T}^d} \norm{D_xT^{-1}}_\xi = D_\xi \to D_0$, there exists $\xi_J > 0$ so that for every $\xi \in (0, \xi_J)$ we have
    \begin{equation*}
      \norm{\restr{D_xT^{-1}}{E^u(x)}}_\xi + \xi \norm{D_x T^{-1}}_\xi \le \nu_u^{-1},
    \end{equation*}
    and, for all $\gamma < J$,
    \begin{equation*}
      \norm{\restr{D_xT}{E^s(x)}}_\xi^{-1} - 2\xi \gamma\norm{D_xT^{-1}}_\xi \ge \nu_s^{-1}.
    \end{equation*}
    In view of the above, whenever $\xi \in (0, \xi_J)$ we therefore have
    \begin{equation*}
      \norm{(\Id - \Gamma^s_{T^{-1}(x)})D_xT^{-1} u}_\xi \le \gamma\nu_s\nu_u^{-1} \norm{\Gamma_{T^{-1}(x)}^s D_xT^{-1}u}_\xi
    \end{equation*}
    for every $u \in \mathcal{C}_\xi(\gamma, x)$. Thus $D_x T^{-1}\mathcal{C}_\xi(x, \gamma) \subseteq \mathcal{C}_\xi\left(T^{-1}x, \gamma \nu_s\nu_u^{-1} \right)$ for every $\xi \in (0, \xi_J)$, $x \in \mathbb{T}^d$ and $\gamma < J$, as required.
  \end{proof}
\end{lemma}

We aim to select $\xi$ so that $C_{\tau,\xi}^{-1} > \Theta_{T,\xi}$, and so that we can apply Lemma \ref{lemma:cone_distortion_3} for some appropriate $J$.
As the adapted metric $\ip_\xi$ uniformly approximates $\ip_0$, we have $\lim_{\xi \to 0} C_{\tau,\xi}^{-1} = C_{\tau,0}^{-1}$. However, upon examining the definition of $\Theta_{T,\xi}$ we observe that the projections $\Gamma^s_x$ depend on $\ip_\xi$, and so the behaviour of $\Theta_{T,\xi}$ as $\xi \to 0$ is not clear.
We address this now.

\begin{lemma}\label{lemma:translate_commutator_bound}
  $\lim_{\xi \to 0} \Theta_{T,\xi} = \Theta_{T,0}$.
  \begin{proof}
    Let $\pi^s_x$ and $\pi_x^u$ be the projections induced by the direct sum $T_x \mathbb{T}^d = E^s(x) \oplus E^u(x)$. We have
    \begin{equation*}
      \norm{\Gamma^s_{x} - \pi^s_x}_\xi = \norm{(\Gamma^s_{x} - \pi^s_x)\pi^s_x}_\xi + \norm{(\Gamma^s_{x} - \pi^s_x)\pi^u_x}_\xi = \norm{(\Gamma^s_{x} - \id)\pi^s_x}_\xi + \norm{\Gamma^s_x \pi^u_x}_\xi.
    \end{equation*}
    Since $\Gamma^s_{x} = \Id$ on $E^s(x)$, $\norm{(\Gamma^s_{x} - \id)\pi^s_x}_\xi = 0$. Let $v \in T_x\mathbb{T}^d$. As $E^s(x)$ and $E^u(x)$ are $\xi$-orthogonal and $\Gamma^s_x$ is an orthogonal projection (both with respect to $\ip_\xi$), we have
    \begin{equation*}
        \norm{\Gamma^s_x \pi^u_x v}_\xi = \sqrt{\abs{\langle \Gamma^s_x \pi^u_x v, \pi^u_x v\rangle_\xi}} \le \sqrt{\xi \norm{\Gamma^s_x \pi^u_x v}_\xi \norm{\pi^u_x v}_\xi} \le \sqrt{\xi} \norm{\pi^u_x}_\xi.
    \end{equation*}
    Thus $\norm{\Gamma^s_{x} - \pi^s_x}_\xi \le \sqrt{\xi} \norm{\pi^u_x}_\xi$.
    Let $P_\xi = \sup_{x \in \mathbb{T}^d} \norm{\pi^u_x}_\xi$.
    For any $x, y \in \mathbb{T}^d$ the triangle inequality yields
    \begin{equation*}\begin{split}
       \abs{\norm{\Gamma^s_{x+y} D_x \tau_y - (D_x \tau_y) \Gamma^s_{x}}_\xi  - \norm{\pi^s_{x+y} D_x \tau_y - (D_x \tau_y) \pi^s_{x}}_\xi} \le &\norm{(\Gamma^s_{x+y} - \pi^s_{x+y})( D_x \tau_y)}_\xi \\
       &+ \norm{(D_x \tau_y) (\Gamma^s_{x} - \pi^s_{x})}_\xi \\
       &\le 2 \sqrt{\xi} C_{\tau,\xi} P_\xi.
    \end{split}\end{equation*}
    It follows that
    \begin{equation}\label{eq:translate_commutator_bound_1}
       \abs{\Theta_{T,\xi}  - \sup_{x,y \in \mathbb{T}^d} \norm{\pi^s_{x+y} D_x \tau_y - (D_x \tau_y) \pi^s_{x}}_\xi} \le 2 \sqrt{\xi} C_{\tau,\xi} P_\xi.
    \end{equation}
    As $\pi^u_x$ is independent of $\xi$, we have $P_\xi \to P_0$.
    Since $E^s(x) \perp E^u(x)$ with respect to $\ip_0$, the projections $\pi^u_x$ and $\pi^s_x$ are orthogonal with respect to $\ip_0$. Thus the uniform convergence of $\ip_\xi$ to $\ip_0$ implies that
    \begin{equation}\label{eq:translate_commutator_bound_2}
      \lim_{\xi \to 0} \sup_{x,y \in \mathbb{T}^d}\norm{\pi^s_{x+y} D_x \tau_y - (D_x \tau_y) \pi^s_{x}}_\xi = \sup_{x,y \in \mathbb{T}^d}\norm{\pi^s_{x+y} D_x \tau_y - (D_x \tau_y) \pi^s_{x}}_0 = \Theta_{T,0}.
    \end{equation}
    Hence, as $C_{\tau, \xi} \to C_{\tau,0}$, by letting $\xi \to 0$ in \eqref{eq:translate_commutator_bound_1} and applying \eqref{eq:translate_commutator_bound_2} we have
    \begin{equation*}
      \lim_{\xi \to 0} \Theta_{T,\xi} = \lim_{\xi \to 0} \sup_{x,y \in \mathbb{T}^d} \norm{\pi^s_{x+y} D_x \tau_y - (D_x \tau_y) \pi^s_{x}}_\xi = \Theta_{T,0},
    \end{equation*}
    as required.
  \end{proof}
\end{lemma}

We now fix, once and for all, the metric that is used in Proposition \ref{prop:translate_bound_cond}.
As $C_{\tau,0}^{-1} > \Theta_{T,0}$, $C_{\tau, \xi} \to C_{\tau,0}$, and, by Lemma \ref{lemma:translate_commutator_bound}, $\Theta_{T,\xi} \to \Theta_{T,0}$, there exists $\mathcal{E} > 0$ so that
\begin{equation}\label{eq:translate_bound_cond_2}
  \inf_{\xi \in [0, \mathcal{E}]} \left( C_{\tau,\xi}^{-1} - \Theta_{T,\xi}\right) > 0, \quad \text{and} \quad \sup_{\xi \in [0, \mathcal{E}]} \Theta_{T, \xi} < \infty.
\end{equation}
We apply Lemma \ref{lemma:cone_distortion_3} with
\begin{equation}\label{eq:translate_bound_cond_3}
  J = 1 + \frac{\sup_{\xi \in  [0, \mathcal{E}]} \Theta_{T, \xi}}{\inf_{\xi \in  [0, \mathcal{E}]} \left(C_{\tau,\xi}^{-1} - \Theta_{T,\xi}\right)}
\end{equation}
to produce an adapted metric $\ip := \ip_\mathcal{E}$, which may replace the metric defined in \cite[Section 3]{gouezel2006banach} after possibly shrinking $\mathcal{E}$ further.

\begin{remark}
  Until the end of this section we only deal with the metric just constructed, and so we drop references to $\xi$ and $\mathcal{E}$ from our notation.
\end{remark}

\paragraph{Step 2: Defining the set of leaves $\tilde{\Sigma}$ and the spaces $\tilde{B}^{p,q}_\C$ and $\tilde{B}^{p-1,q+1}_\C$.}

Our task is now to define a set of leaves $\tilde{\Sigma}$ and spaces $\tilde{B}^{p,q}_\C$ and $\tilde{B}^{p-1,q+1}_\C$ so that the spectral properties of $\LL$ on $B^{p,q}_\C$ and $\tilde{B}^{p,q}_\C$ are identical.
It is necessary to understand how leaves are defined; to this end we reproduce material from the beginning of \cite[Section 3]{gouezel2006banach}. After fixing the metric, in \cite{gouezel2006banach} a small $\kappa > 0$ satisfying various properties is fixed; in particular it is required that $D_xT^{-1}$ expands the vectors in $\mathcal{C}(x,\kappa)$ by at least $\nu_s^{-1}$.
We now choose a smaller value of $\kappa$, as follows. Let $C_\tau = C_{\tau, \mathcal{E}}$ and $\Theta_T = \Theta_{T, \mathcal{E}}$. By  \eqref{eq:translate_bound_cond_2} and \eqref{eq:translate_bound_cond_3}, there exists $\kappa' \le \min\{\kappa, 1/2\}$, $\eta > 1$ such that
\begin{equation}\label{eq:cone_distortion}
  C_\tau^{-1} > \Theta_T( 1+ 2\kappa' \eta), \quad \text{and} \quad  \frac{\Theta_T + (C_\tau + \Theta_T) 2\kappa' \eta}{C_\tau^{-1} - \Theta_T(1 + 2\kappa' \eta)} < J.
\end{equation}
We redefine the original $\kappa$ to be $\kappa'$, noting that this does not alter the validity of any arguments in \cite{gouezel2006banach} (they only require that $\kappa$ be sufficiently small, so we are free to make it as small as we require).
Two components of this construction appear, at first glance, arbitrary: the constant $\eta$ and the inequalities \eqref{eq:cone_distortion}. They appear so that later we may cover translated leaves by the image under some iterate of $T$ by finitely many leaves (Lemma \ref{lemma:covering_leaves}).

As in \cite{gouezel2006banach}, one may construct finitely many $\mathcal{C}^\infty$ charts $\psi_1, \dots, \psi_S$, each respectively defined on $(-r_i, r_i)^d \subseteq \R^d$, so that
\begin{enumerate}[label=(B\arabic*)]
  \item \label{en:B1} $D_0 \psi_i$ is an isometry;
  \item $(D_0 \psi_i)(\R^{d_s} \times \{0\}) = E^s(\psi_i(0))$;
  \item The $\mathcal{C}^{r+1}$ norms of $\psi_i$ and $\psi_i^{-1}$ are bounded by $1 + \kappa$;
  \item \label{en:B4} There exists $c_i \in (\kappa, 2 \kappa)$ such that the cone
  \begin{equation*}
    \mathcal{C}_i = \{ u + v \in \R^d \,| \,u \in \R^{d_s} \times \{0\}, v \in \{0\} \times \R^{d_u}, \norm{v} \le c_i \norm{u} \}
  \end{equation*}
  satisfies the following property: for any $x \in (-r_i, r_i)^d$, $\mathcal{C}(\psi_i(x)) \subseteq (D_x \psi_i)\mathcal{C}_i$ and $(D_{\psi_i(x)}T^{-1})( D_x \psi_i) \mathcal{C}_i \subseteq \mathcal{C}(T^{-1} (\psi_i(x)))$; and
  \item \label{en:B5} $\mathbb{T}^d$ is covered by $\{\psi_i((-r_i/2, r_i/2)^d) \}_{i=1}^S$.
\end{enumerate}

We require that the charts satisfy following additional property concerning the distortion of the cones $\mathcal{C}_i$ under $D_x\psi$ in terms of the stable cones.

\begin{lemma}\label{lemma:cone_distortion_1}
  There exist charts $\{\psi_i\}_{i=1}^S$ satisfying \ref{en:B1}-\ref{en:B5} and so that for any $x \in (-r_i, r_i)^d$ we have
  \begin{equation*}
    (D_x\psi_i) \mathcal{C}_i \subseteq \mathcal{C}\left(\psi_i(x),  \eta c_i\right),
  \end{equation*}
  where $\eta > 1$ is the constant appearing in \eqref{eq:cone_distortion}.
  \begin{proof}
    By compactness it is sufficient to construct for each $y \in \mathbb{T}^d$ a chart $\psi_y : (-r_y,r_y)^d \to \mathbb{T}^d$ that satisfies all of the given requirements and for which $\psi_y(0) = y$.
    As noted immediately before the statement of \cite[Lemma 3.1]{gouezel2006banach}, for each $y \in \mathbb{T}^d$ one can construct a $\psi_y$ satisfying conditions \ref{en:B1}-\ref{en:B5} and so that $\psi_y(0) = y$. Let $\mathcal{C}_y$ denote the corresponding cone and $c_y \in (\kappa, 2 \kappa)$ denote the constant corresponding to $c_i$.
    Let $u +v \in \mathcal{C}_y$ so that $u \in \R^{d_s} \times \{0\}$ and $v \in \{0\} \times \R^{d_u}$. Then
    \begin{equation*}
      \norm{(\Id - \Gamma_y^s)D_0 \psi_y (u+v)} = \norm{(\Id - \Gamma_y^s)D_0 \psi_y v} \le \norm{v} \le c_y \norm{u},
    \end{equation*}
    and
    \begin{equation*}
      \norm{\Gamma_y^s D_0 \psi_y (u+v)} = \norm{\Gamma_y^s D_0 \psi_y u} = \norm{u}.
    \end{equation*}
    Together implying that $(D_0\psi_i)\mathcal{C}_y \subseteq \mathcal{C}\left(y, c_y\right)$.
    Hence, using the compactness of the closed unit ball in $\R^d$ and the uniform continuity of $\psi_y$ on $(-r_y,r_y)^d$, $x \mapsto \Gamma_x^s$ on $\mathbb{T}^d$ and $u \mapsto \norm{u}$ on both the tangent space of $\mathbb{T}^d$ and $(-r_y, r_y)^d$, we may shrink $r_y$ so that $D_x\psi_i \mathcal{C}_y \subseteq \mathcal{C}\left(\psi_y(x),  \eta c_y\right)$ for every $x \in (-r_y, r_y)^d$.
    Thus $\psi_y$ satisfies all the requirements of the lemma, and we may conclude using the compactness of $\mathbb{T}^d$.
  \end{proof}
\end{lemma}

We are now able to construct our modified set of leaves $\tilde{\Sigma}$ and spaces $\tilde{B}^{p,q}_\C$ and $\tilde{B}^{p-1,q+1}_\C$, and prove the claim in Proposition \ref{prop:translate_bound_cond} regarding the spectral properties of $\LL$.
Using the metric $\ip$ constructed in step 1 of the proof of Proposition \ref{prop:translate_bound_cond}, the constant $\kappa$ defined immediately following \eqref{eq:cone_distortion}, and the charts from Lemma \ref{lemma:cone_distortion_1}, we may define the set of leaves $\tilde{\Sigma}$ exactly as in \cite[Section 3]{gouezel2006banach}.
Recall that $p \in \Z^+$ and $q > 0$ satisfy $p + q < r$.
In exactly the same way that the set of leaves $\Sigma$ from \cite{gouezel2006banach} induces spaces $B^{p,q}_\C$ and $B^{p-1,q+1}_\C$ (see \cite[Section 3]{gouezel2006banach}), our set of leaves $\tilde{\Sigma}$ induces spaces $\tilde{B}^{p,q}_\C$ and $\tilde{B}^{p-1,q+1}_\C$.
The proofs of \cite[Lemma 2.2]{gouezel2006banach} and \cite[Theorem 2.3]{gouezel2006banach} hold verbatim for $\LL$ on $\tilde{B}^{p,q}_\C$.
Thus the essential spectral radius of $\LL$ on $\tilde{B}^{p,q}_\C$ is bounded by $\max\{\nu_u^{-p}, \nu_s^{q} \}$.
As per \cite[Remark 2.5]{gouezel2006banach}, the spectral data of $\LL$ associated to eigenvalues outside of the ball of radius $\max\{\nu_u^{-p}, \nu_s^{q} \}$ are the same on $B^{p,q}_\C$ and $\tilde{B}^{p,q}_\C$, and in particular the generalised eigenspaces of all such eigenvalues lie in $B^{p,q}_\C \cap \tilde{B}^{p,q}_\C$.

\paragraph{Step 3: Obtaining the inequality \eqref{eq:translate_bound_cond_1}.}

Key to establishing \eqref{eq:translate_bound_cond_1} is the following lemma, which extends the result from \cite[Lemma 3.3]{gouezel2006banach} to include translated leaves.
Throughout this step of the proof we assume the reader is familiar with the definition of the set of leaves $\tilde{\Sigma}$ from \cite[Section 3]{gouezel2006banach}.

\begin{lemma}\label{lemma:covering_leaves}
  There exists $N,M,C > 0$ such that for any $W \in \tilde{\Sigma}$, with associated full admissible leaf $\tilde{W}$, and $y \in \mathbb{T}^d$ there exists $\{W_s\}_{s=1}^m \subseteq \tilde{\Sigma}$ with $m < M$ so that
  \begin{enumerate}
    \item $T^{-N}(W + y) \subseteq \bigcup_{s=1}^m W_s \subseteq T^{-N}(\tilde{W} + y)$.
    \item There are $\mathcal{C}^{r+1}$ functions $\{\rho_s\}_{i=1}^m$ so that each $\rho_s$ is compactly supported on $W_s$, $\sum_s \rho_s = 1$ on $T^{-N}(W + y)$, and $\norm{\rho_s}_{\mathcal{C}^{r+1}} \le C$.
  \end{enumerate}
\end{lemma}

We will now give a brief, non-technical overview of the strategy for proving Lemma \ref{lemma:covering_leaves}. Most of the proof is dedicated to finding leaves $W_s$ which verify the containment
\begin{equation*}
  T^{-N}(W + y) \subseteq \bigcup_{s=1}^m W_s \subseteq T^{-N}(\tilde{W} + y).
\end{equation*}
Let us consider the case where $y = 0$, which is the subject of \cite[Lemma 3.3]{gouezel2006banach}. The idea in this case is to use the expansion and regularisation of $T^{-1}$ in the stable direction to prove that $T^{-1}(\tilde{W})$ is locally `leaf-like'. One then picks certain subsets of $T^{-1}(\tilde{W})$, and proves they are leaves that cover $T^{-1}(W)$.
The main issue resulting from translation by a non-zero $y$ is the possibility that the translated leaf will not be `leaf-like'. Specifically, while all leaves in $\tilde{\Sigma}$ lie approximately parallel to the stable manifold, the translate of a leaf may well lie very near an unstable manifold.
Our main hypothesis, the inequality \eqref{eq:almost_flat_constants}, will imply that translated leaves do not lie too close to the unstable manifold, which allows for the aforementioned distortion to be corrected using the regularisation of $T^{-n}$ for some large $n$.

We require some preliminary lemmas before proving Lemma \ref{lemma:covering_leaves}. The following gives an estimate of the distortion that stable cones experience under translation, and is where the hypothesis that $C_\tau \Theta_T < 1$ is crucial.

\begin{lemma}\label{lemma:cone_distortion_2}
  If $ \gamma > 0$ satisfies $C_\tau^{-1}  > \Theta_T(1 + \gamma)$ then for each $x,y \in \mathbb{T}^d$ we have
  \begin{equation*}
    (D_x\tau_{y}) \mathcal{C}(x, \gamma) \subseteq \mathcal{C}\left(x + y, \frac{C_\tau \gamma + \Theta_{T}(1 + \gamma)}{C_\tau^{-1} - \Theta_{T}(1 + \gamma)}\right).
  \end{equation*}
  \begin{proof}
    Suppose that $u \in \mathcal{C}(x, \gamma)$. We have
    \begin{equation}\begin{split}\label{eq:cone_distortion_2_1}
      \norm{(\Id - \Gamma^s_{x+y}) (D_x\tau_{y})u} \le \norm{(\Id - \Gamma^s_{x+y}) (D_x\tau_{y})\Gamma^s_x u} + \norm{(\Id - \Gamma^s_{x+y}) (D_x\tau_{y})(\Id -\Gamma^s_x)u}.
    \end{split}\end{equation}
    We estimate the first term on the right-hand side of \eqref{eq:cone_distortion_2_1}. By the triangle inequality we have
    \begin{equation*}
      \norm{(\Id - \Gamma^s_{x+y}) (D_x\tau_{y})\Gamma^s_x u} \le \norm{(\Id - \Gamma^s_{x+y}) \Gamma^s_{x+y} (D_x\tau_{y}) u} + \norm{(\Id - \Gamma^s_{x+y}) (\Gamma^s_{x+y} D_x\tau_{y}  -  D_x\tau_{y} \Gamma^s_x) u}.
    \end{equation*}
    The first term on the right hand side is 0, whereas the second term may be estimated using the definition of $\Theta_T$ (see \eqref{eq:almost_flat_constants}), yielding
    \begin{equation*}
      \norm{(\Id - \Gamma^s_{x+y}) (D_x\tau_{y})\Gamma^s_x u} \le \Theta_T \norm{u}.
    \end{equation*}
    As $u \in \mathcal{C}(x, \gamma)$, we have
    \begin{equation} \label{eq:cone_distortion_2_2}
      \norm{(\Id - \Gamma^s_{x+y}) (D_x\tau_{y})\Gamma^s_x u} \le \Theta_T(1 + \gamma) \norm{\Gamma_x^s u}.
    \end{equation}
    We turn to estimating the second term on the right-hand side of \eqref{eq:cone_distortion_2_1}. Using the definition of $C_\tau$ and as $u \in \mathcal{C}(x, \gamma)$, we have
    \begin{equation}\label{eq:cone_distortion_2_3}
      \norm{(\Id - \Gamma^s_{x+y}) D_x\tau_{y}(\Id -\Gamma^s_x)u} \le \norm{D_x\tau_{y}}\norm{(\Id -\Gamma^s_x)u} \le C_\tau \gamma \norm{\Gamma^s_x u}.
    \end{equation}
    Applying \eqref{eq:cone_distortion_2_2} and \eqref{eq:cone_distortion_2_3} to \eqref{eq:cone_distortion_2_1} yields
    \begin{equation}\label{eq:cone_distortion_2_4}
      \norm{(\Id - \Gamma^s_{x+y}) (D_x\tau_{y})u} \le \left(\Theta_T(1 + \gamma) + C_\tau \gamma \right)\norm{\Gamma^s_x u}.
    \end{equation}
    Alternatively, the reverse triangle inequality yields
    \begin{equation*}
      \norm{\Gamma^s_{x+y} D_x\tau_{y}u} \ge \norm{D_x\tau_{y} \Gamma^s_x u} - \norm{(\Gamma^s_{x+y} D_x\tau_{y} - D_x\tau_{y} \Gamma^s_x )u}
    \end{equation*}
    Using a similar process as in the estimation of \eqref{eq:cone_distortion_2_4}, we obtain
    \begin{equation}\label{eq:cone_distortion_2_5}
      \norm{\Gamma^s_{x+y} (D_x\tau_{y})u} \ge \norm{D_x\tau_{y} \Gamma^s_x u} - \Theta_T \norm{u} \ge (C_\tau^{-1} - \Theta_T(1 + \gamma)) \norm{\Gamma^s_x u}.
    \end{equation}
    By our assumptions $C_\tau^{-1}  > \Theta_T(1 + \gamma)$, and so we may combine \eqref{eq:cone_distortion_2_4} and \eqref{eq:cone_distortion_2_5} to obtain
    \begin{equation*}
    \norm{(\Id - \Gamma^s_{x+y}) (D_x\tau_{y})u} \le \frac{C_\tau \gamma + \Theta_T(1 + \gamma)}{C_\tau^{-1} - \Theta_T(1 + \gamma)} \norm{\Gamma^s_{x+y} D_x\tau_{y}u},
    \end{equation*}
    as required.
  \end{proof}
\end{lemma}

In the following lemma we show that the distortion of the stable cones experience under translation may be corrected by applying $T^{-n}$ for $n$ large.

\begin{lemma}\label{lemma:distorted_leaf_tangent_space}
  Recall the cones $\mathcal{C}_i$ from \ref{en:B4}. There exists $N_1 > 0$ such that for every $x, y \in \mathbb{T}^d$, where $\psi_i^{-1}(x) \in (-r_i,r_i)^d$, we have
  \begin{equation*}
      (D_{x+y}T^{-N_1} )(D_x\tau_y)( D_{\psi_i^{-1}(x)} \psi_i)\mathcal{C}_i \subseteq \mathcal{C}\left(T^{-N_1}(x+y), \kappa \right).
  \end{equation*}
  Moreover, if $\psi_j^{-1}(T^{-N_1}(x+y)) \in (-r_j,r_j)^d$ then
  \begin{equation*}
    (D_{T^{-N_1}(x+y)}\psi_j^{-1})( D_{x+y}T^{-N_1} )(D_x\tau_y) (D_{\psi_i^{-1}(x)} \psi_i)\mathcal{C}_i \subseteq \mathcal{C}_j.
  \end{equation*}
  \begin{proof}
    Recall from \ref{en:B4} that $c_i \le 2\kappa$. Lemma \ref{lemma:cone_distortion_1} implies that $(D_{\psi_i^{-1}(x)}\psi_i) \mathcal{C}_i \subseteq \mathcal{C}\left(x, 2\eta \kappa \right)$.
    By \eqref{eq:cone_distortion} we have $C_\tau^{-1}  > \Theta_T(1 + 2\eta \kappa)$, and so Lemma \ref{lemma:cone_distortion_2} yields
    \begin{equation*}
      (D_x\tau_y )(D_{\psi_i^{-1}(x)} \psi_i)\mathcal{C}_i
      \subseteq \mathcal{C}\left(x+y, \frac{\Theta_T + (C_\tau + \Theta_T) 2\kappa \eta}{C_\tau^{-1} - \Theta_T - \Theta_T 2\kappa \eta}\right).
    \end{equation*}
    Let $N_1 \in \Z^+$ be large enough so that
    \begin{equation*}
      \nu_s^{N_1}\nu_u^{-N_1}\frac{\Theta_T + (C_\tau + \Theta_T) 2\kappa \eta}{C_\tau^{-1} - \Theta_T - \Theta_T 2\kappa \eta} \le \kappa.
    \end{equation*}
    By the definition of the adapted metric $\ip$ from step 1 of the proof of Proposition \ref{prop:translate_bound_cond}, the second inequality in \eqref{eq:cone_distortion}, and Lemma \ref{lemma:cone_distortion_3}, it follows that
    \begin{equation*}
      (D_{x+y}T^{-N_1}) (D_x\tau_y) (D_{\psi_i^{-1}(x)} \psi_i)\mathcal{C}_i \subseteq \mathcal{C}\left(T^{-N_1}(x+y), \kappa \right).
    \end{equation*}
    Since all our estimates are uniform in $x$ and $y$ we obtain the first claim. The second claim follows from \ref{en:B4}.
  \end{proof}
\end{lemma}

Recall that $T^{-1}$ is expansive along leaves in $\tilde{\Sigma}$, since all leaves are approximately parallel to the stable direction of $T$. In the proof of Lemma \ref{lemma:covering_leaves} we will require that a version of this property holds for translated leaves as well.
Up until now we have considered how translation affects the stable cones, and how applying $T^{-1}$ corrects for any distortion in the cones.
In the following lemma we apply the same idea to show that $T^{-n}$ is expansive along translated leaves provided that $n$ is sufficiently large.

\begin{lemma}\label{lemma:translated_leaf_expansion}
  Let $H = \inf_{x \in \mathbb{T}^d} \norm{D_x T}$ and $N_1$ be from Lemma \ref{lemma:distorted_leaf_tangent_space}.
  If $n > N_1$ and $\nu_s^{-n + N_1} H^{-N_1} > 1$ then for any $W \in \tilde{\Sigma}$ and $y \in \mathbb{T}^d$ the map $T^{-n}$ expands distances on $\tilde{W} +y$ by at least $\nu_s^{-n + N_1} H^{-N_1}$.
  \begin{proof}
    Let $\psi_i$ be a chart whose image contains $\tilde{W}$ and for which the tangent space of $\psi_i^{-1}(\tilde{W})$ is contained in $\mathcal{C}_i$.
    Suppose $a,b \in \tilde{W} + y$ and that $\gamma : [0,1] \to T^{-n}(\tilde{W} + y)$ is a distance minimizing geodesic from $T^{-n}(a)$ to $T^{-n}(b)$. Define $\gamma_n := T^n \circ \gamma$ and note that $\gamma_n$ is a differentiable curve from $a$ to $b$ lying in $\tilde{W} + y$. For $n > N_1$ we have
    \begin{equation}\begin{split}
      \label{eq:translated_leaf_expansion_1}
      d_{T^{-n}(\tilde{W}+y)}(T^{-n}(a),T^{-n}(b)) &= \intf_{0}^1 \norm{D_t \gamma} dt \\
      &= \intf_{0}^1 \norm{(D_{(T^{-N_1}\circ \gamma_n)(t)} T^{-n+N_1}) (D_{\gamma_n(t)}T^{-N_1} )(D_t \gamma_n)} dt,
    \end{split}\end{equation}
    where $N_1$ is the constant from Lemma \ref{lemma:distorted_leaf_tangent_space}.
    Since the image of $\gamma_n$ is a closed sub-manifold of $\tilde{W} + y$ and the tangent space of $\tilde{W}$ at $w$ is contained in $(D_{\psi^{-1}(w)}\psi)\mathcal{C}_i$, the image of $D_t \gamma_n$ is contained in $(D_{ \gamma_n(t) - y}\tau_y)( D_{(\psi_i^{-1}\circ \tau_{-y} \circ \gamma_n)(t)} \psi_i)\mathcal{C}_i$.
    Thus, by Lemma \ref{lemma:distorted_leaf_tangent_space} we have $(D_{\gamma_n(t)}T^{-N_1} )(D_t \gamma_n) \subseteq \mathcal{C}(T^{-N_1}(\gamma_n(t)), \kappa)$.
    As $DT^{-n}$ expands vectors in stable cones by at least $\nu_s^{-n}$ we may bound \eqref{eq:translated_leaf_expansion_1} as follows
    \begin{equation*}\begin{split}
      d_{T^{-n}(\tilde{W}+y)}(T^{-n}(a),T^{-n}(b))
      &\ge \nu_s^{-n + N_1} \intf_{0}^1 \norm{ (D_{\gamma_n(t)}T^{-N_1} )(D_t \gamma_n)} dt\\
      &\ge \nu_s^{-n + N_1} H^{-N_1} \intf_{0}^1 \norm{ D_t \gamma_n} dt \ge \nu_s^{-n + N_1} H^{-N_1} d_{\tilde{W} + y}(a,b).
    \end{split}\end{equation*}
    Hence $T^{-n}$ expands distances in $\tilde{W}+y$ by a factor of at least $\nu_s^{-n + N_1} H^{-N_1}$ provided that $n > N_1$ and $\nu_s^{-n + N_1} H^{-N_1} > 1$.
  \end{proof}
\end{lemma}

The following lemma quantifies the regularisation that leaves experience under $T^{-1}$, and is a strengthening of \cite[Lemma 3.1]{gouezel2006banach}. Whereas the previous results were concerned with the regularisation of the first derivative of the leaves (via the contraction of stable cones), the forthcoming result concerns the regularisation of the higher derivatives of leaves.

\begin{lemma}\label{lemma:leaf_regularisation}
  For $L > 0$ and $i = 1,\dots,S$ let $G_i(L)$ be the set defined immediately before \cite[Lemma 3.1]{gouezel2006banach}, and let
  \begin{equation*}
    R(L) := \inf \{ L' : (\psi_{j}^{-1}\circ T^{-1}\circ\psi_i)(W) \in G_j(L') \text{ for every } W \in G_i(L) \text{ and } i,j = 1, \dots, S \}.
  \end{equation*}
  For every $K$ sufficiently large the following holds: after possibly refining the charts $\{\psi_i\}_{i=1}^S$ from Lemma \ref{lemma:cone_distortion_1}, for each $L > 0$ there exists $N(L) \in \Z^+$ so that for each $n \ge N(L)$ we have $R^{n}(L) \le K$.
  \begin{proof}
    The finiteness of $R(L)$ and the fact that $R(L) < L$ for $L$ sufficiently large follow from \cite[Lemma 3.1]{gouezel2006banach}.
    Our more general claim is a classical consequence of the uniform hyperbolicity of $T$ and the regularisation of the associated graph transform, so we will only sketch the ingredients of the proof.

    Suppose that $W \in G_{i}(L)$ is the graph of $\chi: \cl{B(x, A\delta)} \to (-r_i,r_i)^{d_u}$.
    As outlined at the beginning of \cite[Section 6.4.b]{katok1997introduction}, using the exponential map, \cite[Lemma 6.2.7]{katok1997introduction} and after possibly refining the set of charts $\{\psi_i\}_{i=1}^S$ so that each $r_i$ is sufficiently small, one may apply the arguments from \cite[Theorem 6.2.8]{katok1997introduction} (refer to steps 3 and 4 of the proof of \cite[Theorem 6.2.8]{katok1997introduction} for context, and to step 5 for the relevant argument) to conclude that $(\psi_{j}^{-1}\circ T^{-1}\circ\psi_i)(W)$ is the graph of some map $\chi': U \subseteq (-r_j,r_j)^{d_s} \to (-r_j,r_j)^{d_u}$.
    Due to the uniform convergence of the graph transform as outlined in step 5 of the proof of \cite[Theorem 6.2.8]{katok1997introduction}, there exists some $L' < L$ so that $\norm{\chi'}_{\mathcal{C}^{r+1}} \le L'$ for every such $\chi$ provided that $L$ is sufficiently large (i.e. bigger than $K$). Hence $R(L)$ exists and satisfies $R(L) < L$ for $L$ large enough.
    Further examining the proof of \cite[Theorem 6.2.8]{katok1997introduction} yields the stronger claim that for sufficiently large $K$ we have
    \begin{equation*}
      \sup_{L \ge K} \frac{R(L)}{L} < 1,
    \end{equation*}
    which immediately yields the required claim.
  \end{proof}
\end{lemma}

Before we prove Lemma \ref{lemma:covering_leaves} we recall a quantitative version of the inverse function theorem.

\begin{lemma}[{\cite[XIV \S1 Lemma 1.3]{lang2012real}}]\label{lemma:inverse_func}
  Let $E$ be a Banach space, $U \subseteq E$ be open, and $f \in \mathcal{C}^1(U,E)$. Assume $f(0) = 0$ and $f'(0) = \Id$. Let $r \ge 0$ and assume that $\cl{B_r(0)} \subseteq U$. Let $s \in (0,1)$, and assume that
  \begin{equation*}
    \norm{f'(z) - f'(x)}\le s
  \end{equation*}
  for every $z,x \in \cl{B_r(0)}$. If $y \in E$ and $\norm{y} \le (1- s)r$, then there exists a unique $x \in \cl{B_r(0)}$ such that $f(x) = y$.
\end{lemma}

\begin{proof}[The proof of Lemma \ref{lemma:covering_leaves}]
  Let $A$ be the constant appearing in \cite[equation (3.1)]{gouezel2006banach}, and let $\delta$ be the constant defined immediately afterwards.
  Let $W \in \tilde{\Sigma}$. Denote by $\tilde{W}$ the associated full admissible leaf, and by $\chi: \cl{B(x, A\delta)} \to (-2r_i/3,2r_i/3)^{d_u}$ the map defining $\tilde{W}$ i.e. $\tilde{W} = \psi_i \circ (\Id, \chi)(\cl{B(x, A\delta)})$.
  Fix $z \in B(x, \delta)$ and note that $B(z ,(A-1)\delta) \subseteq B(x,A \delta)$.
  For any $n \in \Z^+$ and $y \in \mathbb{T}^d$ let $\ell(n,y)$ be an index for which $T^{-n}(\psi_i(z, \chi(z)) + y) \in \psi_{\ell(n,y)} \left((-r_{\ell(n,y)}/2, r_{\ell(n,y)}/2)^d\right)$ (recall \ref{en:B5}).
  Let $\pi^s : \R^{d} \to \R^{d_s}$ be the projection onto the first $d_s$ components, and $\pi^u : \R^{d} \to \R^{d_u}$ the projection onto the last $d_u$ components.
  Note that $(\psi_{\ell(N,y)}^{-1} \circ T^{-N} \circ \tau_y \circ \psi_i \circ (\Id, \chi))(B(x, A\delta))$ is the union of finitely many disjoint, path-connected subsets; let $Q_N \subseteq B(x, A\delta)$ denote the pre-image under $\psi_{\ell(N,y)}^{-1} \circ T^{-N} \circ \tau_y \circ \psi_i \circ (\Id, \chi)$ of the particular subset containing $(\psi_{\ell(N,y)}^{-1} \circ T^{-N} \circ \tau_y \circ \psi_i \circ (\Id, \chi))(z)$.
  Define $F_{N} := \psi_{\ell(N,y)}^{-1} \circ T^{-N} \circ \tau_y \circ \psi_i \circ \restr{(\Id, \chi)}{Q_N}$.
  We will show that for sufficiently large $N$ one can use $F_N$ to construct an admissible leaf $W_z$ so that $\psi_i(z) \in T^{n}(W_z)$ and $T^{n}(W_z) \subseteq \tilde{W}+y$.

  \paragraph{Step I: The invertibility of $\pi^s \circ F_N$ in a neighbourhood of $(\pi^s \circ F_N)(z)$.}
  Let $N_1$ be the constant from Lemma \ref{lemma:distorted_leaf_tangent_space} and recall $H=\inf_{x \in \mathbb{T}^d} \norm{D_x T}$ from Lemma \ref{lemma:translated_leaf_expansion}.
  As remarked in the proof of \cite[Lemma 3.3]{gouezel2006banach} both $\psi_i^{-1}$ and $\psi_{\ell(N,y)}$ are $(1 + \kappa)$-Lipschitz.
  Let $N_2 > N_1$ be such that $\nu_s^{-n + N_1} H^{-N_1} > 1$ whenever $n > N_2$. By Lemma \ref{lemma:translated_leaf_expansion}, if $n > N_2$ then $T^{-n}$ expands distances on $\tilde{W}+y$ by at least $\nu_s^{-n + N_1} H^{-N_1}$.
  It is clear that $\tau_{y}^{-1}$ is $C_{\tau}$-Lipschitz by the definition of $C_{\tau}$.
  It is clear that $d((\Id,\chi)(a),(\Id,\chi)(b)) \ge d(a,b)$ for every $a,b \in \cl{B(z,A \delta)}$.
  Using the above estimates to bound the Lipschitz constant of $F_N^{-1}$ for $N > N_2$ we obtain
  \begin{equation*}
    d(F_N(a),F_N(b)) \ge \nu_s^{-N + N_1}\left(H^{N_1} C_\tau (1+ \kappa)^2 \right)^{-1} \quad \forall a,b \in Q_N.
  \end{equation*}
  As in \cite[Lemma 3.3]{gouezel2006banach} we have $\abs{\pi^s(v)} \ge (1 + c_{\ell(N,y)}^2)^{-1/2} \abs{v}$ whenever $v \in \mathcal{C}_{\ell(N,y)}$.
  Since $\sup_{i} c_{i} < 2 \kappa$ and, for $N > N_1$, the tangent space of $F_N$ is contained in $\mathcal{C}_{\ell(N,y)}$, for every $N > N_2$ we have
  \begin{equation}
    d((\pi_s \circ F_N)(a),(\pi_s \circ F_N)(b)) \ge \frac{\nu_s^{-N + N_1}}{H^{N_1}C_{\tau}( 1+ \kappa)^2 \sqrt{1 + 4 \kappa^2}},
  \end{equation}
  provided that $a,b$ are sufficiently close.
  Since $\nu_s < 1$ there exists $N_3 \ge N_2$ so that for each $N \ge N_3$ the map $\pi_s \circ F_N$ locally expands distances by at least
  \begin{equation*}
    \frac{\nu_s^{-N + N_1}}{H^{N_1}C_{\tau}( 1+ \kappa)^2  \sqrt{1 + 4 \kappa^2}} > \frac{A}{A-1},
  \end{equation*}
  from which it follows that $D_w(\pi^s \circ F_N)^{-1}$ exists for every $(\pi^s \circ F_N)(w) \in Q_N$ and satisfies
  \begin{equation}\label{eq:covering_leaves_1}
    \norm{D_w(\pi^s \circ F_N)^{-1}} \le \nu_s^{N - N_1}(H^{N_1}C_{\tau}( 1+ \kappa)^2  \sqrt{1 + 4 \kappa^2}) < \frac{A-1}{A}.
  \end{equation}
  We will now obtain a lower bound on the size of $Q_N$. Note that $(\Id, \chi)$ is $(\sqrt{1 + \kappa^2})$-Lipschitz due to the tangent space of $\chi$ being a subset of $\mathcal{C}_i$.
  Let $P = \sup_{x \in \mathbb{T}^d} \norm{D_x T^{-1}}$.
  From these estimates, as well as those in the previous paragraph, we may conclude that $F_N$ is $(P^N C_\tau (1+\kappa)^2 \sqrt{1 + \kappa^2})$-Lipschitz.
  Let
  \begin{equation*}
    L_N := \min\left\{(A-1)\delta, (P^N C_\tau (1+\kappa)^2 \sqrt{1 + \kappa^2})^{-1} \min_{j} r_j/2 \right\}.
  \end{equation*}
  We will prove that $B(z,L_N) \subseteq Q_N$.
  If $w \in (-r_i,r_i)^{d_s} \cap \cl{B(z, L_N)}$ then $w \in \cl{B(x, A\delta)}$ i.e. $w$ is in the domain of $\chi$. Moreover,
  \begin{equation}\begin{split}\label{eq:covering_leaves_2}
    d((T^{-N} \circ \tau_y \circ \psi_i \circ (\Id, \chi))(w)&, (T^{-N} \circ \tau_y \circ \psi_i \circ (\Id, \chi))(z)) \\
    &\le P^N C_\tau (1+\kappa) \sqrt{1+\kappa^2}\norm{w-z}.
  \end{split}\end{equation}
  Recall that $F_N(z) \in (-r_{\ell(N,y)}/2, r_{\ell(N,y)}/2)^{d}$ and so $\psi_{\ell(N,y)}$ is defined on $B\left(F_N(z), \min_j \frac{r_j}{2}\right)$. Since $\psi_{\ell(N,y)}$ is $(1+\kappa)$-Lipschitz and a bijection (onto its range), we have
  \begin{equation}\label{eq:covering_leaves_3}
    \psi_{\ell(N,y)}\left(B\left(F_N(z), \min_j \frac{r_j}{2}\right)\right) \supseteq B\left((T^{-N} \circ \tau_y \circ \psi_i \circ (\Id, \chi))(z), (1+ \kappa)^{-1}\min_j \frac{r_j}{2} \right).
  \end{equation}
  From \eqref{eq:covering_leaves_2} and \eqref{eq:covering_leaves_3} we deduce that if $w \in B(z, L_N)$ then $F_N(w)$ is defined and in $(-r_{\ell(N,y)}, r_{\ell(N,y)})^{d_s}$.
  Moreover, $F_n(B(z, L_N))$ is path-connected, being the image of a path-connected set under a continuous function, and so $B(z, L_N) \subseteq Q_N$.
  Let $S_N: B(0, L_N) \to \R^{d_s}$ be defined by
  \begin{equation*}
    S_N(w) = (D_z(\pi^s \circ F_N))^{-1} \cdot \left((\pi^s \circ F_N)(w + z) - (\pi^s \circ F_N)(z)\right).
  \end{equation*}
  Our goal is to apply Lemma \ref{lemma:inverse_func} to $S_N$, and then deduce the existence of $(\pi^s \circ F_N)^{-1}$ on some neighbourhood of $(\pi^s \circ F_N)(z)$ that is not too small, but this will take some work.
  For any $a,b \in B(0, L_N)$ we have
  \begin{equation}\label{eq:covering_leaves_4}\begin{split}
    \norm{D_aS_N - D_bS_N} &\le \norm{(D_z(\pi^s \circ F_N))^{-1}} \norm{D_{a+z}(\pi^s \circ F_N) - D_{b+z}(\pi^s \circ F_N)} \\
    &\le \norm{(D_z(\pi^s \circ F_N))^{-1}} \norm{D_{a+z} F_N - D_{b+z} F_N} \\
    &\le \norm{a-b} \sup_{w \in Q_N} \norm{D^2_w F_N} : = \norm{a-b} J_N,
  \end{split}\end{equation}
  where we have used \eqref{eq:covering_leaves_1} and the fact that $w \mapsto D_w \pi^s$ is constant and a contraction.
  Note that
  \begin{equation}\label{eq:covering_leaves_5}
    \frac{1- \sqrt{1 - 8A \delta J_{N_3} \norm{D_z(\pi^s \circ F_{N_3})^{-1}} }}{4 J_{N_3}} = \frac{2A \delta \norm{D_z(\pi^s \circ F_{N_3})^{-1}} }
    {1+ \sqrt{1 - 8A \delta J_{N_3} \norm{D_z(\pi^s \circ F_{N_3})^{-1}}}}.
  \end{equation}
  In the definition of $\tilde{\Sigma}$ we may assume that $\delta$ is as small as we like. Thus, in view of \eqref{eq:covering_leaves_1} and \eqref{eq:covering_leaves_5}, by choosing $\delta$ sufficiently small we may guarantee that
  \begin{equation}\label{eq:covering_leaves_6}
    0 < \frac{1- \sqrt{1 - 8A \delta J_{N_3} \norm{D_z(\pi^s \circ F_{N_3})^{-1}} }}{4 J_{N_3}} < L_{N_3}.
  \end{equation}
  If \eqref{eq:covering_leaves_6} holds then there exists $s \in (0,L_{N_3})$ so that
  \begin{equation}\label{eq:covering_leaves_7}
    ( 1 - 2J_{N_3} s)s \ge A \delta \norm{D_z(\pi^s \circ F_{N_3})^{-1}}.
  \end{equation}
  To summarise, we have proven the following
  \begin{enumerate}
    \item $S_{N_3}$ is well-defined on $\cl{B(0, s)}$ as $s < L_{N_3}$;
    \item $S_{N_3}(0) = 0$ and $S_N'(0) = \Id$; and
    \item By \eqref{eq:covering_leaves_4} we have $\norm{D_aS_{N_3} - D_bS_{N_3}} \le 2 s J_{N_3}$ for every $a,b \in \cl{B(0, s)}$.
  \end{enumerate}
  Thus we may apply Lemma \ref{lemma:inverse_func} to $S_{N_3}$ on $\cl{B(0, s)}$ to conclude the existence of an inverse $S_{N_3}^{-1}$ that is defined on $\cl{B(0, ( 1 - 2J_{N_3} s)s)}$.
  Using the definition of $S_{N_3}$ and \eqref{eq:covering_leaves_7} we recover the existence of an inverse $(\pi^s \circ F_{N_3})^{-1}$ on $\cl{B((\pi^s \circ F_{N_3})(z), A \delta)}$.

  \paragraph{Step II: The definition and properties of leaves covering $T^{-N}(\tilde{W}+y)$.}
  We may define a map $\chi_0: \cl{B((\pi^s \circ F_{N_3})(z), A \delta)} \to \R^{d_u}$ by
  \begin{equation*}
    \chi_0 = \pi^u \circ F_{N_3} \circ (\pi^s \circ F_{N_3})^{-1}.
  \end{equation*}
  Note that the graph of $\chi_0$ is a subset of $F_{N_3}(B(z,s))$ by construction.
  Since the tangent space of $F_{N_3}$ is contained in $\mathcal{C}_{\ell(N_3,y)}$ it follows that $\norm{D\chi_0} \le c_{\ell(N_3,y)}$.
  Hence, for $w \in \cl{B((\pi^s \circ F_{N_3})(z), A \delta)}$, we have
  \begin{equation}\label{eq:covering_leaves_8}
    \norm{\chi_0((\pi^s \circ F_{N_3})(z)) - \chi_0(w)} \le c_{\ell(N_3,y)} A\delta.
  \end{equation}
  Recall from the line before \eqref{eq:cone_distortion} that $\kappa < 1/2$, and from \ref{en:B4} that $c_{\ell(N_3,y)} < 2\kappa$.
  Thus, as $A \delta < \min_j r_j/6$ (see the sentence following \cite[equation (3.1)]{gouezel2006banach}), from \eqref{eq:covering_leaves_8} we have
  \begin{equation*}
    \norm{\chi_0((\pi^s \circ F_{N_3})(z)) - \chi_0(w)} < \min_j r_j/6.
  \end{equation*}
  Since $(\pi^s \circ F_{N_3})(z) \in (-r_{\ell(N_3,y))}/2, r_{\ell(N_3,y))}/2)^{d_s}$, it follows that $\chi_0(w) \in (-2r_{\ell(N_3,y)}/3, 2r_{\ell(N_3,y)}/3)^{d_u}$. Thus the image of $\chi_0$ is a subset of $(-2r_{\ell(N_3,y)}/3, 2r_{\ell(N_3,y)}/3)^{d_u}$.
  Since the $\mathcal{C}^{r+1}$ norm of $F_{N_3}$ may be bounded independently of $y \in \mathbb{T}^d$, $z \in B(x,\delta)$ and $W \in \tilde{\Sigma}$, by the inverse function theorem there exists some absolute $Y$ so that for any $\chi_0$ produced by the construction just carried out we have $\norm{\chi_0}_{\mathcal{C}^{r+1}} \le Y$.
  Thus the graph of $\chi_0$ belongs to $G_{\ell(N_3, y)}(Y)$ (recall the definition of the sets $G_i(K)$ from \cite[Section 3]{gouezel2006banach}).

  The issue at this stage is that $\norm{\chi_0}_{\mathcal{C}^{r+1}}$ may not be bounded by the constant $K$ set in \cite[Lemma 3.1]{gouezel2006banach} and so may not define a leaf in $\tilde{\Sigma}$. Instead, we show that $\chi_0$ may be covered by the image of higher-regularity leaves under some iterate of $T$. We note that the following construction is very similar to the one in the proof of \cite[Lemma 3.3]{gouezel2006banach}.

  For each $j \in \Z^+$ we define $\chi_j$ inductively as follows, starting with $j = 1$.
  As the graph of $\chi_{j-1}$ is in $G_{\ell(N_3 + j-1, y)}(R^{j-1}(Y))$, by Lemma \ref{lemma:leaf_regularisation} the image of $\psi_{\ell(N_3 + j, y)}^{-1} \circ T^{-1} \circ \psi_{\ell(N_3 + j-1, y)} \circ (\Id, \chi_{j-1})$ is in $G_{\ell(N_3 + j, y)}(R^{j}(Y))$ and is therefore the graph of a map $\chi_j'$ which contains $\psi_{\ell(N_3 + j, y)}^{-1} \circ T^{-(N_3 + j)} \circ \tau_y \circ \psi_i \circ (z, \chi(z))$.
  Using the expansivity of $T^{-1}$ as in \cite[Lemma 3.3]{gouezel2006banach}, one deduces that the domain of $\chi_j'$ contains the set
  \begin{equation}\label{eq:covering_leaves_9}
    \cl{B(\psi_{\ell(N_3 + j, y)}^{-1} \circ T^{-(N_3 + j)} \circ \tau_y \circ \psi_i \circ (z, \chi(z)), A \delta)}.
  \end{equation}
  We define $\chi_j$ to be the restriction of $\chi_j'$ to \eqref{eq:covering_leaves_9}.
  By Lemma \ref{lemma:leaf_regularisation} we have $\chi_{N_4} \in G_{\ell(N_3 + N_4, y)}(K)$, where $N_4$ denotes the constant $N(Y)$ given by Lemma \ref{lemma:leaf_regularisation} and $K$ is the constant from \cite[Lemma 3.1]{gouezel2006banach}.
  Thus the image of
  \begin{equation*}
    \cl{B(\psi_{\ell(N_3 + N_4, y)}^{-1} \circ T^{-(N_3 + N_4)} \circ \tau_y \circ \psi_i \circ (z, \chi(z)), \delta)}
  \end{equation*}
  under $\psi_{\ell(N_3 + N_4, y)} \circ (\Id, \chi_{N_4})$ is a leaf in $\tilde{\Sigma}$.

  \paragraph{Step III: Concluding.}
  We may apply this construction to any $z \in B(x, \delta)$ to produce a leaf $W_z \in \tilde{\Sigma}$ such that $T^{N_3 + N_4}(W_z) \subset \tilde{W} + y$. Moreover, the constants $N_3$ and $N_4$ are independent of $y \in \mathbb{T}^d$, $z \in B(x,\delta)$  and $W \in \tilde{\Sigma}$. Set $N = N_3 + N_4$.
  By varying $z$ we observe that the set of such leaves covers $T^{-N}(W + y)$.
  As in the end of \cite[Lemma 3.3]{gouezel2006banach}, the claim that the number of leaves required to cover $T^{-N}(W + y)$ may be bounded independently of $W$ and $y$ follows from \cite[Theorem 1.4.10]{hormander1990analysis}, as too does the existence of
  partitions of unity satisfying all of the required properties.
  Hence we have verified all the conclusions of Lemma \ref{lemma:covering_leaves}.
\end{proof}

By using Lemma \ref{lemma:covering_leaves} and adapting arguments from \cite[Section 6]{gouezel2006banach} we may now complete the proof of Proposition \ref{prop:translate_bound_cond} by proving \eqref{eq:translate_bound_cond_1}.
Let $N$ be the constant from Lemma \ref{lemma:covering_leaves}. We will show that $\norm{\tau_y \LL^N}_{p,q}$ can be bounded independently of $y$. The argument for bounding $\norm{\tau_y \LL^N}_{p-1,q+1}$ is identical.
Moreover, we will only derive the inequality \eqref{eq:translate_bound_cond_1} for the spaces $\tilde{B}^{p,q}$ and $\tilde{B}^{p-1,q+1}$, since it is straightforward to then derive the corresponding inequality for their complexifications.
We begin by bounding $\norm{\tau_y \LL^N }^{-}_{0,q}$.
Let $h \in \mathcal{C}^r(\mathbb{T}^d, \R)$, $W \in \tilde{\Sigma}$, $\varphi \in \mathcal{C}^{q}_0(W, \R)$ satisfy $\norm{\varphi}_{\mathcal{C}^{q}} \le 1$, and $y \in \mathbb{T}^d$.
Let $J_W \tau_{-y}$ denote the Jacobian of $\tau_{-y} : W + y \to W$.
Then
\begin{equation}\label{eq:translate_bound_cond_4}
  \intf_W \tau_y \LL^N h \cdot \varphi d\Omega = \intf_{W+y} \LL^N h \cdot \varphi \circ \tau_{-y} \cdot J_W \tau_{-y} d\Omega.
\end{equation}
Recall that $\LL^N h = (h \abs{\det DT^N}^{-1}) \circ T^{-N}$ and let $J_{W+y}T^N$ denote the Jacobian of $T^{N}: T^{-N}(W+y) \to W+y$.
Let $\{W_s\}_{s=1}^j$, $\{\rho_s\}_{s=1}^j$ satisfy the conclusion of Lemma \ref{lemma:covering_leaves}.
By changing coordinates and applying Lemma \ref{lemma:covering_leaves} to \eqref{eq:translate_bound_cond_4} we obtain
\begin{equation}\begin{split}\label{eq:translate_bound_cond_5}
  \int_{W+y} \LL^N h &\cdot \varphi \circ \tau_{-y} \cdot J_W \tau_{-y} \, \mathrm{d}\Omega \\
  &= \intf_{T^{-N}(W+y)} h \cdot \abs{\det DT^N}^{-1} \cdot (\varphi \circ \tau_{-y} \circ T^N) \cdot (J_W \tau_{-y}) \circ T^N \cdot J_{W+y}T^N d\Omega \\
  &= \sum_{s=1}^j \intf_{W_s} h \cdot \abs{\det DT^N}^{-1} \cdot (\varphi \circ \tau_{-y} \circ T^N) \cdot (J_W \tau_{-y}) \circ T^N\cdot \rho_s \cdot J_{W+y}T^N d\Omega.
\end{split}\end{equation}
By the definition of $\norm{\cdot}_{0,q}$ the final expression in \eqref{eq:translate_bound_cond_5} is bounded above by
\begin{equation*}
  \norm{h}_{0,q} \sum_{s=1}^j \norm{\abs{\det DT^N}^{-1} \cdot (\varphi \circ \tau_{-y} \circ T^N) \cdot (J_W \tau_{-y}) \circ T^N \cdot \rho_s \cdot J_{W+y}T^N}_{\mathcal{C}^q(W_s)}.
\end{equation*}
Recall that $q \le r-1$. Since $T$ is a $\mathcal{C}^r$ diffeomorphism and $\mathbb{T}^d$ is compact, it follows that the $\mathcal{C}^q(W_s)$ norms of $J_{W+y}T^N$ and $\abs{\det DT^N}^{-1}$ are bounded independently of $W$ and $y$.
Using continuity and compactness, we observe that $\sup_{x,y \in \mathbb{T}^d}\norm{D^k_x\tau_{y}}$ is finite for every positive integer $k$.
Thus the $\mathcal{C}^q(W_s)$ norm of $\varphi \circ \tau_{-y} \circ T^N$ is bounded independently of $y$ and $\varphi$ (provided $\norm{\varphi}_{\mathcal{C}^q(W_s)} \le 1$).
Similarly, it is clear that the $\mathcal{C}^q(W_s)$ norm of $(J_W \tau_{-y}) \circ T^N$ is bounded independently of $W$ and $y$.
Recall from Lemma \ref{lemma:covering_leaves} that $S$ and $\mathcal{C}^q$ norms of each $\rho_s$ are both bounded independently of $W$ and $y$.
Hence, as $\norm{\cdot}_{\mathcal{C}^q(W_s)}$ is sub-multiplicative, there exists some $C_{0,q} > 0$ such that
\begin{equation*}
  \sum_{s=1}^j \norm{\abs{\det DT^N}^{-1} \cdot (\varphi \circ \tau_{-y} \circ T^N) \cdot \rho_s \cdot J_{W+y}T^N}_{\mathcal{C}^q(W_s)} \le C_{0,q}
\end{equation*}
for every choice of $W, y$ and $\varphi$ with  $\norm{\varphi}_{\mathcal{C}^q} \le 1$.
Thus, by taking supremums of the terms in \eqref{eq:translate_bound_cond_4} we have
\begin{equation}\label{eq:translate_bound_cond_6}
  \norm{\tau_{y} \LL^N h}_{0,q}^{-} \le C_{0,q} \norm{h}_{0,q}.
\end{equation}
We turn to bounding $\norm{\tau_{y} \LL^N h}_{k,q+k}^{-}$ for $0 < k \le p$. Let $h$, $W$ and $y$ be as before.
Suppose that $\varphi \in \mathcal{C}^{q+k}_0(W, \R)$ satisfies $\norm{\varphi}_{\mathcal{C}^{q+k}} \le 1$ and that $\{v_i\}_{i=1}^k \subseteq \mathcal{V}^r(W)$ is such that $\norm{v_i}_{\mathcal{C}^r} \le 1$.
Then
\begin{equation*}
  \intf_W (v_1 \dots v_k) (\tau_y \LL^N h) \cdot \varphi d\Omega = \intf_W \tau_y (\tilde{v}_{1,y} \dots \tilde{v}_{1,k} \LL^N h) \cdot \varphi d\Omega,
\end{equation*}
where $\tilde{v}_{i,y}(x) := (D_{x-y}\tau_y) v_i(x-y)$.
With $J_W \tau_{-y}$, $\{W_s\}_{s=1}^j$ and $\{\rho_s\}_{s=1}^j$ as before we have
\begin{equation}\begin{split}\label{eq:translate_bound_cond_7}
  \int_W \tau_y& (\tilde{v}_{1,y} \dots \tilde{v}_{1,k} \LL^N h) \cdot \varphi \, \mathrm{d} \Omega= \intf_{W+y}  (\tilde{v}_{1,y} \dots \tilde{v}_{1,k} \LL^N h) \cdot \varphi \circ \tau_{-y}  \cdot J_W \tau_{-y} d\Omega \\
   &= \sum_{s=1}^j \intf_{T^N(W_s)} (\tilde{v}_{1,y} \dots \tilde{v}_{1,k} \LL^N h) \cdot \varphi \circ \tau_{-y} \cdot \rho_s \circ T^{-N} \cdot J_W \tau_{-y} d\Omega.
\end{split}\end{equation}
Since the $\mathcal{C}^{q+k}$ norms of $\varphi \circ \tau_{-y}$ and $J_W \tau_{-y}$ are bounded independently of $\varphi$, $y$ and $W$, we may replace $\varphi \circ \tau_{-y} \cdot J_W \tau_{-y}$ by some $\phi \in \mathcal{C}^{q+k}_0(W, \R)$ with $\mathcal{C}^{q+k}$ norm bounded independently of $\varphi$, $y$ and $W$.
Additionally, the $\mathcal{C}^{r}$ norm of each $\tilde{v}_{i,y}$ may be bounded independently of $y$ and $W$ due to the $\sup_{x,y \in \mathbb{T}^d} \norm{D_x^\ell\tau_y}$ being finite for each positive integer $\ell$.
Upon replacing $\varphi \circ \tau_{-y} \cdot J_W \tau_{-y}$ with $\phi$, the expression on the right hand side of \eqref{eq:translate_bound_cond_7} is exactly in the form of \cite[(6.4)]{gouezel2006banach}.
Using the arguments from \cite[Lemma 6.3]{gouezel2006banach}, one then obtains a bound of the form
\begin{equation*}
  \abs{\intf_W (v_1 \dots v_k) (\tau_y \LL^N h) \cdot \varphi d\Omega} \le C_{p,q} \norm{h}_{p,q} + C_{p-1,q+1} \norm{h}_{p-1,q+1},
\end{equation*}
for some $C_{p,q}, C_{p-1,q+1}> 0$ that are independent of $h$, $W$, $k$, $y$, $\varphi$ and each $v_i$.
By the definition of $\norm{\cdot}_{k,q+k}^{-}$, and as $\norm{\cdot}_{p,q}$ dominates $\norm{\cdot}_{p-1,q+1}$, we therefore have
\begin{equation}\label{eq:translate_bound_cond_8}
  \norm{\tau_{y}\LL^N h}_{k,q+k}^{-} \le (C_{p,q} + C_{p-1,q+1}) \norm{h}_{p,q}.
\end{equation}
The required bound is obtained by considering \eqref{eq:translate_bound_cond_6}, \eqref{eq:translate_bound_cond_8} and the definition of $\norm{\cdot}_{p,q}$.
Thus we have established \eqref{eq:translate_bound_cond_1}, which completes the proof of Proposition \ref{prop:translate_bound_cond}.

\section{Properties of the perturbed cat map}\label{sec:perturbed_cat}
In this appendix we consider the maps $T_\delta : \mathbb{T}^2 \to \mathbb{T}^2$, $\delta \in \R$, defined by
\begin{equation*}
  T_\delta(x_1, x_2) = (2x_1 + x_2, x_1 + x_2) + \delta(\cos(2\pi x_1), \sin(4 \pi x_2 +1)),
\end{equation*}
and prove the claim made in Proposition \ref{thm:perturbed_cat_props} that for $0\le \delta < 0.0108$ the map $T_\delta$ is an Anosov diffeomorphim and satisfies the conditions of Proposition \ref{prop:translate_bound_cond}.
The proof is broken up into Lemmas \ref{lemma:perturbed_cat_diffeo}, \ref{lemma:perturbed_cat_anosov} and \ref{lemma:fejer_cond_for_perturbed_cat}.

Throughout this section we denote the Euclidean norm on $\R^2$ (and the associated operator norm) by $\wnorm{\cdot}$, and the usual Euclidean inner product by $\ip$.
We begin by proving that $T_\delta$ is a diffeomorphism for sufficiently small $\delta$ by using a quantitative version of the inverse function theorem (Lemma \ref{lemma:inverse_func}).

\begin{lemma}\label{lemma:perturbed_cat_diffeo}
  If $\delta \in [0, 0.0108)$ then $T_\delta$ is a diffeomorphism.
  \begin{proof}
    We have
    \begin{equation*}
      D_{(x_1,x_2)}T_\delta =
      \begin{pmatrix}
        2 - 2\pi \delta \sin(2\pi x_1) & 1 \\
        1 & 1 + 4 \pi \delta  \cos(4\pi x_2 +1)\\
      \end{pmatrix},
    \end{equation*}
    and
    \begin{equation}\label{eq:perturbed_cat_diffeo_1}
      \det D_{(x_1,x_2)}T_\delta =1 + 8\pi \delta \cos(4\pi x_2 +1) - 2\pi \delta \sin(2\pi x_1) - 8 \pi^2 \delta^2 \sin(2\pi x_1)\cos(4\pi x_2 +1).
    \end{equation}
    In particular, $D_{(0,0)}T_\delta$ is invertible if $\abs{\delta} < 1/(8\pi)$.
    Denote by $\bar{T_\delta}: \R^2 \to \R^2$ the lifting of $T_\delta$. For $z \in \R$ let $\widetilde{z}$ denote the equivalence class containing $z$ in $\mathbb{T} = \R / \Z$.
    Define $R_{\delta}: \R^2 \to \R^2$ by
    \begin{equation*}
      R_{\delta}(y_1,y_2) := (D_{(0,0)}T_\delta)^{-1} \left( \bar{T}_\delta(y_1,y_2) - \bar{T}_\delta(0,0)\right).
    \end{equation*}
    Note that $R_{\delta}(0,0) = 0$ and $D_{(0,0)}R_{\delta} = \Id$. We will estimate
    \begin{equation*}
      \wnorm{D_{(y_1,y_2)}R_{\delta} - D_{(w_1,w_2)}R_{\delta}} = \wnorm{(D_{(0,0)}T_\delta)^{-1} \left(D_{( \widetilde{y_1}, \widetilde{y_2})}T_\delta - D_{(\widetilde{w_1}, \widetilde{w_2})}T_\delta\right)}.
    \end{equation*}
    We clearly have
    \begin{equation}\label{eq:perturbed_cat_diffeo_2}
      \wnorm{D_{( \widetilde{y_1}, \widetilde{y_2})}T_\delta - D_{(\widetilde{w_1}, \widetilde{w_2})}T_\delta} \le 8\pi\abs{\delta}.
    \end{equation}
    Use the fact that the Frobenius norm dominates the Euclidean operator norm, we have
    \begin{equation}\label{eq:perturbed_cat_diffeo_3}
      \wnorm{D_{(0,0)}T_\delta^{-1}} = \frac{\wnorm{D_{(0,0)}T_\delta}}{\abs{\det D_{(0,0)}T_\delta}}
      \le \frac{\sqrt{6 + (1+ 4\pi\delta \cos(1))^2}}{\abs{\det D_{(0,0)}T_\delta}}.
    \end{equation}
    Thus, by \eqref{eq:perturbed_cat_diffeo_1}, \eqref{eq:perturbed_cat_diffeo_2} and \eqref{eq:perturbed_cat_diffeo_3},
    \begin{equation*}
      \wnorm{D_{(y_1,y_2)}R_{\delta} - D_{(w_1,w_2)}R_{\delta}} \le \frac{ 8\pi \abs{\delta}\sqrt{6 + (1+ 4\pi\delta \cos(1))^2}}{1-8 \pi \delta} := s.
    \end{equation*}
    If $\abs{\delta} < 0.0108$ then $s < 1$.
    Then $R_{\delta}$ verifies the conditions of Lemma \ref{lemma:inverse_func} and has an inverse $R_{\delta}: B(0, (1-s)r) \to \R^2$, where $r > 0$. Since there is no dependence on $r$ in the above procedure, we may extend $R_{\delta}^{-1}$ to $\R^2$.
    Thus $\bar{T}_\delta$ is invertible, and so $T_\delta$ is invertible too. It is standard that $T_\delta^{-1}$ and $T_\delta$ have the same smoothness.
  \end{proof}
\end{lemma}

Let $x \in \mathbb{T}^2$. The eigenvalues of $D_x T_0$ are $\lambda = \frac{3 - \sqrt{5}}{2} < 1$ and $\lambda^{-1} = \frac{3 + \sqrt{5}}{2} > 1$.
Let $\widetilde{E}^s(x)$ be the span of $\left(1, \frac{-\sqrt{5} - 1}{2}\right)$ and $\widetilde{E}^u(x)$ be the span of $\left(1, \frac{\sqrt{5} - 1}{2}\right)$. Note that $\widetilde{E}^s(x)$ and $\widetilde{E}^u(x)$ are the eigenspaces of $\lambda$ and $\lambda^{-1}$, respectively.
It is trivial that the spaces $\widetilde{E}^u(x)$ and $\widetilde{E}^s(x)$ depend continuously on $x$, and that $\widetilde{E}^u(x) \oplus \widetilde{E}^s(x) = T_x \mathbb{T}^2$ for every $x \in \mathbb{T}^2$.
Let $\Pi^{u}, \Pi^{s} : \R^2 \to \R^2$ denote the orthogonal projections onto $\widetilde{E}^u(x)$ and $\widetilde{E}^s(x)$, respectively. Since $T_0$ is symmetric, $\widetilde{E}^u(x) \perp \widetilde{E}^s(x)$ and so $\Id - \Pi^{s} = \Pi^u$.
For $\alpha > 0$ define
\begin{equation*}
  K^u_\alpha(x) = \{ v \in T_x \mathbb{T}^2 : \wnorm{\Pi^{s} v  } \le \alpha \wnorm{\Pi^{u} v}\}, \text{ and }
  K^s_\alpha(x) = \{ v \in T_x \mathbb{T}^2 : \wnorm{\Pi^{u} v  } \le \alpha \wnorm{\Pi^{s} v}\}.
\end{equation*}
To prove that $T_\delta$ is an Anosov diffeomorphism it remains to prove that $\mathbb{T}^2$ is a hyperbolic set for $T_\delta$.
We do this by verifying the conditions of the following result, which we have adapted to our setting for simplicity.

\begin{proposition}[{\cite[Proposition 5.4.3]{brin2002introduction}}]\label{prop:invariant_cones_anosov}
  If there exists $\alpha > 0$ such that for every $x \in \mathbb{T}^2$ we have
  \begin{enumerate}[label=(A\arabic*)]
    \item \label{en:A1} $D_x T_\delta (K^u_\alpha(x)) \subseteq K^u_\alpha(T_\delta(x)) $ and $D_x T_\delta^{-1} (K^s_\alpha(x)) \subseteq K^s_\alpha(T_\delta^{-1}(x)) $; and
    \item \label{en:A2} $\wnorm{D_x T_\delta v } < \wnorm{v}$ for $v \in K^s_\alpha(x) \setminus \{ 0\}$ and $\wnorm{D_x T_\delta^{-1} v } < \wnorm{v}$ for $v \in K^u_\alpha(x) \setminus \{ 0\}$.
  \end{enumerate}
  Then there are constants $\nu_{u,\delta} > 1$ and $ 0 < \nu_{s,\delta} < 1$, and for each $x \in \mathbb{T}^2$, subspaces $E^s_\delta(x)$ and $E^u_\delta(x)$ such that
  \begin{enumerate}
    \item $T_x \mathbb{T}^2 = E^s_\delta(x) \oplus E^u_\delta(x)$;
    \item $D_x T_\delta^{-1} (E^s_\delta(x)) = E^s(T_\delta^{-1}(x))$ and $D_x T_\delta (E^u_\delta(x)) = E^u(T_\delta(x))$;
    \item $\wnorm{\restr{D_x T_\delta}{E^s_\delta(x)}} \le \nu_{s,\delta} $ and $\wnorm{\restr{D_x T_\delta^{-1}}{E^u_\delta(x)}} \le \nu_{u,\delta}^{-1}$; and
    \item $E^s_\delta(x) \subseteq K^s_\alpha(x)$ and $E^u_\delta(x) \subseteq K^u_\alpha(x)$.
  \end{enumerate}
  In particular, $T_\delta$ is Anosov.
\end{proposition}

\begin{lemma}\label{lemma:perturbed_cat_anosov}
If $\delta \in [0, 0.0108)$ then $T_\delta$ is Anosov and the conclusion of Proposition \ref{prop:invariant_cones_anosov} holds with $\alpha = 0.11872$.
\begin{proof}
We first diagonalise $D_xT_0$ as $R_\theta^{-1} \Lambda R_\theta$, where $\Lambda$ is a $2\times 2$ diagonal matrix with the vector $[\lambda, 1/\lambda]$ on the diagonal, and $R_\theta$ is clockwise rotation by angle
$\theta=\tan^{-1}((1-\sqrt{5})/2)$.
Note that $D_x T_\delta = D_x T_0 + \Delta_x$ where $\Delta_x : \R^2 \to \R^2$ is defined by the matrix
  \begin{equation*}
   \Delta_x= \begin{pmatrix}
      -2\pi \delta \sin(2\pi x_1) & 0\\
      0 & 4 \pi \delta  \cos(4\pi x_2 +1)\\
    \end{pmatrix}.
  \end{equation*}
We use the shorthand $\delta_1=-2\pi \delta \sin(2\pi x_1)$ and $\delta_2=4 \pi \delta  \cos(4\pi x_2 +1)$.
In order to satisfy the second part of (A1) of Proposition \ref{prop:invariant_cones_anosov}, we require that $\Lambda+R_\theta\Delta_xR_\theta^{-1}$  preserve $\mathcal{K}^u_\alpha:=\{(\beta, \gamma)^\top\in\mathbb{R}^2: |\beta|\le \abs{\gamma}\alpha\}$.
One may confirm that
\begin{equation}
\label{rot}
    R_\theta\Delta_xR_\theta^{-1} =
    \begin{pmatrix}
      \delta_1\cos^2\theta+\delta_2\sin^2\theta & (1/2)\sin(2\theta)(\delta_2 - \delta_1)\\
    (1/2)\sin(2\theta)(\delta_1-\delta_2) & \delta_2\cos^2\theta+\delta_1\sin^2\theta \\
    \end{pmatrix}.
  \end{equation}
Multiplying $\Lambda+R_\theta\Delta_xR_\theta^{-1}$ with the vectors $(\alpha, 1)^\top$ and $(-\alpha, 1)^\top$ we see that a sufficient condition to preserve $\mathcal{K}^u_\alpha$ is that $$\frac{(\lambda+\delta')\alpha+\delta'}{1/\lambda-\delta'-\delta'\alpha}\le\alpha,$$
where $\delta'=\max\{\sup_x\delta_1,\sup_x\delta_2\}=4\pi\delta$.
Since $1/\lambda-\delta'-\delta'\alpha > 0$ for $0 \le \delta \le0.1862$ we may rearrange the above in terms of $\delta$ to obtain
\begin{equation}\label{eq:perturbed_cat_anosov_1}
  \delta\le\frac{\alpha(1/\lambda-\lambda)}{4\pi(\alpha+1)^2}.
\end{equation}

Since $\delta \in [0,0.0108)$, by Lemma \ref{lemma:perturbed_cat_diffeo} the map $T_\delta$ is a diffeomorphism.
To satisfy the first part of (A1) of Proposition \ref{prop:invariant_cones_anosov}, using the notation above, we
note that
\begin{equation}
\label{inv}
    (D_xT_\delta)^{-1} =
   (1/\det(D_xT_\delta)) \begin{pmatrix}
      1+\delta_2 & -1 \\
      -1 & 2+\delta_1\\
    \end{pmatrix},
  \end{equation}
and that for the purposes of cone preservation we need not consider the determinant factor.
Therefore, $(D_xT_\delta)^{-1}=(1/\det(D_xT_\delta))((D_xT_0)^{-1}+\Delta'_x)$, where
\begin{equation}
\label{Delprime}
    \Delta'_x =
    \begin{pmatrix}
      \delta_2 & 0 \\
      0 & \delta_1\\
    \end{pmatrix}.
  \end{equation}
The cone preservation condition will be implied by the preservation of $\mathcal{K}^s_\alpha:=\{(\gamma, \beta)^\top\in\mathbb{R}^2: |\beta|\le \abs{\gamma} \alpha\}$ by $D_xT_\delta^{-1}$.
Multiplying $\Lambda^{-1}+R_\theta\Delta'_xR_\theta^{-1}$ with the vectors $(1, \alpha)^\top$ and $(1, -\alpha)^\top$ yields an identical set of calculations to those for $\mathcal{K}^u_\alpha$, resulting in the same bound for $\delta$ as in \eqref{eq:perturbed_cat_anosov_1}.
Substituting $\alpha=0.11872$ into this bound yields a numerical upper bound for $\delta$ of 0.0169, which is larger than the value reported in the proposition statement.

To verify (A2) we demonstrate contraction for elements of $\mathcal{K}_\alpha^s$;  the same contractions occur in the original (unrotated) cones $K^s_\alpha(x)$ and $K^u_\alpha(x)$ under $D_xT_\delta$ and $D_xT^{-1}_\delta$, respectively.
Writing $\Lambda+R_\theta\Delta_xR_\theta^{-1}=\left(
                                               \begin{array}{cc}
                                                 \lambda+a & b \\
                                                 c &\lambda^{-1}+d \\
                                               \end{array}
                                             \right)$
and multiplying by the unit vector $(1/\sqrt{1+\beta^2})(1, \beta)^\top$ for $-\alpha\le\beta\le \alpha$, the square of the norm of this vector is $((\lambda+a)^2+b^2\beta^2+2(\lambda+a)b\beta+c^2+(\lambda^{-1}+d)^2\beta^2+2c\beta(\lambda^{-1}+d))/(1+\beta^2)$.
We require the above expression to be strictly less than 1 for contraction.
Grouping terms to obtain a quadratic in $\beta$ we wish to show \begin{equation}\label{betapoly}
\beta^2(b^2+(\lambda^{-1}+d)^2-1)+2\beta((\lambda+a)+c(\lambda^{-1}+d))+((\lambda+a)^2+c^2-1)<0
\end{equation}
for $-\alpha\le\beta\le\alpha$.
This quadratic has a local minimum since $\lambda^{-1} \ge 1 + \abs{d}$ whenever $\delta \in [0,0.1288)$; therefore the maxima are at $\beta= \pm \alpha$.
Using the fact that $\max\{ \abs{a}, \abs{b} , \abs{c} , \abs{d}\} \le \delta'$ one may readily check that contraction occurs at $\beta=\pm\alpha$ for $\delta' \in [0, 0.6734)$ or, equivalently, for $\delta \in [0, 0.0536)$.

The contraction of vectors in $\mathcal{K}_\alpha^u$ under $D_xT_\delta^{-1}$ follows similarly.
With the notation above, one easily verifies $1-(5/2)\delta'-(\delta')^2/2\le\det(D_xT_\delta)$.
Replacing the two `1's in \eqref{betapoly} with the factor $1-(5/2)\delta'-(\delta')^2/2$, one verifies as above that the polynomial \eqref{betapoly} has positive leading term for $\delta \in [0,0.1288)$ and is negatively valued for $-\alpha\le\beta\le \alpha$ provided that $\delta \in [0, 0.0293)$.
Thus vectors in $\mathcal{K}_\alpha^u$ are contracted under $D_xT_\delta^{-1}$ for $\delta$ in the advertised range.

As we have verified all the conditions of Proposition \ref{prop:invariant_cones_anosov} for $T_\delta$ whenever $\delta \in [0,0.0108)$, it follows that $T_\delta$ is Anosov for $\delta$ in the same range.
\end{proof}
\end{lemma}

To complete the proof of Proposition \ref{thm:perturbed_cat_props} it now suffices to prove that $T_\delta$ satisfies the conditions of Proposition \ref{prop:translate_bound_cond}.
Denote by $\pi^s_x$ (resp. $\pi^u_x$) the projection onto $E^s_\delta(x)$ along $E^u_\delta(x)$ (resp. $E^u_\delta(x)$ along $E^s_\delta(x)$).
Let $w_s$ and $w_u$ be the unit vectors in the rays defined by $ \left(1, \frac{-\sqrt{5} - 1}{2}\right)$ and $\left(1, \frac{\sqrt{5} - 1}{2}\right)$, respectively.
Let $w_u(x)$ and $w_s(x)$ be the unit vectors in $E^s_\delta(x)$ and $E^u_\delta(x)$ for which $\langle w_u(x) , w_u \rangle> 0 $ and  $\langle w_s(x) , w_s \rangle> 0$. For $v \in \R^2$ we denote by $v^\perp$ the vector obtained by rotating $v$ anticlockwise by $\pi/2$ about the origin. In particular, $w_s^\perp = w_u$ and $w_u^\perp = -w_s$.
For $v \in T_x \mathbb{T}^2$ let
\begin{equation*}
  \norm{v}_0 = \sqrt{\wnorm{\pi^s_x v}^2 + \wnorm{\pi^u_x v}^2}.
\end{equation*}
We can recover a Riemannian metric from $\norm{\cdot}_0$ by the polarisation identity.
By Proposition \ref{prop:invariant_cones_anosov} we have $\wnorm{\restr{D_x T_\delta}{E^s_\delta(x)}} \le \nu_{s,\delta}$ and $\wnorm{\restr{D_x T_\delta^{-1}}{E^u_\delta(x)}} \le \nu_{u,\delta}^{-1}$.
Thus the metric induced by $\norm{\cdot}_0$ satisfies \ref{en:M2} and is adapted.
In the following two lemmas we collect some useful inequalities, before proving that $T_\delta$ satisfies the conditions of Proposition \ref{prop:translate_bound_cond} for all $\delta \in [0,0.0108)$ in Lemma \ref{lemma:fejer_cond_for_perturbed_cat}. The first such bound follows from basic trigonometry.

\begin{lemma}\label{lemma:new_cone_approx}
  If $v_i \in K^s_\alpha(x)$, $i=1,2$ with $\wnorm{v_i} =1$ and $\langle v_1, v_2\rangle > 0$, then $\wnorm{v_1 - v_2} \le \frac{2\alpha}{\sqrt{1 + \alpha^2}}$.
  Similar statements hold for $v_i \in K^u_\alpha(x)$.
\end{lemma}

\begin{lemma}\label{lemma:perturbed_cat_zero_norm}
  If $\alpha < 1$, then for every $v \in T_x \mathbb{T}^2$ we have
  \begin{equation*}
    \frac{\sqrt{(1-\alpha^2)^{2} -\beta(\alpha)}}{1-\alpha^2} \abs{v} \le \norm{v}_0 \le \frac{\sqrt{(1-\alpha^2)^{2} +\beta(\alpha)}}{1-\alpha^2} \abs{v},
  \end{equation*}
  where
  \begin{equation*}
    \beta(\alpha) := \sqrt{2}\alpha(3\alpha + \sqrt{2 + 3\alpha^2})\sqrt{1 + \alpha\sqrt{2+3\alpha^2}}.
  \end{equation*}
  \begin{proof}
    Let $v \in T_x \mathbb{T}^2$ with $\wnorm{v} = 1$. By writing $v = \pi^s_{x} v + \pi^u_{x} v$ we find that
    \begin{equation}\label{eq:perturbed_cat_zero_norm_1}
      \norm{v}_0 = \sqrt{\wnorm{\pi^s_{x} v}^2 + \wnorm{\pi^u_{x} v}^2} = \sqrt{ 1 - 2 \langle \pi^s_{x} v, \pi^u_{x} v \rangle }.
    \end{equation}
    One verifies that $\pi^s_{x} v = \frac{\langle w_u(x)^\perp, v \rangle}{\langle w_u(x)^\perp, w_s(x) \rangle} w_s(x)$ and $\pi^u_{x} v = \frac{\langle w_s(x)^\perp , v\rangle}{\langle w_s(x)^\perp, w_u(x)\rangle } w_u(x)$.
    Hence
    \begin{equation}\label{eq:perturbed_cat_zero_norm_2}
      \abs{\langle \pi^s_{x} v, \pi^u_{x} v \rangle} = \abs{\frac{\langle w_s(x), w_u(x)\rangle \langle w_s(x)^\perp, v\rangle \langle w_u(x)^\perp, v\rangle }{\langle w_s(x)^\perp, w_u(x)\rangle \langle w_u(x)^\perp, w_s(x)\rangle}}.
    \end{equation}
    Since $w_s(x) \in K^s_\alpha(x)$ and $w_u(x) \in K^u_\alpha(x)$, basic trigonometry yields
    \begin{equation}\label{eq:perturbed_cat_zero_norm_3}
      \abs{\langle w_s(x), w_u(x)\rangle} \le \cos\left(\frac{\pi}{2} - 2\tan^{-1}(\alpha)\right) = 2\sin(\tan^{-1}(\alpha))\cos(\tan^{-1}(\alpha))= \frac{2\alpha}{1+ \alpha^2}.
    \end{equation}
    Alternatively, as $w_s(x)^\perp \in K^u(x)$,
    \begin{equation}\label{eq:perturbed_cat_zero_norm_4}
      \abs{\langle w_s(x)^\perp, w_u(x)\rangle} \ge \cos(2\tan^{-1}(\alpha)) = \cos^2(\tan^{-1}(\alpha)) - \sin^2(\tan^{-1}(\alpha)) = \frac{1 - \alpha^2}{1 + \alpha^2}.
    \end{equation}
    The same argument yields the same lower bound for $\abs{ \langle w_u(x)^\perp, w_s(x)\rangle}$.
    Writing $v = Aw_s(x) + Bw_s(x)^\perp$, we have
    \begin{equation*}
      \langle w_s(x)^\perp, v\rangle \langle w_u(x)^\perp, v\rangle  = B (-A + \langle w_u(x)^\perp + w_s(x), v\rangle).
    \end{equation*}
    Since $w_s(x)$ and $w_s(x)^\perp$ are orthogonal and $\abs{v} =   1$, we have $\abs{B} = \sqrt{1 - A^2}$. On the other hand, by Lemma \ref{lemma:new_cone_approx} and Cauchy-Schwarz we have $\abs{\langle w_u(x)^\perp + w_s(x), v\rangle} \le \frac{2\alpha}{\sqrt{1+\alpha^2}}$. Upon substituting we have
    \begin{equation}\label{eq:perturbed_cat_zero_norm_5}
      \abs{\langle w_s(x)^\perp, v\rangle \langle w_u(x)^\perp, v\rangle}  \le \sqrt{1 - A^2} \left(\abs{A} + \frac{2\alpha}{\sqrt{1+\alpha^2}}\right).
    \end{equation}
    We may bound the right-hand side of \eqref{eq:perturbed_cat_zero_norm_5} from above by differentiating with respect to $A$ and solving the resulting quadratic equation (noting that we only have to consider the case where $A \ge 0$ due to the symmetry about $A = 0$ in the right hand side of \eqref{eq:perturbed_cat_zero_norm_5}). In particular, \eqref{eq:perturbed_cat_zero_norm_5} is maximised when
    \begin{equation*}
      \abs{A} = \frac{-\alpha + \sqrt{2 + 3\alpha^2}}{2\sqrt{1+\alpha^2}},
    \end{equation*}
    which, when substituted into \eqref{eq:perturbed_cat_zero_norm_5}, yields
    \begin{equation}\label{eq:perturbed_cat_zero_norm_6}
      \abs{\langle w_s(x)^\perp, v\rangle \langle w_u(x)^\perp, v\rangle} \le \frac{(3\alpha + \sqrt{2 + 3\alpha^2})\sqrt{1  + \alpha \sqrt{2 + 3 \alpha^2}}}{2\sqrt{2}(1 + \alpha^2)}.
    \end{equation}
    Applying \eqref{eq:perturbed_cat_zero_norm_3}, \eqref{eq:perturbed_cat_zero_norm_4} and \eqref{eq:perturbed_cat_zero_norm_6} to \eqref{eq:perturbed_cat_zero_norm_2} yields
    \begin{equation}\label{eq:perturbed_cat_zero_norm_7}
      \abs{\langle \pi^s_{x} v, \pi^u_{x} v \rangle} \le
      \frac{\alpha(3\alpha + \sqrt{2 + 3\alpha^2})\sqrt{1  + \alpha \sqrt{2 + 3 \alpha^2}}}{\sqrt{2}(1 - \alpha^2)^2}.
    \end{equation}
    Hence, upon substituting \eqref{eq:perturbed_cat_zero_norm_7} into \eqref{eq:perturbed_cat_zero_norm_1} we obtain
    \begin{equation*}
      \frac{\sqrt{(1-\alpha^2)^{2} -\beta(\alpha)}}{1-\alpha^2} \abs{v} \le \norm{v}_0 \le \frac{\sqrt{(1-\alpha^2)^{2} +\beta(\alpha)}}{1-\alpha^2} \abs{v}
    \end{equation*}
    as announced.
  \end{proof}
\end{lemma}

\begin{lemma}\label{lemma:fejer_cond_for_perturbed_cat}
    In the setting of Lemma \ref{lemma:perturbed_cat_anosov} and Proposition \ref{prop:invariant_cones_anosov} we have $C_{\tau,0} \Theta_{T_\delta,0} <1$.
  \begin{proof}
    Since $\alpha < 1$ the conclusion of Lemma \ref{lemma:perturbed_cat_zero_norm} holds.
    Since $D_x \tau_y$ is an isometry with respect to the usual Riemannian metric on the tangent space of $\mathbb{T}^2$, by Lemma \ref{lemma:perturbed_cat_zero_norm} we have
    \begin{equation*}
      \sup_{\norm{v}_0 \le 1} \norm{D_x \tau_y v}_0 \le \sqrt{
      \frac{(1-\alpha^2)^{2} +\beta(\alpha)}
          {(1-\alpha^2)^{2} -\beta(\alpha)}}.
    \end{equation*}
    Thus
    \begin{equation}
      C_{\tau,0} \le \sqrt{
      \frac{(1-\alpha^2)^{2} +\beta(\alpha)}
          {(1-\alpha^2)^{2} -\beta(\alpha)}}.
    \end{equation}
    We turn to bounding $\Theta_{T_\delta,0}$. Since $\pi^s_x$ and $\pi^u_x$ are complementary orthogonal projections with respect to the inner product associated to $\norm{\cdot}_0$ we have
    \begin{equation*}
      \Theta_{T_\delta,0} = \sup_{x, y \in \mathbb{T}^2} \norm{\pi^s_{x+y} D_x \tau_y - (D_x \tau_y) \pi^s_{x}}_0.
    \end{equation*}
    By using Lemma \ref{lemma:perturbed_cat_zero_norm} in the same way as when bounding $C_{\tau,0}$ we find that
    \begin{equation*}
      \Theta_{T_\delta,0} \le \sqrt{\frac{(1-\alpha^2)^{2} +\beta(\alpha)}{(1-\alpha^2)^{2} -\beta(\alpha)}}
      \sup_{x, y \in \mathbb{T}^2} \wnorm{\pi^s_{x+y} D_x \tau_y - (D_x \tau_y) \pi^s_{x}}.
    \end{equation*}
    Let $v \in T_x \mathbb{T}^2$ with $\wnorm{v} = 1$. Recalling the definition of $\pi^s_x$ and then applying the triangle and Cauchy-Schwarz inequalities we have
    \begin{equation}\label{eq:fejer_cond_for_cat_1}\begin{split}
      \lvert(\pi^s_{x+y} D_x \tau_y &- (D_x \tau_y) \pi^s_{x})v\rvert \\
      = &\wnorm{\frac{\langle w_u(x+y)^\perp, v \rangle}{\langle w_u(x+y)^\perp, w_s(x+y) \rangle} w_s(x+y) - \frac{\langle w_u(x)^\perp, v \rangle}{\langle w_u(x)^\perp, w_s(x) \rangle} w_s(x)} \\
      \le &\wnorm{\frac{\langle w_u(x+y)^\perp, v \rangle}{\langle w_u(x+y)^\perp, w_s(x+y) \rangle} w_s(x+y) - \frac{\langle w_u(x)^\perp, v \rangle}{\langle w_u(x+y)^\perp, w_s(x+y) \rangle} w_s(x+y)} \\
      & + \wnorm{\frac{\langle w_u(x)^\perp, v \rangle}{\langle w_u(x+y)^\perp, w_s(x+y) \rangle} w_s(x+y) -
      \frac{\langle w_u(x)^\perp, v \rangle}{\langle w_u(x)^\perp, w_s(x) \rangle} w_s(x+y)}\\
      & + \wnorm{\frac{\langle w_u(x)^\perp, v \rangle}{\langle w_u(x)^\perp, w_s(x) \rangle} w_s(x+y) -
      \frac{\langle w_u(x)^\perp, v \rangle}{\langle w_u(x)^\perp, w_s(x) \rangle} w_s(x)}\\
      \le &\frac{\abs{\langle w_u(x+y)^\perp, v \rangle - \langle w_u(x)^\perp, v \rangle}}{\abs{\langle w_u(x+y)^\perp, w_s(x+y) \rangle}}\\
      & + \frac{\abs{\langle w_u(x+y)^\perp, w_s(x+y) \rangle - \langle w_u(x)^\perp, w_s(x) \rangle }}{\abs{\langle w_u(x+y)^\perp, w_s(x+y) \rangle \langle w_u(x)^\perp, w_s(x) \rangle}} \\
      & + \frac{\wnorm{w_s(x+y) -  w_s(x)}}{\abs{\langle w_u(x)^\perp, w_s(x) \rangle}} .
    \end{split}\end{equation}
    We will bound the various terms on the right hand side of \eqref{eq:fejer_cond_for_cat_1}.
    By Lemma \ref{lemma:new_cone_approx} we have
    \begin{equation}\label{eq:fejer_cond_for_cat_2}
       \abs{\langle w_u(x+y)^\perp, v \rangle - \langle w_u(x)^\perp, v \rangle} \le \frac{2\alpha}{\sqrt{1 + \alpha^2}}, \text{ and }  \abs{w_s(x+y) -  w_s(x)} \le \frac{2\alpha}{\sqrt{1 + \alpha^2}}.
    \end{equation}
    By definition we have $\langle w_u(x)^\perp, w_s(x)\rangle \le 0$ for every $x \in \mathbb{T}^2$.
    More precisely, by using Cauchy-Schwarz and \eqref{eq:perturbed_cat_zero_norm_4} we find that
    \begin{equation*}
      -1 \le \langle w_u(x)^\perp, w_s(x)\rangle \le \frac{\alpha^2 - 1}{1+ \alpha^2}.
    \end{equation*}
    Hence,
    \begin{equation}\label{eq:fejer_cond_for_cat_3}
      \abs{\langle w_u(x+y)^\perp, w_s(x+y) \rangle - \langle w_u(x)^\perp, w_s(x) \rangle } \le \frac{\alpha^2}{1 + \alpha^2}.
    \end{equation}
    Using \eqref{eq:fejer_cond_for_cat_3} and \eqref{eq:perturbed_cat_zero_norm_4} to bound the second term of \eqref{eq:fejer_cond_for_cat_1}, we obtain
    \begin{equation}\label{eq:fejer_cond_for_cat_4}
      \frac{\abs{\langle w_u(x+y)^\perp, w_s(x+y) \rangle - \langle w_u(x)^\perp, w_s(x) \rangle }}{\abs{\langle w_u(x+y)^\perp, w_s(x+y) \rangle \langle w_u(x)^\perp, w_s(x) \rangle}} \le \frac{\alpha^2(1 + \alpha^2)}{(1 - \alpha^2)^2}.
    \end{equation}
    Using \eqref{eq:fejer_cond_for_cat_2} and \eqref{eq:perturbed_cat_zero_norm_4} to bound the first term of \eqref{eq:fejer_cond_for_cat_1}, \eqref{eq:fejer_cond_for_cat_3} and \eqref{eq:perturbed_cat_zero_norm_4} to bound the second term,
    and the bound \eqref{eq:fejer_cond_for_cat_4} yields
    \begin{equation*}\begin{split}
      \wnorm{(\pi^s_{x+y} D_x \tau_y - (D_x \tau_y) \pi^s_{x})v} &\le \frac{2\alpha\sqrt{1+\alpha^2}}{1 - \alpha^2} + \frac{\alpha^2(1 + \alpha^2)}{(1 - \alpha^2)^2} + \frac{2\alpha\sqrt{1+\alpha^2}}{1-\alpha^2}\\
      & \le \frac{4\alpha(1-\alpha^2)\sqrt{1+\alpha^2} +\alpha^2(1+\alpha^2)}{(1 - \alpha^2)^2}.
    \end{split}\end{equation*}
    Thus,
    \begin{equation}\label{eq:fejer_cond_for_cat_5}
      \Theta_{T_\delta,0} \le \left(\sqrt{\frac{(1-\alpha^2)^{2} +\beta(\alpha)}{(1-\alpha^2)^{2} -\beta(\alpha)}}\right) \left(\frac{4\alpha(1-\alpha^2)\sqrt{1+\alpha^2} +\alpha^2(1+\alpha^2)}{(1 - \alpha^2)^2}\right).
    \end{equation}
    Combining \eqref{eq:fejer_cond_for_cat_1} and \eqref{eq:fejer_cond_for_cat_5} yields
    \begin{equation*}
      C_{\tau,0} \Theta_{T_\delta,0} \le \left(\frac{(1-\alpha^2)^{2} +\beta(\alpha)}{(1-\alpha^2)^{2} -\beta(\alpha)} \right)\left(\frac{4\alpha(1-\alpha^2)\sqrt{1+\alpha^2} +\alpha^2(1+\alpha^2)}{(1 - \alpha^2)^2}\right).
    \end{equation*}
    So if $\alpha < 0.11872$ then $C_{\tau,0} \Theta_{T_\delta,0} < 1$.
  \end{proof}
\end{lemma}

\bibliographystyle{siam}
\bibliography{/bibliography}

\begin{thebibliography}{10}

\bibitem{Aimino2015}
{\sc R.~Aimino and S.~Vaienti}, {\em A note on the large deviations for
  piecewise expanding multidimensional maps}, in Nonlinear Dynamics New
  Directions: Theoretical Aspects, Springer International Publishing, 2015,
  pp.~1--10.

\bibitem{bahsoun2016rigorous}
{\sc W.~Bahsoun, S.~Galatolo, I.~Nisoli, and X.~Niu}, {\em Rigorous
  approximation of diffusion coefficients for expanding maps}, Journal of
  Statistical Physics, 163 (2016), pp.~1486--1503.

\bibitem{bahsoun2018variance}
{\sc W.~Bahsoun, I.~Melbourne, and M.~Ruziboev}, {\em Variance continuity for
  {L}orenz flows}, arXiv preprint arXiv:1812.08998,  (2018).

\bibitem{BKL02}
{\sc M.~Blank, G.~Keller, and C.~Liverani}, {\em Ruelle--{P}erron--{F}robenius
  spectrum for {A}nosov maps}, Nonlinearity, 15 (2002), p.~1905.

\bibitem{brin2002introduction}
{\sc M.~Brin and G.~Stuck}, {\em Introduction to Dynamical Systems}, Cambridge
  University Press, 2002.

\bibitem{buckholtz2000hilbert}
{\sc D.~Buckholtz}, {\em Hilbert space idempotents and involutions},
  Proceedings of the American Mathematical Society, 128 (2000), pp.~1415--1418.

\bibitem{butterley2007smooth}
{\sc O.~Butterley and C.~Liverani}, {\em Smooth {A}nosov flows: Correlation
  spectra and stability}, Journal of Modern Dynamics, 1 (2007), pp.~301--322.

\bibitem{cfstability}
{\sc H.~Crimmins and G.~Froyland}, {\em Stability and approximation of
  statistical limit laws for multidimensional piecewise expanding maps}, arXiv
  preprint arXiv:1808.09524,  (2018).

\bibitem{DJ99}
{\sc M.~Dellnitz and O.~Junge}, {\em On the approximation of complicated
  dynamical behavior}, SIAM Journal on Numerical Analysis, 36 (1999),
  pp.~491--515.

\bibitem{Fernando2018}
{\sc K.~Fernando and P.~Hebbar}, {\em Higher order asymptotics for large
  deviations}, arXiv preprint arXiv:1811.06793,  (2018).

\bibitem{F95}
{\sc G.~Froyland}, {\em Finite approximation of {Sinai-Bowen-Ruelle} measures
  for {Anosov} systems in two dimensions}, Random and Computational Dynamics, 3
  (1995), pp.~251--264.

\bibitem{froyland2014detecting}
{\sc G.~Froyland, C.~Gonz{\'a}lez-Tokman, and A.~Quas}, {\em Detecting isolated
  spectrum of transfer and {K}oopman operators with {F}ourier analytic tools},
  Journal of Computational Dynamics, 1 (2014), pp.~249--278.

\bibitem{gouezel2015limit}
{\sc S.~Gou{\"e}zel}, {\em Limit theorems in dynamical systems using the
  spectral method}, vol.~89 of Proceedings of Symposia in Pure Mathematics,
  American Mathematical Society, 2015, pp.~161--193.

\bibitem{gouezel2006banach}
{\sc S.~Gou{\"e}zel and C.~Liverani}, {\em Banach spaces adapted to {A}nosov
  systems}, Ergodic Theory and dynamical systems, 26 (2006), pp.~189--217.

\bibitem{hennion2001limit}
{\sc H.~Hennion and L.~Herv{\'e}}, {\em Limit theorems for Markov chains and
  stochastic properties of dynamical systems by quasi-compactness}, vol.~1766,
  Springer Science \& Business Media, 2001.

\bibitem{hirsch2012differential}
{\sc M.~Hirsch}, {\em Differential Topology}, Graduate Texts in Mathematics,
  Springer New York, 2012.

\bibitem{hormander1990analysis}
{\sc L.~H{\"o}rmander}, {\em The Analysis of Linear Partial Differential
  Operators I}, Springer, 1990.

\bibitem{pollicott2017rigorous}
{\sc O.~Jenkinson, M.~Pollicott, and P.~Vytnova}, {\em Rigorous computation of
  diffusion coefficients for expanding maps}, To appear in J. Stat. Phys.,
  (2018).

\bibitem{kato1966perturbation}
{\sc T.~Kato}, {\em Perturbation theory for linear operators}, Grundlehren der
  mathematischen Wissenschaften, Springer Berlin Heidelberg, 1966.

\bibitem{katok1997introduction}
{\sc A.~Katok and B.~Hasselblatt}, {\em Introduction to the modern theory of
  dynamical systems}, vol.~54, Cambridge university press, 1997.

\bibitem{katok1989differentiability}
{\sc A.~Katok, G.~Knieper, M.~Pollicott, and H.~Weiss}, {\em Differentiability
  and analyticity of topological entropy for {A}nosov and geodesic flows},
  Inventiones mathematicae, 98 (1989), pp.~581--597.

\bibitem{katznelson2002introduction}
{\sc Y.~Katznelson}, {\em An introduction to harmonic analysis}, Cambridge
  University Press, 2002.

\bibitem{keller1999stability}
{\sc G.~Keller and C.~Liverani}, {\em Stability of the spectrum for transfer
  operators}, Annali della Scuola Normale Superiore di Pisa-Classe di Scienze,
  28 (1999), pp.~141--152.

\bibitem{Kifer_1974}
{\sc J.~I. Kifer}, {\em {On} {small} {random} {perturbations} {of} {some}
  {smooth} {dynamical} {systems}}, Mathematics of the {USSR}-Izvestiya, 8
  (1974), pp.~1083--1107.

\bibitem{lang2012real}
{\sc S.~Lang}, {\em Real and Functional Analysis}, Graduate Texts in
  Mathematics, Springer New York, 2012.

\bibitem{orey1989}
{\sc S.~Orey and S.~Pelikan}, {\em Deviations of trajectory averages and the
  defect in {P}esin's formula for {A}nosov diffeomorphisms}, Transactions of
  the American Mathematical Society, 315 (1989), pp.~741--753.

\bibitem{ratner1973central}
{\sc M.~Ratner}, {\em The central limit theorem for geodesic flows on
  $n$-dimensional manifolds of negative curvature}, Israel Journal of
  Mathematics, 16 (1973), pp.~181--197.

\bibitem{roe1988elliptic}
{\sc J.~Roe}, {\em Elliptic operators, topology and asymptotic methods},
  Chapman and Hall/CRC, 1988.

\bibitem{rousseau1983theoreme}
{\sc J.~Rousseau-Egele}, {\em Un th{\'e}oreme de la limite locale pour une
  classe de transformations dilatantes et monotones par morceaux}, The Annals
  of Probability,  (1983), pp.~772--788.

\bibitem{ruelle1997differentiation}
{\sc D.~Ruelle}, {\em Differentiation of {SRB} states}, Communications in
  Mathematical Physics, 187 (1997), pp.~227--241.

\bibitem{ruelle2008differentiation}
\leavevmode\vrule height 2pt depth -1.6pt width 23pt, {\em Differentiation of
  {SRB} states for hyperbolic flows}, Ergodic Theory and Dynamical Systems, 28
  (2008), pp.~613--631.

\bibitem{dynamicmode2019}
{\sc J.~Slipantschuk, O.~F. Bandtlow, and W.~Just}, {\em Dynamic mode
  decomposition for analytic maps}, arXiv preprint arXiv:1905.09266,  (2019).

\bibitem{wormell2017spectral}
{\sc C.~Wormell}, {\em Spectral {G}alerkin methods for transfer operators in
  uniformly expanding dynamics}, Numerische Mathematik, 142 (2019),
  pp.~421--463.

\bibitem{young1986stochastic}
{\sc L.-S. Young}, {\em Stochastic stability of hyperbolic attractors}, Ergodic
  Theory and Dynamical Systems, 6 (1986), pp.~311--319.

\bibitem{young2002srb}
\leavevmode\vrule height 2pt depth -1.6pt width 23pt, {\em What are {SRB}
  measures, and which dynamical systems have them?}, Journal of Statistical
  Physics, 108 (2002), pp.~733--754.

\end{thebibliography}

\end{document}